\documentclass[11pt]{amsart}
\usepackage{amsmath,amsfonts,amssymb}
\usepackage{mathtools}
\usepackage[margin=2.1cm]{geometry}
\usepackage[english]{babel} 
\usepackage{mathrsfs}
\usepackage{textcomp}
\usepackage{fancyhdr}
\usepackage{eso-pic}
\usepackage{bbold}
\usepackage{color}
\usepackage{xcolor}
\usepackage{textgreek}
\usepackage{comment}
\usepackage{stmaryrd}
\usepackage[utf8]{inputenc}   
\usepackage[T1]{fontenc}
\usepackage[shortlabels]{enumitem}
\newtheorem{theorem}{Theorem}[section]
\newtheorem{lemma}[theorem]{Lemma}

\theoremstyle{definition}
\newtheorem{definition}[theorem]{Definition}

\newtheorem{rema}[theorem]{Remark}
\newtheorem{coro}[theorem]{Corollary}
\newtheorem{prop}[theorem]{Proposition}
\numberwithin{equation}{section}

\newcommand{\abs}[1]{\lvert#1\rvert}

\let\dem=\proof
\let\enddem=\endproof

\DeclarePairedDelimiter\norm\lVert\rVert

\DeclareMathOperator\dvg{div}
\DeclareMathOperator\rot{curl}
\DeclareMathOperator{\adj}{Adj}

\newcommand\R{\mathbb{R}}
\newcommand\N{\mathbb{N}}
\newcommand\D{\mathbb{D}}

\newcommand\dpt{\partial_t}

\newcommand\cC{\mathscr{C}}

\newcommand\loc{\text{loc}}
\definecolor{mycolor}{HTML}{D35400}

\definecolor{hide}{HTML}{A0AC81}
\definecolor{myred}{HTML}{8A1538}

\usepackage[colorlinks=true,citecolor=blue,urlcolor=cyan]{hyperref}
\usepackage{cleveref}
\crefformat{rema}{#2 Remark~#1#3}{}
\crefformat{lemma}{#2Lemma~#1#3}{}
\crefformat{theo}{#2Theorem~#1#3}{}
\crefformat{appendix}{#2Appendix~#1#3}{}
\crefformat{section}{#2Section~#1#3}{}
\crefformat{def}{#2Definition~#1#3}{}
\crefformat{prop}{#2Proposition~#1#3}{}
\crefformat{coro}{#2Corollary~#1#3}{}
\usepackage{todonotes}

\begin{document}

\title[]{ Global-in-time Well-Posedness   of the Compressible {N}avier-{S}tokes Equations with Striated Density }

\author[X. Liao]{Xian Liao}
\address{(X. Liao) Karlsruhe Institut of Technology, Karlsruhe, Germany}
\curraddr{}
\email{xian.liao@kit.edu}
\thanks{}

\author[S. M. Zodji]{Sagbo Marcel Zodji}
\address{(S. M. Zodji) Institut de Mathématiques de Jussieu Paris Rive Gauche, Université Paris Cité, Paris, France}
\email{marcel.zodji@imj-prg.fr}
\thanks{}


\date{\today}

\dedicatory{}


\begin{abstract}
We first show local-in-time well-posedness of the compressible Navier-Stokes equations, assuming  striated regularity while no other smoothness or smallness conditions on the initial  density.
With these local-in-time solutions served as blocks, 
for  \textit{less} regular initial data where the vacuum is permitted,  the global-in-time well-posedness  follows 
 from the energy estimates and the propagated striated regularity of the density function, 
if the    bulk viscosity coefficient is large enough in the two dimensional case.
The global-in-time well-posedness holds also true in the three dimensional case,
  provided with large bulk viscosity coefficient together with small initial energy.
This solves the density-patch problem in the exterior domain for the compressible model with $W^{2,p}$-Interfaces. 
Finally,    the singular incompressible limit toward the inhomogenous incompressible model when the bulk viscosity coefficient tends to infinity is obtained. 
\end{abstract}

\maketitle
\section{Introduction}
In this paper, we establish the existence and uniqueness of global-in-time weak solutions of  compressible  viscous flows,   and at the same time we investigate the dynamics of  density-interfaces   in dimension $d \in \{2,3\}$. More precisely, we consider the following compressible Navier-Stokes equations describing the motion of  compressible  viscous fluids:
\begin{gather}\label{ep0.1}
    \begin{cases}
        \dpt \rho +\dvg (\rho u)=0,\\
        \dpt (\rho u)+\dvg (\rho u\otimes u)+ \nabla P(\rho)= \mu \Delta u+(\mu+\lambda)\nabla \dvg u.
    \end{cases}
\end{gather}
Here $\mu>0$ represents the dynamic viscosity, and $\lambda>0$ stands for the kinetic  viscosity. 
In the present paper, $\mu$ is some fixed positive constant   while the constant $\lambda$ may become very large.
For notational simplicity, we introduce the so-called bulk viscosity coefficient 
    \begin{equation*}\label{nu}\nu=2\mu+\lambda,
    \end{equation*}
    which tends to infinity when $\lambda\rightarrow\infty$.
      We always assume that our fluids are (strictly) viscous:
\begin{equation*}\label{assumption:nu}
\nu\geqslant \underline{\nu}>0,
\end{equation*} 
where $\underline{\nu}$ is a fixed positive constant.

 In the above, $t\geq 0$, $x\in \R^d$, $d=2,3$ denote the time and space variables respectively.
 The notations $\rho= \rho(t,x)\geqslant 0$ and $u=u(t,x)\in \R^d$ represent, respectively, the density and velocity 
of the compressible fluid, which serve as the unknowns in the problem. Meanwhile, $P=P(\rho)$ is a given smooth (in this paper we assume $P\in \cC^2(\R,\R)$) and increasing function  of the density (that is, $P'(\rho)>0$). The system \eqref{ep0.1} is supplemented with initial data 
\begin{gather}\label{ep4.30}
    (\rho, \rho u)_{|t=0}=(\rho_0, \rho_0 u_0), 
\end{gather}
which satisfy 
\begin{gather}\label{ep3.16}
    \rho_0\geqslant  0, \quad \rho_0\in L^\infty(\R^d;[0,\infty)),\quad  \rho_0-\widetilde \rho \in L^2(\R^d;\R),
    \quad u_0\in H^1(\R^d;\R^d),
\end{gather}
where $\widetilde\rho>0$ is some given positive equilibrium state of the density. 

\subsection{Striated regularity}\label[section]{Subs:TR}
We   assume further striated regularity with respect to a given nondegenerate family of vector fields for  the initial density $\rho_0$ in this paragraph.

We first   introduce some notations, based on \cite{danchin2020well}.
For some $p\in (d , \infty)$, $\mathbb L^{\infty,p}(\R^d; \R^d)$ denotes the vector space of bounded   vector fields  with gradients in $L^p(\R^d; \R^{d\times d})$.
From now on we denote  the Lebesgue  spaces $L^p(\R^d;\R^n)$   resp. Sobolev spaces $H^s(\R^d;\R^n)$ with $p\in [1,\infty]$, $s\in\R$ and $n\in \mathbb{N}^*$, simply by $L^p(\R^d)$ or $L^p$, resp. $H^s(\R^d)$ or $H^s$ with an abuse of notations. 
We have defined
\begin{gather*}
    \mathbb L^{\infty,p}(\R^d) =\left\{Y\in L^\infty(\R^d)\;\Big| \; 
    \norm{Y}_{\mathbb L^{\infty,p}(\R^d)}:= \norm{Y}_{L^\infty(\R^d)} +\norm{\nabla Y}_{L^p(\R^d)}<\infty\right\}.
\end{gather*} 
We define the norm $\|\cdot\|_{\mathbb L^{\infty,p}}$ for a family of vector fields $\mathcal{Y}= (Y_1, Y_2,\dots, Y_m)\subset \mathbb L^{\infty,p}(\R^d)$, $m\in \N$ as
\begin{gather*}
    \norm{\mathcal Y}_{\mathbb L^{\infty,p}(\R^d)}:= \sup_{1\leqslant \upsilon\leqslant m}\norm{Y_{\upsilon}}_{\mathbb L^{\infty,p}(\R^d)}.
\end{gather*}

\begin{definition}[Nondegeneracy]\label{def}
    Let $\mathcal Y=(Y_1, Y_2,\dots , Y_m)\subset \mathbb L^{\infty,p}(\R^d)$ be a family of  $m$ vector fields with $m\geqslant d-1$ and $p\in (d,\infty)$. We say that $\mathcal Y$ is \emph{nondegenerate}
    if it satisfies the following property 
    \begin{gather*}
    I(\mathcal Y):=\inf_{x\in \R^d}\sup_{\Upsilon\in \Upsilon^m_{d-1}}\bigg|\bigwedge^{d-1} Y_\Upsilon(x)\bigg|^{\tfrac{1}{d-1}} >0.
    \end{gather*}
    Above $\Upsilon \in \Upsilon^m_{d-1}$ means that $\Upsilon=(\upsilon_1, \upsilon_2,\cdots, \upsilon_{d-1})$ with each $\upsilon_i\in \{1,\cdots,m\}$ and $\upsilon_i<\upsilon_j$ for $i<j$, 
    $Y_{\Upsilon}:=(Y_{\upsilon_1}, Y_{\upsilon_2}, \dots, Y_{\upsilon_{d-1}})$, while the symbol $\displaystyle\bigwedge ^{d-1} Y_{\Upsilon}$ stands for the unique element of $\R^d$ such that 
    \begin{gather*}
         \left(\bigwedge ^{d-1} Y_{\Upsilon}\right) \cdot Z=\det (Y_{\upsilon_1}, Y_{\upsilon_2}, \cdots, Y_{\upsilon_{d-1}}, Z),\quad \forall\, Z \in \R^d.
    \end{gather*}
\end{definition}

\begin{definition}[Striated regularity with respect to a nondegenerate family of vector fields]\label[def]{def:SR}
Let $Y\in \mathbb{L}^{\infty, p}(\R^d)$, $p\in (d,\infty)$ be  a       (single) vector field, and $\mathcal Y=(Y_{1}, Y_2, \dots , Y_m)\subset \mathbb L^{\infty,p}(\R^d)$ be a \emph{nondegenerate} family of vector fields with $m\geqslant  d-1$.
\begin{enumerate}[(a)]
    \item  
A function $g\in L^\infty(\R^d)$ is said to be of class $L^p(\R^d)$
    along $Y$, if
    $$g\in \mathbb L^{ p}_{Y}(\R^d):=\{g\in L^\infty(\R^d)\,| \dvg (g Y)\in L^p(\R^d)\}.$$ 
    We define the 
derivative of the function $g$ along $Y$ as follows 
\begin{gather*}
    \partial_Y g
    := \dvg (g Y)- g \dvg Y,
\end{gather*}
and hence we can equivalently define $$\mathbb L^{ p}_{Y}(\R^d)=\{g\in L^\infty(\R^d)\, | \,\partial_Y g\in L^p(\R^d)\}.$$
\item   A function $g\in L^\infty(\R^d)$ is said to be of  class $L^p(\R^d)$ along the family $\mathcal{Y}$, if
    \[
    g\in \mathbb{L}_{\mathcal Y}^p(\R^d):=\bigcap_{1\leqslant \upsilon\leqslant m} \mathbb L^{p}_{Y_{\upsilon}}(\R^d),
    \]
   and we equip the space $\mathbb L^{p}_{\mathcal Y}(\R^d)$    with the following norm
    \begin{gather*}\label{striatedregularity}
    \norm{g}_{\mathbb L^p_{\mathcal Y}(\R^d)}:= \dfrac{1}{I(\mathcal Y)}\sup_{1\leqslant \upsilon\leqslant m}\left[\norm{g}_{L^\infty(\R^d)}\norm{Y_\upsilon}_{\mathbb L^{\infty,p}(\R^d)}+\norm{\dvg(g Y_\upsilon)}_{L^p(\R^d)}\right],
    \end{gather*}
    which is equivalent to the norm with $\dvg(gY_\upsilon)$ above replaced by $\partial_{Y_\upsilon}g$.
\end{enumerate} 
\end{definition}
 
We now continue with the assumption of the initial density $\rho_0$ given in \eqref{ep4.30}-\eqref{ep3.16} associated with the compressible Navier-Stokes equations \eqref{ep0.1}. We assume further that 
\begin{gather}\label{ep4.29}
  \rho_0\in \mathbb L^p_{\mathcal X_0}(\R^d),
\end{gather}
where $\mathcal{X}_0= (X_{0,1}, \dots X_{0,m})\subset \mathbb L^{\infty,p}(\R^d) $ is a given \emph{nondegenerate} family of vector fields for some $m\geqslant  d-1$ and $p\in (d,\infty)$.
\begin{rema}[Initial density of density-patch type]\label[rema]{rem:patch}
    It's interesting to notice that   the initial density of the form 
      \begin{equation}\label{rho0,patch}
      \rho_0= \alpha \mathbb{1}_{D_0}+ \widetilde \rho \mathbb{1}_{D_0^c}, \quad \alpha\geqslant  0,
      \end{equation}
      satisfies the assumptions for $\rho_0$ in  \eqref{ep3.16}-\eqref{ep4.29}, if $D_0$ is a $W^{2,p}(\R^d)$ (with $p>d$) bounded, simply connected domain in $\R^d$.
      Indeed, \eqref{ep4.29} holds for
       a nondegenerate (divergence-free) family of vector fields $\mathcal{X}_0= (X_{0,1}, \dots X_{0,m})\subset \mathbb L^{\infty,p}(\R^d)$ which is tangent to $\partial D_0$
      \footnote{ Indeed, for $d=2$   the existence of such a nondegenerate family of tangential vector fields is obvious since we can take $X_{0,1}=\begin{pmatrix}
           \partial_{x_2}f\\ -\partial_{x_1}f
       \end{pmatrix}=:\nabla^\perp f$ to be the tangent vector field close to   $\partial D_0$ with $f\in W^{2,p}(\R^2)$  and $f|_{\partial D_0}=0$ and $\nabla f|_{\partial D_0}\neq 0$, while $X_{0,2}=\nabla^\perp(\chi x_1)$ to be a  non-zero  vector field   with $\chi$ a smooth cutoff function away from the  boundary, see e.g. \cite[(1.10)]{paicu2020striated} (with $m=3$).
       The existence result for $d=3$ with $m=5$ follows from the similar idea, see e.g. \cite[Proposition 3.2]{GSR}, where $X_{0,1}=\begin{pmatrix}
           0\\-\partial_{x_3}f\\ \partial_{x_2}f
       \end{pmatrix}$, $X_{0,2}=\begin{pmatrix}
           \partial_{x_3}f\\0\\-\partial_{x_1}f
       \end{pmatrix}$, $X_{0,3}=\begin{pmatrix}
           -\partial_{x_2}f\\ \partial_{x_1}f\\0
       \end{pmatrix}$ are generated by the function $f\in W^{2,p}(\R^3)$ with $f|_{\partial D_0}=0$ and $\nabla f|_{\partial D_0}\neq 0$, while $X_{0,4}=\begin{pmatrix}
           \partial_{x_3}(\chi x_3)\\0\\-\partial_{x_1}(\chi x_3)
       \end{pmatrix}$, $X_{0,5}=\begin{pmatrix}
           -\partial_{x_2}(\chi x_1)\\ \partial_{x_1}(\chi x_1)\\0
       \end{pmatrix}$ forms a nondegenerate family away from $\partial D_0$ with $\chi$ a smooth cutoff function away from the boundary $\partial D_0$.}, and this means that the initial density given by \eqref{rho0,patch} persists tangential regularity with respect to the boundary $\partial D_0$.
\end{rema}
\subsection{Statement of the main results}\vspace{0.3cm}
The purpose of this paper is threefold.
\begin{enumerate}[1.]
    \item We establish the local-in-time well-posedness of the system \eqref{ep0.1} for   \textit{positive} density function with striated regularity, under the compatibility condition. We thus remove the smallness condition required on the density in  Danchin, Fanelli, Paicu's paper 
    \cite{danchin2020well}.
    
    \item These local-in-time solutions become  global-in-time unique  solutions of the Cauchy problem  \eqref{ep0.1}-\eqref{ep4.30}-\eqref{ep3.16}-\eqref{ep4.29}, if
    \begin{itemize}
        \item $d=2$, and the bulk viscosity coefficient is large enough $\nu\geq \nu_0$ with $\nu_0$ depending on the norms of the initial data given in \eqref{ep4.30}-\eqref{ep3.16}. 
        This result is inspired by   the work by Danchin and Mucha \cite{danchin2023compressible}.
        \item $d=3$,   the initial energy is small and the bulk viscosity coefficient $\nu\geq\nu_0$ is large enough.
        Here although $\|\rho_0-\widetilde\rho\|_{L^2(\R^d)}$ is assumed to be small, $\rho_0$ may have large variation in $L^\infty(\R^d)$. This result supplements the local-in-time well-posedness work \cite{danchin2020well} with global-in-time well-posedness result and the work by Shibata and Zhang \cite{shibata2023global} with less regular  initial data. 
    \end{itemize}
    \item Additionally, by letting the bulk   viscosity tend to infinity $\nu\rightarrow\infty$, we establish a singular limit toward the incompressible  inhomogeneous model \textit{on the whole space}, in the spirit of Danchin and Mucha's work \cite{danchin2023compressible} where the considered domain has finite measure. 
\end{enumerate}

 \subsubsection{Local-in-time well-posedness and continuation criterion}
We begin by providing the statement of the local-in-time result, which technically further assumes the strict positivity of the initial density function and the compatibility condition on the initial data. 
\begin{theorem}[Local-in-time well-posedness and continuation criterion]\label[theo]{intro:local}
    We consider the Cauchy problem of the compressible Navier-Stokes equations \eqref{ep0.1} supplemented with the initial data \eqref{ep4.30} verifying \eqref{ep3.16} and \eqref{ep4.29}. We further assume the strict positivity of the initial density and the compatibility condition as follows
\begin{gather}\label{ep4.44}
    0<\underline\rho\leqslant \rho_0(x) \quad\text{and}\quad \mu \Delta u_0+ (\mu+\lambda)\nabla\dvg u_0-\nabla P(\rho_0)\in L^2(\R^d).
\end{gather}

     Then, there exists a time $T>0$ and a unique solution $(\rho,u)$ to the Cauchy problem \eqref{ep0.1}-\eqref{ep4.30}, satisfying the following properties:
    \begin{enumerate}[1)]
        \item Energy bounds: $u\in \cC([0,T], H^1(\R^d))$, $\dot u\in \cC([0,T], L^2(\R^d))$,   $\sqrt{\sigma}\nabla \dot u, \sigma \ddot u\in L^\infty((0,T), L^2(\R^d))$,\\
        and
         $\nabla \dot u,\;\sqrt{\sigma}\ddot u,\; \sigma \nabla\ddot u \in L^2((0,T)\times \R^d)$.
         
        Here and in what follows we use the notations
         \begin{align}\label{sigma}
        &\sigma=\sigma(t):=\min\{1,t\},
        \quad
         \dot v:=\bigl(\dpt  +  u\cdot \nabla \bigr) v,
        \quad\ddot v:= \bigl(\dpt  +  u\cdot \nabla \bigr) \dot v.
    \end{align} 

     \item Striated regularity: For all $0<t<T$, we have $\rho(t)\in \mathbb L_{\mathcal X(t)}^p(\R^d)$, where $\mathcal X(t)=(X_{\upsilon}(t))_{1\leqslant \upsilon\leqslant m}\subset \mathbb L^{\infty,p}(\R^d)$ is the nondegenerate family of vector fields transported by the fluid flow, in the sense that
           each vector field $X_\upsilon(t)$, $1\leqslant \upsilon\leqslant m$  solves uniquely the following Cauchy problem 
        \begin{equation}\label{ep4.33}
           \left\{ \begin{array}{l}
                \dpt X_\upsilon +u\cdot \nabla X_\upsilon= \partial_{X_\upsilon} u,\\
                {X_{\upsilon}}_{|t=0}= X_{0,\upsilon}.
            \end{array}\right.
        \end{equation} 
       Here the directional derivative was defined in  \cref{def:SR}: $\partial_{X_\upsilon} u^j=\dvg(u^j X_{\upsilon})-u^j\dvg X_{\upsilon}$, $1\leqslant j\leqslant d$.
       
        The velocity field is Lipschitz continuous when integrated in time   and persists also the striated regularity for positive times:
        \begin{itemize}
            \item For $d=2$ or for  $d=3$ and $3<p\leqslant 6$,   $ \nabla u\in L^2((0,T), \mathbb L^p_{\mathcal X}(\R^d))$;
        \item For $d=3$\; and \;$6<p<\infty$, 
          $\sigma^{\tfrac{3}{4}-\tfrac{1}{r}-\tfrac{3}{2p}}\nabla u\in  L^r((0,T),  \mathbb L^p_{\mathcal X}(\R^3))$
          and \;$\sigma^{\tfrac{3}{4}-\tfrac{1}{r}}\nabla \dot u\in L^r((0,T), L^3(\R^3))$, for any $2\leqslant r\leqslant \infty$.
        \end{itemize}  

        \item Continuation criterion: If  $(\rho,u)$ is the solution   defined up to a maximal time $T^*>0$ and
    \begin{align}
        \limsup_{t\to T^*}\left\{ \norm{\mathcal X(t)}_{\mathbb L^{\infty,p}(\R^d)}+\dfrac{1}{I(\mathcal X(t))}
        +\left\|\dfrac{1}{\rho(t)}\right\|_{L^\infty(\R^d)}
        +\left\|\rho(t)\right\|_{L^\infty(\R^d)}\right.\notag\\
        \left.+ \norm{\partial_{\mathcal X(t)}\rho(t)}_{L^p(\R^d)}+ \norm{\nabla u(t)}_{L^2(\R^d)}+\norm{\dot u(t)}_{L^2(\R^d)}\right\}<\infty,\label{intro:thblowup}
    \end{align}
    then $T^*=\infty$.
    \end{enumerate} 
\end{theorem}
The solution of \cref{intro:local} is constructed in the spirit of a recent contribution of the second author \cite{zodji2023well}, which deals with the more involved case of density-dependent viscosity coefficient. Thus, we only present a sketch of the proof of 
\cref{intro:local} in \cref{local}.
\begin{rema} 
\begin{enumerate} [(a)]
    \item This result supplements the contribution by Danchin, Fanelli and Paicu  \cite{danchin2020well} by removing the smallness condition on the  density deviation. Unlike the maximum regularity argument used in \cite{danchin2020well}, which requires a critical regularity for one part of the initial velocity, our method relies on the change into Lagrangian coordinates along with energy estimation methods.
    \item The compatibility condition in \eqref{ep4.44},
$\mu \Delta u_0+ (\mu+\lambda)\nabla \dvg u_0- \nabla P(\rho_0)\in L^2(\R^d) $
expresses the continuity of the stress tensor, and
does not require (explicitly) smoothness of the density.
The parabolic effect of the momentum equations ensures that this condition holds true at positive times even for less regular initial data, see \cite{hoff1995global}.
\item  The velocity field persists further regularity  property which are stated in \cref{coro1} below, thanks to the decomposition of the velocity gradient \eqref{ep4.19}.
\end{enumerate}
\end{rema}

\subsubsection{Definitions of energy functionals}\label[section]{subss:def}
 The global-in-time well-posedness result will follow from the above local-in-time well-posedness, continuation criterion and a series of energy estimates.
 We define in this subsection the relevant  energy functionals.

 Recall the positive density equilibrium $\widetilde\rho$ in \eqref{ep3.16}.
 We define  first
 \begin{align}
     &\hbox{The pressure equilibrium }\widetilde P:=P(\widetilde\rho),\label{tildeP}\\
     &\hbox{The $\rho$-dependent functions } 
H_l(\rho)= \rho \int_{\widetilde\rho}^\rho s^{-2}\abs{P(s)-\widetilde P}^{l-1}(P(s)-\widetilde P)ds, \quad l\in [1,\infty),\label{ep1.8}\\
     &\hbox{The pressure deviation }
     G(t,x) :=\bigl(P(\rho)\bigr)(t,x)-\widetilde P, \label{g}\\
     &\hbox{The effective flux }F(t,x)=\nu (\dvg u)(t,x)- G(t,x).\label{F}
 \end{align}  

 We define the associated energy function of the compressible Navier-Stokes equations \eqref{ep0.1}
 \begin{equation}\label{intro:E}
     E(t)= \int_{\R^d} \left[\rho \dfrac{\abs{u}^2}{2}+ H_1(\rho)\right](t,x)dx+\int_0^t \left[\mu \norm{\nabla u(t')}_{L^2(\R^d)}^2+ (\mu+\lambda)\norm{\dvg u(t')}_{L^2(\R^d)}^2\right]dt',
 \end{equation}
which consists of 
    \begin{itemize}
        \item The kinetic energy $\frac12\|\sqrt{\rho}u(t)\|_{L^2(\R^d)}^2$;
        \item The potential energy $\displaystyle\int_{\R^d}\bigl(H_1(\rho)\bigr)(t,x) dx$ with $H_1(\rho)$ defined in \eqref{ep1.8}: $\displaystyle H_1(\rho)=\rho\int_{\widetilde\rho}^{\rho} \frac{P(s)-\widetilde P}{s^2}ds$;
        \item The energy dissipation $\mu\|\nabla u\|_{L^2((0,t)\times\R^d)}^2+(\mu+\lambda)\|\dvg u\|_{L^2((0,t)\times\R^d)}^2$.
    \end{itemize}
The energy $E(t)$ is conserved for regular enough solutions of \eqref{ep0.1}.
    
Recall the notations in \eqref{sigma}, and we introduce   energy functionals of higher order
  \begin{equation}\label{ep4.42}
        \begin{array}{l}
            \displaystyle
            \mathcal{A}_1(t)=\frac{\mu}2\norm{\nabla u(t)}_{L^2(\R^d)}^2+ \frac{\mu+\lambda}2\norm{\dvg u(t)}_{L^2(\R^d)}^2+\int_0^t \norm{\sqrt{\rho}\dot u(t')}_{L^2(\R^d)}^2dt',\\
            \displaystyle
            \mathcal{A}_2(t)= \sigma(t)\norm{\sqrt{\rho}\dot u(t)}_{L^2(\R^d)}^2+\int_0^t \sigma(t')\left[\mu\norm{\nabla \dot u(t')}_{L^2(\R^d)}^2+\dfrac{\mu+\lambda}{\nu^2} \norm{\dot F(t')}_{L^2(\R^d)}^2\right]dt'.
        \end{array}
    \end{equation}
   The hierarchy of energy functionals $E(t), \mathcal{A}_1(t), \mathcal{A}_2(t)$ encode $L^2(\R^d)$-norm, $\dot H^1(\R^d)$-norm for $u(t)$ and (time-weighted) $L^2(\R^d)$-norm for the material derivative $\dot u(t)$, respectively.
    Although trivially $|\dvg u|\leqslant d|\nabla u|$, we will make efforts to get the (large)   viscosity coefficient $\lambda$ before $\dvg u$ in the definition of $E, \mathcal{A}_1$, such that intuitively $\dvg u\rightarrow 0$ as $\lambda\rightarrow\infty$  if $E,\mathcal{A}_1$ is bounded uniformly in time. 
 The review of their history  can be found in  \cref{Subs:Review} below.

Recall the initial data \eqref{ep4.30} and we denote $G_0(x)=G(0,x)=\bigl(P(\rho_0)\bigr)(x)-\widetilde P$.
    For notational simplicity we denote
    \begin{align}
        &\hbox{The first  initial energy }
        E_0:= E(0)=\int_{\R^d} \left[\rho_0 \dfrac{\abs{u_0}^2}{2}+ H_1(\rho_0)\right](x)\,dx,\label{E0}\\
        &\hbox{The total  initial energy }
        E_0^\nu:=E_0+ \mu\norm{\nabla u_0}_{L^2(\R^d)}^2+ \nu \norm{\dvg u_0}_{L^2(\R^d)}^2+ \dfrac{1}{\nu}\norm{G_0}_{L^2(\R^d)}^2,\label{E0-lambda}\\
        &\hbox{The upper bound of the initial density }
         \rho_0^\ast:=\sup_{x\in \R^d}\rho_0(x).\label{rho0ast}
    \end{align}
    We observe that for initial data given in \eqref{ep4.30}-\eqref{ep3.16}, 
    \begin{align*}
        E_0\leq C(\rho_0^\ast) \|(\rho_0-\widetilde\rho, u_0)\|_{L^2(\R^d)}<+\infty, \quad E_0^\nu\leq C(\mu,  \underline{\nu}, \rho_0^\ast)(E_0+  \|\nabla u_0\|_{L^2(\R^d)})+\nu\|\dvg u_0\|_{L^2(\R^d)}^2<+\infty.
    \end{align*}
    We aim to bound $ \mathcal{A}_1(t), \mathcal{A}_2(t)$ globally in time in terms of $E_0^\nu, \rho_0^\ast$, if $\nu\geq\nu_0$ is large enough (and if the initial energy is small enough for $d=3$).
  The following  quantity $\mathcal{A}_3(t)$ captures the striated regularity of the density function along the family of vector fields $\mathcal X(t)=(X_{\upsilon}(t))_{1\leqslant \upsilon\leqslant m}$ transported by the flow as in \eqref{ep4.33}
      \begin{gather}\label{ep4.41}
            \mathcal{A}_3(t)= \norm{\mathcal{X}(t)}_{\mathbb L^{\infty,p}(\R^d)}+\sup_{1\leqslant \upsilon\leqslant m}
            \norm{(\partial_{ X_\upsilon}\rho)(t)}_{L^p(\R^d)}. 
        \end{gather} 
        It is straightforward to see that  $\mathcal{A}_3(t)$ grows exponentially in $\|\nabla u\|_{L^1_tL^\infty}$.
        We aim to show that the striated regularity encoded in $\|\log \mathcal{A}_3\|_{L^1((0,t))}$ together with the energy functionals $\mathcal{A}_1(t), \mathcal{A}_2(t)$ controls   $\|\nabla u\|_{L^1_tL^\infty}$.
        Gronwall's inequality hence implies the exponential-in-time control   of  $\|\nabla u\|_{L^1_tL^\infty}$.
\subsubsection{Global-in-time well-posedness}
We now state our global-in-time result for less regular initial data on which the assumption \eqref{ep4.44} is not assumed. 
\begin{theorem}\label[theo]{th1} 
Assume   the  Cauchy problem \eqref{ep0.1}-\eqref{ep4.30}-\eqref{ep3.16}-\eqref{ep4.29} and the following conditions
 \begin{align}
     &  \hbox{either }   d=2\hbox{ and }  \nu\geqslant \nu_0,\label{nu0}\\
&\hbox{or }       d=3,\, p\in (3,6),\, 
 E_0^\nu  E_0\leqslant c  \hbox{ and } \nu \geqslant \nu_0, 
\label{c0}
        \end{align}
        where $c$ is a fixed constant depending only on $\mu, \underline{\nu}$ while  $\nu_0$ is a constant depending additionally on the initial norms: $E_0, \|\nabla u_0\|_{L^2(\R^d)},
        \rho_0^\ast$.
Then the Cauchy problem   has a unique global-in-time  solution $(\rho,u)$   verifying 
        \begin{enumerate}[1)]
            \item  Energy bounds: For all $t\geqslant 0$, we have 
            \begin{gather}\label{ep4.32}
            \begin{cases}\vspace{0.25cm}\displaystyle
                E(t)+ \mathcal{A}_1(t)+ \mathcal{A}_2(t) 
                \leqslant   C_0^\nu,\\
                \displaystyle
                \norm{\rho(t)-\widetilde \rho}_{L^\infty(\R^d)}^2\leqslant \norm{\rho_0-\widetilde \rho}_{L^\infty(\R^d)}^2 +  C_0^\nu,
            \end{cases} 
        \end{gather}
        where the constant $C_0^\nu$ depends   on $\mu, \underline{\nu}$, $\rho_0^\ast$, and (superlinearly) on $E_0^\nu$. 
        \item  Striated regularity: For all $t\geqslant 0$, $\rho(t)\in \mathbb L^p_{\mathcal X(t)}(\R^d)$, where $\mathcal X(t)=(X_{\upsilon}(t))_{1\leqslant \upsilon\leqslant m}\subset \mathbb L^{\infty,p}(\R^d)$ is a nondegenerate family of vector fields defined to solve the   Cauchy problem \eqref{ep4.33}.
        
        Moreover $\nabla u\in L^1_\loc ([0,\infty), L^\infty(\R^d))$ with the following estimates:
        \begin{gather}\label{A3}
          \hspace{-.5cm}  \begin{cases}\displaystyle
                \mathcal{A}_3(t)\leqslant   \mathcal{A}_3(0) \exp\Bigl( 
    C_0 \int_0^t\Bigl[1+\sqrt{t}+\norm{\nabla u(t')}_{L^\infty(\R^d)}dt'\Bigr]\Bigr),
    \displaystyle
          \int_0^t \norm{\nabla u(t')}_{L^\infty(\R^d)}dt'\leqslant C_0\bigl(1+\frac{\mathcal{A}_3(0)}{I(\mathcal{X}_0)}\bigr)\exp(C_0t), \\ 
          \displaystyle
        \int_0^t \norm{\partial_{\mathcal X(t')}\nabla u(t')}_{L^p(\R^d)}dt'\leqslant C_0 (1+t+t\mathcal{A}_3(t)) \mathcal{A}_3(t),
            \end{cases}
        \end{gather}     
        where   $C_0$ depends  on $\mu, \underline{\nu}, m,  d, p, \rho_0^\ast, E_0^\nu$.
        \end{enumerate}
\end{theorem}
\begin{rema} [Bounds for $\dvg u$]
We have assumed some uniform bounds (with respect to $\nu$) for $\dvg u_0$ implicitly: The conditions in \eqref{nu0} and \eqref{c0} imply that
\[
\nu \norm{\dvg u_0}_{L^2(\R^d)}^2\leqslant 
\begin{cases}
    E_0^\nu<\infty,  \quad &\text{if}\quad d=2,\\
    E_0^\nu\min\left\{1,\; \dfrac{c}{E_0}\right\} <\infty,  \quad &\text{if}\quad d=3.
\end{cases}
\]
 This boundedness is propagated over time:
\[
\nu \norm{\dvg u(t)}_{L^2(\R^d)}^2\leqslant  C_0^\nu.
\]
\end{rema}
\cref{th1} and \cref{rem:patch} imply immediately
\begin{coro}[Density patch problem in the exterior domain]
The Cauchy problem \eqref{ep0.1}-\eqref{ep4.30} with initial density of density-patch type \eqref{rho0,patch} and $u_0\in H^1(\R^d)$, under the assumption \eqref{nu0} or \eqref{c0}, has a unique global-in-time solution $(\rho,u)$ with $\rho(t)$ persisting tangential regularity with respect to the boundary $\partial D_t$ which is transported by the flow of $u$ and keeps its $W^{2,p}(\R^d)$-regularity. 
\end{coro}

\begin{rema} 
     If $\alpha>0$, deriving a uniformly positive lower bound for the density is straightforward (see \cite[Section 2.1 (Step 6)]{danchin2023exponential}). 
     This results in an exponential-in-time decay  of the jump in the density $\rho(t)$ across $\partial D_t$, as observed in \cite{hoff2002dynamics, hoff2008lagrangean}.
\end{rema}
Intuitively, thanks to the uniform bound in \eqref{ep4.32}: $\mathcal{A}_1(t)\leqslant C_0^\nu$, letting $\nu\rightarrow\infty$  yields a couple $(\varrho, v)$ that satisfies the incompressible inhomogeneous model:
\begin{gather}\label{ep4.35}
    \begin{cases}
        \dpt \varrho +\dvg (\varrho v)=0,\\
        \dpt (\varrho v)+ \dvg (\varrho v\otimes v)+\nabla \Pi -\mu \Delta v=0,\\
        \dvg v=0.
    \end{cases}
\end{gather} 
\begin{coro}[Incompressible limit]\label[coro]{theo1}
    Let   $(\rho_0,u_0)$ be the initial data given in \eqref{ep4.30} verifying  \eqref{ep3.16}, \eqref{ep4.29}  and  $\dvg u_0=0$. 
    Let $(\rho^{(\nu)}, u^{(\nu)})$ be the corresponding unique solution constructed in \cref{th1}, under the assumption \eqref{nu0} or \eqref{c0}.
    
    Then $(\rho^{(\nu)}, u^{(\nu)})_\nu$ converges weakly-* to $(\varrho, v)$ in $L^\infty((0,\infty)\times \R^d)\times L^\infty((0,\infty), H^1(\R^d))$ as $\nu$ goes to infinity, and $(\varrho, v)$ solves (uniquely)
    the inhomogeneous, incompressible model \eqref{ep4.35} with initial data $(\rho_0, u_0)$ in the distribution sense. Moreover, we have 
    \begin{gather}
        \begin{cases}
            \dvg u^{(\nu)} =\mathcal O(\nu^{-1/2}) \quad \text{ in } \quad L^2\cap L^\infty((0,\infty), L^2(\R^d)),\\
            \dpt (\rho^{(\nu)} u^{(\nu)})+\dvg (\rho^{(\nu)} u^{(\nu)} \otimes u^{(\nu)})-\nabla F^{(\nu)} -\mu \Delta u^{(\nu)}= \mathcal{O}(\nu^{-1/2})\;\; \text{in} \;\; L^\infty((0,\infty), \dot H^{-1}(\R^d)),
        \end{cases}
    \end{gather}
    where $F^{(\nu)}=\nu\dvg u^{(\nu)}-G^{(\nu)}$ with
    $G^{(\nu)}=P(\rho^{(\nu)})-\widetilde P$.
\end{coro}
The proofs of \cref{th1} and \cref{theo1} are presented  in \cref{approximate},   based on the a priori estimates   in \cref{sketch} and their proofs in \cref{proofs}.
\begin{rema}
    This result in \cref{theo1} is a partial continuation of the work by Danchin and Mucha \cite{danchin2016compressible,Danchin_2019,danchin2023compressible}, and Danchin and Wang \cite{danchin2023exponential}, and stands, as far as we know, as the first one dealing with \emph{discontinuous} initial data \emph{in the whole space}. 
    We notice that, except for the work \cite{danchin2016compressible} dealing with the whole space case and initial data in the critical Besov space, the other  studies rely heavily on the assumption that the domain has finite measure. The extension to the whole space, especially for $d=2$, is not obvious, and it requires some refined computations, e.g. the compensated result by Coifman, Lions, Meyer, Semmes in  \cite{coifman1993compensated}.
\end{rema} 

\subsection{Review of known results}\label[section]{Subs:Review}

Classical solutions for the Navier-Stokes equations \eqref{ep0.1} with regular initial data are known to exist, since the work by Nash \cite{nash1962probleme}, Itaya \cite{ itaya1970existence,itaya1971cauchy}, Solonnikov \cite{solonnikov1980solvability}, Tani \cite{tani1976existence}, just to cite a few examples. These solutions are defined up to a positive time which depends on the (norms of) initial data.  The first result addressing the global-in-time well-posedness of classical solutions is provided by Matsumura and Nishida \cite{matsumura1980initial} for small initial data in $L^1(\mathbb{R}^3)\cap H^3(\mathbb{R}^3)$. Nowadays, global-in-time classical solutions are known to exist under smallness assumption on the initial data in critical Besov space \cite{charve2010global,chen2010global,haspot2011existence}.
 
\noindent\textbf{Weak solutions and estimates for $E(t)$ and $G$.}
Similar to the solutions constructed by Leray \cite{leray1934mouvement}  for the \emph{incompressible} Navier-Stokes equations, there are well-established results that investigate the existence of global-in-time \emph{weak} solutions for the compressible Navier-Stokes equations \eqref{ep0.1}, with finite initial energy. The first result was obtained by P.-L. Lions \cite{lions1996mathematical}, followed by Feireisl, Novotn\'y, Petzeltov\'a \cite{feireisl2001existence}, for pressure laws of the form $P(\rho) = a \rho^\gamma$, $a>0$, with some limitations on $\gamma$.
These weak solutions verify the following classical energy inequality:
\begin{gather}\label{EE:intro}
    E(t)   \leqslant E(0)=E_0,
\end{gather}
where the functional $E$ has been  given in  \eqref{ep4.42}.

The introduction of the (generalized) specific energy $H_l(\rho)$,   $l\in [1,\infty)$, in \eqref{ep1.8} helps (technically) to estimate the pressure deviation $G$.
As observed in e.g. \cite{bresch2021extension}, the so-defined $H_l(\rho)$ is non-negative: $H_l(\rho)\geqslant  0$, since the pressure $P(\rho)$ is an increasing function of the density.

For the classical case $l=1$, $H_1(\rho)$ appears in the definition of $E(t)$, which is integrable in space uniformly in time due to \eqref{EE:intro}: $\displaystyle\int_{\R^d} \bigl(H_1(\rho)\bigr)(t,x)\,dx
    \leqslant   E(t)=    E_0$.
Consequently, under the a priori assumption
$$\rho(t,x)\leqslant \rho^\ast,$$
 we have the estimates for $G$ uniformly in time by the energy $E_0$ below
\begin{gather}\label{ep4.46}
    \sup_{[0,t]}\norm{G}_{L^q(\R^d)}^q\leqslant C^\ast \sup_{t'\in [0,t]}\int_{\R^d} H_1(\rho(t',x))\,dx
    \leqslant   C^\ast  E_0, \quad \text{with}\quad q\in [2,\infty),
\end{gather}
where the constant $C^\ast$ depends only on $\rho^\ast$ and $q$.
In the above, the first inequality follows from the definition of $H_1(\rho)$   in \eqref{ep1.8}.

General $ H_l(\rho) $, $l\geqslant  1$, as a function of $\rho$, satisfies the following ordinary differential equation 
\[
\rho H'_l(\rho)- H_l (\rho)= \abs{P(\rho)-\widetilde P}^{l-1} (P(\rho)-\widetilde P),
\]
and hence, by virtue of the mass equation $\eqref{ep0.1}_1$, the function 
$\bigl( H_l(\rho)\bigr)(t,x)$ satisfies the following time evolutionary equation 
\begin{equation*}\label{ep1.32}
    \dpt H_l(\rho)+\dvg (H_l(\rho)u)+ \abs{P(\rho)-\widetilde P}^{l-1} (P(\rho)-\widetilde P)\dvg u=0,
\end{equation*}
which is in the same spirit of the renormalized continuity equation appearing in e.g. \cite{lions1996mathematical}.
By view of the definitions \eqref{g} and \eqref{F}, it is equivalent to
\begin{equation*}\label{ep3.19}
    \dpt H_l(\rho)+\dvg (H_l(\rho)u)+ \frac1\nu |G|^{l+1}  
    =-\frac{1}{\nu}\abs{G}^{l-1}G F.
\end{equation*}
Consequently,   integrating the above in  space yields, after H\"older's inequality, the following  
  \begin{gather}\label{ep1.40}
 \dfrac{d}{dt}\| H_l(\rho)(t,x)\|_{L^1(\R^d)}+\dfrac{1}{\nu}\norm{G}_{L^{l+1}(\R^d)}^{l+1}\leqslant \dfrac{1}{\nu}\norm{G}_{L^{l+1}(\R^d)}^{l} \norm{F}_{L^{l+1}(\R^d)}^{l+1},
 \end{gather}
and hence by Young's inequality and integration in time,  $G$ can be controlled by $F$ in the following way 
\begin{gather}\label{ep1.9}
      \dfrac{1}{\nu} \norm{G}_{L^{l+1}((0,t)\times \R^d)}^{l+1} \leqslant  2\norm{H_l(\rho_0)(x)}_{L^1(\R^d)}+ \dfrac{C}{\nu} \norm{F}_{L^{l+1}((0,t)\times \R^d)}^{l+1},\quad \forall t\in (0,\infty).
\end{gather}

\vspace{0.2cm}
\noindent\textbf{Density patch problem.}
In the last three decades, there has been   growing interests in exploring the properties of weak solutions to models arising from fluid mechanics that enable tracking down discontinuities of some quantities such as density or vorticity. 
We refer to the   density patch problem for incompressible models stated in \cite{lions1996mathematicalv1}:
\emph{Consider the incompressible model \eqref{ep4.35} in two dimension with initial density as the characteristic function of some regular domain $D_0\in \R^2$:   $\rho_0=\mathbb{1}_{D_0}$. The   density-patch problem asks whether or not, at positive times, the density is still some characteristic function $\mathbb{1}_{D(t)}$ with the domain $D(t)\subset \R^2$ preserving the initial regularity of $D_0$}.
This problem is almost solved for incompressible models, even for density-dependent viscosity (under some smallness assumption) or higher Sobolev regularity of $D_0$, see  \cite{danchin2012lagrangian, danchin2013incompressible, danchin2019incompressible, danchin2017persistence, denisova2001evolution, Denisova2008,gancedo2023global, liao2019global, liao2016global, liao2019globallow}.

However, for a similar problem in the context of compressible fluids, there are not so many results. On one hand, the global classical solutions constructed by Matsumura and Nishida, or in critical Besov space, are too strong in a way that they do not allow for discontinuous solutions. On the other hand, the weak solutions constructed by P.-L. Lions or Feireisl, Novotn\'y, Petzeltov\'a only require that the initial energy is finite, allowing for discontinuous density. However, the velocity is relatively weak, with $\nabla u\in L^2((0,\infty)\times \mathbb{R}^d)$, and this is insufficient to track down discontinuities in the density.
A natural idea is to construct weak solutions in a class that allows for tracking down the discontinuity of the interface.
The first result addressing this issue is, as far as we know, \cite{hoff2002dynamics} by Hoff, where the author considered an initial density with H\"older regularity on both sides of a suitable curve, allowing for jumps across this curve.
The initial curve is transported by the flow of the velocity into a curve that maintains its initial regularity. The density also remains H\"older continuous on both sides of the transported curve, and moreover, its jump through the latter decays exponentially over time. This result pertains only in the case of linear pressure law and 
small bulk  viscosity.
Recently, these restrictions were removed in \cite{zodji2023discontinuous}, even in the more 
challenging case of density-dependent viscosity. \cref{th1} is thus added to this list,  in the constant viscosity setting, with domains having Sobolev regularity, and the density can be large in $L^\infty(\R^d)$, unlike the cited results.

\subsubsection{Hoff's strategy}
We review briefly some key concepts in Hoff's works \cite{hoff1995global,hoff1995strongpoly,hoff2002dynamics,hoff2008lagrangean}.

\noindent\textbf{Energy functionals.}
In \cite{hoff1995global}, Hoff introduced the following energy functionals which can be compared with our definitions in \eqref{ep4.42}
\begin{gather*}
\begin{cases}\displaystyle
    \mathcal{A}_1^H(t)=\sup_{[0,t]}\sigma\norm{\nabla u}_{L^2(\R^d)}^2+\int_0^t \sigma(t') \norm{\sqrt{\rho}\dot u(t')}^2_{L^2(\R^d)}d t',\\
    \displaystyle\mathcal{A}_2^H(t)= \sup_{[0,t]}\sigma^d\norm{\sqrt{\rho}\dot u}_{L^2(\R^d)}^2+\int_0^t \sigma^d(t') \norm{\nabla \dot u(t')}_{L^2(\R^d)}^2 dt',\\
    \displaystyle\mathcal{B}(t)= \sup_{[0,t]} \norm{\rho-\widetilde\rho}_{L^\infty(\R^d)}^2,
\end{cases}
\end{gather*}
where the time weight $\sigma$ and the material derivative $\dot u$ are defined as in \eqref{sigma}.
He provides bounds for these functionals by requiring that the initial velocity is small in $L^2(\R^d)$ but can be large in $L^{2^d}(\R^d)$. Additionally, he requires that  the initial density is bounded away from zero and bounded from above, along with some technical assumptions. 

\noindent\textbf{Effective flux and vorticity.}
Hoff's computations, mainly while propagating the lower and upper bounds of the density, rely strongly on   the \emph{effective viscous flux} $F$ given in \eqref{F}: $F=\nu\dvg u-G.$
It plays a crucial role by connecting the momentum equations and the mass equation, which was discovered by Hoff and  Smoller in \cite{hoff1985solutions}.
It presents its power also in  the study of the propagation of oscillations in \cite{serre1991variations}, and in the constructions of weak solutions in \cite{feireisl2001existence,hoff1995strongpoly,lions1996mathematical}. 

In fact, recall the momentum equations $\eqref{ep0.1}_2$, which can be written by virtue of the mass conservation law $\eqref{ep0.1}_1$ as
$$\rho \dot u-\mu\Delta u-(\mu+\lambda)\nabla\dvg u+\nabla \bigl(P(\rho)-\widetilde P\bigr) =0.$$ 
We apply the divergence operator to it to obtain the Poisson equation for $F$ as follows
\begin{gather}\label{ep4.2}
\Delta F= \dvg (\rho \dot u).
\end{gather}
Similarly we can apply the curl operator to the momentum equations, to derive the Poisson equation for   the vorticity, $\rot u,$ as follows
\begin{gather}\label{ep4.3}
\mu\Delta \rot u= \rot (\rho \dot u).
\end{gather}

Consequently, the regularity of the material derivative of the velocity $\dot u$, as provided by functionals $\mathcal{A}_1^H$ and $\mathcal{A}_2^H$, allows the effective flux $F$  and the vorticity $\rot u$  to be regular at positive time, even for  rough density. This means that there is some cancellation between the divergence of the velocity and the pressure at positive times. In particular, the fact that $F\in L^{8/3}((1,\infty), L^\infty(\R^d))$ allows him to propagate the $L^\infty(\R^d)$ estimate for the density. 

Thanks to this observation, under smallness condition on the initial data, Hoff proved existence of global weak solutions for the system \eqref{ep0.1} with a linear pressure law in a first paper \cite{hoff1995global}. He later considered pressure laws of the form $P(\rho)= a\rho ^\gamma$, with $\gamma>1$ in a second paper \cite{hoff1995strongpoly}, in which, again, the effective flux played a crucial role in proving compactness for the density.

\noindent\textbf{Velocity gradient expression involving Riesz operators.}
In order to study the dynamics of discontinuous surfaces,  Hoff in \cite{hoff2002dynamics} used the following decomposition of the velocity gradient:
\begin{gather}\label{ep4.40}
    \mu \nabla u= -(-\Delta)^{-1} \nabla (\rho \dot u)+ \dfrac{\mu+\lambda}{\nu}\mathcal{R}\mathcal{R} F+ \dfrac{\mu}{\nu}\mathcal{R}\mathcal{R} G
    =:\mu\nabla \widetilde u+ \mu\nabla u_G,
\end{gather}
which is nothing but a rewriting of the above momentum equations, where $\mathcal{R}_j= (\frac{1}{i}\partial_j) (-\Delta)^{-\tfrac{1}{2}}$, $1\leqslant j\leqslant d$,  are the Riesz operators.

By assuming more  regularity on the velocity $u_0\in H^\beta(\R^2)$, he reduces the singularity of time weighs in the definitions of functionals  $\mathcal{A}_1^H$ and $\mathcal{A}_2^H$.
Namely, in dimension two,  the time weights $\sigma$ and $\sigma^2$ are replaced, respectively, by $\sigma^{1-\beta}$ and $\sigma^{2-\beta}$.
Thus, $\nabla \widetilde u$ and the effective flux $F$ belong to  
$L^1_\loc ([0, \infty), \cC^\alpha(\R^2))$ for all $0<\alpha<\beta$. 
With the help of the regularity of $F$, Hoff propagated the \textit{piecewise} H\"older regularity of the density, resulting   \textit{piecewise} H\"older continuity of $\nabla u_G$ on both sides of a time-dependent curve. This time-dependent curve is the transport of an initial suitable curve with some geometric assumptions, and only provided with bounded velocity gradient can the structure of the density and of the curve be propagated.

However, since Riesz operators \emph{fail} to be continuous on $L^\infty(\R^d)$, additional regularity must be assumed on the density to obtain $\mathcal{R}\mathcal{R} G\in L^\infty(\R^d)$.  
In \cite{hoff2008lagrangean}, Hoff and Santos observed that in the configuration of the previous works (see \cite{hoff1995global,hoff1995strongpoly}), the rough part of the velocity gradient $\nabla u_G$ belongs to $L^\infty((0,\infty), BMO(\R^d))$. 
In this case, the  initial interface $\gamma_0\in \cC^\alpha$, $\alpha>0$,  is transported to an interface $\gamma_t \in \cC^{\alpha_t}$ at time $t>0$, with $\alpha_t$ decaying exponentially to zero in time. 
  
Hence to propagate interface regularity (more than continuity) requires a Lipschitz velocity. For the incompressible model with constant viscosity, this regularity is directly obtained from energy computations and interpolations. In contrast, for the compressible case with discontinuous density, the problem is more delicate, and the issue is to find an appropriate functional space framework such that even-order Riesz operators maps into $L^\infty(\R^d)$. 

\subsubsection{The strategy by use of tangential regularity}
Apart from the tools used in \cite{hoff2002dynamics, zodji2023discontinuous} to handle the rough part of the velocity gradient, there exists another framework that allows for the same. 
It is referred  as  tangential/striated regularity space, which goes back to Chemin's study  (see  e.g. \cite{chemin1991mouvement,chemin1993persistance})   of the vortex patch problem for the ideal incompressible model.
See also \cite{bertozzi1993} for another interesting geometric proof for the persistence of regularity in the vortex patch problem. 
Chemin's idea has been  further developed to higher dimensional cases in \cite{Raphael_Danchin_Persistance, GSR}, to the inhomogeneous case in \cite{fanelli2012},
as well as   the density patch problem of the inhomogeneous incompressible Navier-Stokes model in e.g.
\cite{  liao2019global, liao2016global, liao2019globallow, paicu2020striated}.
However, there are very few results in this direction for the compressible case. 
To the best of our knowledge, the only one work in the literature is \cite{danchin2020well} by Danchin, Fanelli, and Paicu. They establish the local-in-time well-posedness of the compressible equations \eqref{ep0.1} with a striated initial density, and we  now delve into a brief discussion of their methods. From the momentum equations $\eqref{ep0.1}_2$,  they express the velocity as:
\begin{gather}\label{ep4.43}
    u= w-\nabla (I_d-\Delta)^{-1} G,
\end{gather}
where $w$ solves a  parabolic equation with source term in some suitable space $L^r((0,T),L^p(\R^d))$. They employ maximal regularity tools to establish Lipschitz bounds for $w$. Meanwhile, Lipschitz bound for the second term of the velocity's expression \eqref{ep4.43}, associated with the pressure, is obtained through tangential regularity estimates.
The maximal regularity argument requires smallness assumption on the density in $L^\infty(\R^d)$, and the global-in-time result is still missing.
Toward this, we establish local-in-time well-posedness of the system \eqref{ep0.1} without imposing any smallness condition on the initial data (see \cref{intro:local}), and  global-in-time well-posedness (see \cref{th1}) without any smallness assumption of the initial density fluctuation in $L^\infty(\R^d)$, and the vacuum is allowed.
This is accomplished through a coupling mechanism that involves the effective flux. By achieving this objective, we propagate the Sobolev regularity of interfaces over time. 

\vspace{0.2cm}
\noindent\textbf{Incompressible limit.}
We aim also to establish an incompressible limit in the spirit of the work of Danchin and Mucha \cite{danchin2023compressible}. Let us look briefly at this question. 
The work by Matsumura and Nishida \cite{matsumura1980initial} paved the way for attempts to relax the assumptions on the initial data. Despite reducing the regularity assumption to critical Besov space or even Lebesgue space, the condition of smallness is frequently encountered in the literature. In their work \cite{danchin2016compressible}, Danchin and Mucha introduced a new framework that enables them to bypass the smallness condition on the initial data, namely, replacing the  smallness in the initial data by large enough    bulk viscosity coefficient.
In particular, as the bulk viscosity $\nu\rightarrow\infty$, the solution converges to a limit that satisfies the incompressible model. 
This has been down for initial data in critical Besov space. 
For less regular initial data, they work on the torus in \cite{Danchin_2019,danchin2023compressible}, where they rely technically on
  the finite-measure of the domain, particularly on the logarithmic interpolation inequality, which proves to be crucial in handling vacuum states in \cite{danchin2023compressible}.
We also refer to the work by Danchin and Wang \cite{danchin2023exponential} where exponential decay rate of the solutions of the compressible model on torus has been investigated. 
However, the exponential decay  does not generally hold in the whole space. For instance, the work by Hu and Wu \cite{hu2020optimal} provides lower bound for the norms of solutions in certain cases. 
We obtain  similar results to those in \cite{danchin2023compressible}  in the presence of vacuum  on the whole space (see \cref{theo1}), where we apply technically, in dimension two.   The extension to the whole space, particularly in two dimension, 
the compensated result by Coifman, Lions, Meyer, Semmes \cite{coifman1993compensated}.


\subsection*{Outline of the paper}
The rest of the paper is structured as follows. In  next section  \cref{sec1}, we give the proofs of \cref{th1} and \cref{theo1}, provided with the validity of \cref{intro:local}, whose 
proof is postponed in  \cref{local}.
A useful density-weighted interpolation inequality is established in \cref{sec:inter}.

\section{Proof of the main results}\label[section]{sec1}
This section is devoted to the proofs of  \cref{th1} and \cref{theo1}, which goes from a priori estimates for  solutions of the Navier-Stokes equations \eqref{ep0.1} to the proof of the compactness of approximate solutions. It is divided into three parts. In the first one, \cref{sketch}, we summarise all key ideas with brief explanations and give the a priori estimates in a series of lemmas. Technical details and the proofs of these lemmas are presented in the second part \cref{proofs}.
As we will see in the final part of the proof in \cref{approximate}, the existence of a local-in-time solution (without any smallness condition in the density) is by no means obvious, and it is the purpose of \cref{local}. 
The regularity of the (local-in-time) solution is sufficient in order to use $u$ and $\dot u$ as test functions in the subsequent computations to get energy estimates.
\subsection{Proof ideas and statements of lemmas}\label[section]{sketch}
In this section, we give the main ideas of the proof of  \cref{th1}.
We state the energy estimates for  the solutions of the compressible Navier-Stokes equations \eqref{ep0.1}  with initial data \eqref{ep4.30} satisfying \eqref{ep3.16}.
The tangential regularity  \eqref{ep4.29} is assumed when we show the boundedness of the Lipschitz-norm of the velocity vector field as a second step.
Recall the   definitions of  the energy functionals
$$E(t), \mathcal{A}_1(t), \mathcal{A}_2(t), \mathcal{A}_3(t)$$ 
together with the notations $\widetilde P, H_l(\rho), G, F$ and $E_0,E_0^\nu,\rho_0^\ast $, given in \cref{subss:def}.

 In the literature (see e.g. \cite{lions1996mathematical} where finally only the \textit{energy inequality} \eqref{EE:intro} was   established for  weak solutions)   the following \textit{a priori energy equality} for $E(t)$ was  shown for strong solutions
\begin{gather}\label{ep1.1}
    E(t)= E(0)=:E_0.
\end{gather} 
More precisely, it follows from  taking the scalar product of the momentum equation $\eqref{ep0.1}_2$ with the velocity $u$ and then integrating in time and space.
This energy balance \eqref{ep1.1} is going to be used freely in the proof, and we aim to show the estimates for $\mathcal{A}_1, \mathcal{A}_2, \mathcal{A}_3$.

In the following  we state step by step
\begin{itemize}
    \item In \cref{subsubsA1A2rho}: Energy estimates for $\mathcal{A}_1, \mathcal{A}_2$ together with the boundedness of the density deviation $\|\rho-\widetilde\rho\|_{L^\infty_{t,x}}$.

    Under the assumption that 
the density 
  is a priori upper bounded 
  \begin{equation}\label{assumption:rho,bdd}
      0\leqslant \rho(t,x)\leqslant \rho^\ast,
  \end{equation}
  for some $\rho^\ast>0$,   
we  show first (local-in-time) a priori \emph{energy estimates} for $\mathcal{A}_1, \mathcal{A}_2$ (see \cref{Aprioriestimates1} and \cref{Aprioriestimates2}) and then a  \textit{boundedness of the density} in terms of $\mathcal{A}_1, \mathcal{A}_2$ (see \cref{Aprioriestimates3}) for   solutions of the Cauchy problem   \eqref{ep0.1}-\eqref{ep4.30}-\eqref{ep3.16}. 
Under the assumption \eqref{nu0} or \eqref{c0}, that is, in the case of either \textit{large bulk viscosity coefficient for $d=2$} or with \textit{small initial energy and large bulk visocosity coefficient for $d=3$},  a bootstrap argument implies the \textit{global-in-time} a priori energy estimates for $\mathcal{A}_1, \mathcal{A}_2$, and density bound estimate (see \cref{lemma4}). 

\item In \cref{subsubsA3nablau}: The striated regularity estimate for $\mathcal{A}_3$ together with the boundedness of the velocity gradient $\|\nabla u\|_{L^1_tL^\infty_x}$. 

With the  estimates in   \cref{subsubsA1A2rho} at hand, we turn to the striated regularity for the density function $\mathcal{A}_3(t)$ for solutions of the Cauchy problem \eqref{ep0.1}-\eqref{ep4.30}-\eqref{ep3.16}-\eqref{ep4.29}, which finally implies the Lipschitz-continuity of the velocity field
(see \cref{Aprioriestimates4}), thanks to the $L^\infty$-estimates for the double Riesz-operators provided with extra striated regularity (see  \cref{Linftyestimate}).
\end{itemize}

\subsubsection{A priori estimates for $\mathcal{A}_1, \mathcal{A}_2$ and $\|\rho-\widetilde\rho\|_{L^\infty_{t,x}}$}\label[section]{subsubsA1A2rho}

In order to derive higher-order energy estimates for the velocity $u$ and its material derivative 
$$\dot u:=(\partial_t+u\cdot\nabla) u,$$ 
we use first $\dot u$ as a test function in the weak formulation of the momentum equation 
$\eqref{ep0.1}_2$ to establish bounds for $\mathcal{A}_1$.  
The functional $\mathcal{A}_2$ emerges when, first rewriting the momentum equation $\eqref{ep0.1}_2$ with the effective flux $F$, and then  applying the operator 
$\dpt \;\cdot \;+ \dvg (\; \cdot\; u)$ 
to the resulting equation before testing it with $\dot u$. 

In two dimension, the following estimates are valid for these functionals $\mathcal{A}_1, \mathcal{A}_2$.
\begin{lemma}[Preliminary energy estimates for $d=2$]\label[lemma]{Aprioriestimates1}
    Assume that $d=2$ 
    and \eqref{assumption:rho,bdd}.
    Then the following a priori bounds hold true for the functionals $\mathcal{A}_1$ and $\mathcal{A}_2$: 
    \begin{align}
       &\mathcal{A}_1(t)\leqslant C^\ast \left(   E_0^\nu + \frac{1}{\nu^{3/2} } \mathcal{A}_1(t) (E_0+\mathcal{A}_1(t))   \right)\exp\left(C^* E_0\right),
       \label{Aprioriestimates1:A1}\\
       & \mathcal{A}_2(t)\leqslant C^\ast\left( (E_0+\frac{1}{\nu^4}E_0^2)
    +   (1+ E_0+ \mathcal{A}_1(t)) \mathcal{A}_1(t) \right),
   \label{Aprioriestimates1:A2}
    \end{align}
where the constant $C^\ast$ depends on $\mu,\underline{\nu}$ and (increasingly) on   $\rho^*$.
 \end{lemma}

The proof of \cref{Aprioriestimates1} is established through refined computations, and the compensated result by Coifman \emph{et al.} \cite{coifman1993compensated} turns out to be crucial for achieving a uniform bound with respect to $\lambda$. We refer to  \cref{proofs1} below for the detailed proof. 

For $d=3$, the following estimates hold true for functionals $\mathcal{A}_1$ and $\mathcal{A}_2$. 
\begin{lemma}[Preliminary energy estimates for $d=3$]\label[lemma]{Aprioriestimates2}
 Assume that $d=3$ 
 and \eqref{assumption:rho,bdd}. Then, the  following estimates hold true for the functionals $\mathcal{A}_1$ and $\mathcal{A}_2$:
    \begin{align}
       &\mathcal{A}_1(t)\leqslant C^\ast \bigl(E_0^\nu+ \dfrac{1}{\nu^{\frac23}}E_0^{\tfrac{1}{3}}\bigr) + CE_0\mathcal{A}_1(t)^2,
      \label{Aprioriestimates2:A1}\\
       &\mathcal{A}_2(t)\leqslant C^\ast\Bigl( \frac{1}{\nu}E_0^{\frac13}+E_0+E_0^2+(1 +\mathcal{A}_1(t)^2)\mathcal{A}_1(t)\Bigr).
 \label{Aprioriestimates2:A2}
\end{align}
Here $C$ depends on $\mu, \underline{\nu},$ and $C^\ast$ depends on $\mu, \underline{\nu}$ and (increasingly) on $\rho^\ast$.
\end{lemma}

The proof is given in \cref{proofs2}. Let us point out that the computations in \cite{Danchin_2019,danchin2023compressible,danchin2023exponential} depend heavily on the fact that the domain has finite measure. \cref{Aprioriestimates1}
and \cref{Aprioriestimates2} are the first to provide bounds for the solution $(\rho, u)$ uniformly with respect to $\lambda$ (large) in the whole space, with only bounded density. 

Based on the above estimates, it turns out that the functionals $\mathcal{A}_1$ and $\mathcal{A}_2$ are under control as long as the density is upper-bounded. Therefore, the next step is devoted to estimating the upper bound of the density, whose proof   is given in \cref{Proofs3}. 

\begin{lemma}[Density upper bound in terms of energies]\label[lemma]{Aprioriestimates3}
 Assume 
 \eqref{assumption:rho,bdd}. 
 Then the following bounds hold true for the density:
 \begin{align}\label{Aprioriestimates:rho}
      \norm{\rho-\widetilde\rho}_{L^\infty([0,t]\times\R^d)} 
     &\leqslant \norm{\rho_0-\widetilde\rho}_{L^\infty(\R^d)}
    \nonumber \\
     &+ \dfrac{C^\ast}{\nu^{\frac13}}\times\left\{\begin{array}{ll}
        \big(1+E_{0}^{\tfrac{1}{18}}\big)
   \big( E_0^{\tfrac{1}{6}}+\nu^{\tfrac{1}{6}}\mathcal{A}_1(t)^{\tfrac{1}{6}}\big)\big(\mathcal{A}_1(t)^{\tfrac{1}{3}}+\mathcal{A}_2(t)^{\tfrac{1}{3}}\big),   &  d=2, \\
     \big(\mathcal{A}_1(t)^{\tfrac{1}{2}}+\mathcal{A}_2(t)^{\tfrac{1}{2}}\big),    & d=3.
     \end{array}\right.
 \end{align}
\end{lemma}

Finally, we notice that for $d=2$ there is a small factor $\frac{1}{\nu}$ (or its positive powers) before $\mathcal{A}_1(t), \mathcal{A}_2(t)$ in the estimates \eqref{Aprioriestimates1:A1} and $\eqref{Aprioriestimates:rho}$, while $\mathcal{A}_2(t)$ can be bounded by $\mathcal{A}_1(t)$ and $\rho^\ast$ by \eqref{Aprioriestimates1:A2}. We can close the estimates in \cref{Aprioriestimates1}, \cref{Aprioriestimates2} and \cref{Aprioriestimates3} by a bootstrap argument as
in e.g. \cite{bresch2021extension, danchin2023compressible}, which is not repeated here.
\begin{lemma}[Global-in-time estimates under the assumption \eqref{nu0}  or   \eqref{c0}]\label[lemma]{lemma4}
There exist $c$ depending only on  $\mu, \underline{\nu}$ and $\nu_0\geqslant \underline{\nu}$ depending on $\mu, \underline{\nu}, E_0, \|\nabla u_0\|_{L^2(\R^d)}, 
\rho_0^\ast$   such that 
    \begin{enumerate}
        \item If $d=2$ and   $\nu\geqslant \nu_0$, then
        \[
        \mathcal{A}_1(t)+\mathcal{A}_2(t)\leqslant   C_0^\nu \quad \text{and}\quad \norm{\rho-\widetilde\rho}_{L^\infty([0,t]\times\R^2)}\leqslant \norm{\rho_0-\widetilde\rho}_{L^\infty(\R^2)}+   (C_0^\nu)^{\tfrac{1}{2}}.
        \]
        \item  If  $d=3$, 
       $ E_0^\nu E_0  \leqslant c$   and $\nu\geqslant \nu_0, $
       then    
        \[
        \mathcal{A}_1(t)+\mathcal{A}_2(t)\leqslant   C_0^\nu\quad \text{and}\quad \norm{\rho-\widetilde\rho}_{L^\infty([0,t]\times\R^3)}\leqslant  \norm{\rho_0-\widetilde\rho}_{L^\infty(\R^3)}+  (C_0^\nu)^{\tfrac{1}{2}}.
        \]
        Above, $C_0^\nu$ depends    on $\mu$,  $\underline{\nu}$, $\rho_0^\ast$ and (superlinearly) on $E_0^\nu$.
    \end{enumerate}
\end{lemma}

\subsubsection{A priori estimates for $\mathcal{A}_3$ and $\|\nabla u\|_{L^1_t L^\infty_x}$}\label{subsubsA3nablau}
Now we have \cref{lemma4}, which gives the (uniform) bounds of functionals $\mathcal{A}_1$, $\mathcal{A}_2$ and $\rho$.
We use the notation $C_0$ below 
to denote some time-independent constant depending on the initial data as follows:
\begin{equation}\label{C0}
C_0=C_0 (\mu, \underline{\nu}, m, d,p, \rho_0^\ast, E_0^\nu),
\end{equation}
where $m,p$ appears in the initial condition \eqref{ep4.29}.
$C_0$ may vary from lines to lines and bounds in particular $\mathcal{A}_1, \mathcal{A}_2, \rho$.
The next step is dedicated to translating these bounds   into the tangential regularity estimates for the density, together with the Lipschitz norm of the velocity. 
As the tangential regularity   $\mathcal{A}_3$ can be transported by Lipschitz continuous flow, we sketch the idea to show Lipschitz continuity of $u$   as follows.

 We first recall the following decomposition of the velocity gradient 
 \begin{align}
\nabla u&= \nabla \widetilde u+\nabla u_G \label{ep4.36} \\
&:=\Bigl( -\dfrac{1}{\nu} \mathcal{R} \mathcal{R} (-\Delta)^{-1}\dvg (\rho \dot u)-\dfrac{1}{\mu}\mathcal{R}\mathcal{R}(-\Delta)^{-1}\cdot\rot (\rho \dot u)\Bigr)
+\Bigl( \dfrac{1}{\nu} \mathcal{R} \mathcal{R}G\Bigr).\nonumber
\end{align} 
where $\mathcal{R}_j= (\tfrac{1}{i}\partial_j)(-\Delta)^{-\tfrac{1}{2}}$, with $1\leqslant j\leqslant d$, is the Riesz transform and $G=P(\rho)-\widetilde P$.
Indeed,  we notice that the following expression 
\begin{equation}\label{Deltau}
\Delta u^j = \partial_j \dvg u+ \partial_k \rot_{jk} u,\quad \text{with}\quad \rot_{jk} u=\partial_k u^j-\partial_j u^k,\quad j,k=1,\cdots,d,
\end{equation}
and from \eqref{F}, \eqref{ep4.2} and \eqref{ep4.3} we have 
\begin{gather}\label{divu,F,curlu}
\dvg u=\dfrac{1}{\nu} \left(F+ G\right),\; \quad  F= -(-\Delta)^{-1}\dvg (\rho \dot u)\quad \text{and}\quad \mu \rot u= -(-\Delta)^{-1}  \rot (\rho \dot u).
\end{gather}
Hence the velocity gradient can be expressed as in \eqref{ep4.36}:
\begin{align*}
\nabla u&= -\nabla (-\Delta)^{-1} \Delta u\nonumber\\
        &=-\nabla (-\Delta)^{-1} \nabla \dvg u- \nabla (-\Delta)^{-1}\dvg (\rot u)\nonumber\\
        &=-\dfrac{1}{\nu}\nabla (-\Delta)^{-1} \nabla \left( F+ G\right)- \dfrac{1}{\mu}\nabla (-\Delta)^{-1}\dvg (\mu\rot u)\nonumber\\
&=\left( -\dfrac{1}{\nu} \mathcal{R} \mathcal{R} (-\Delta)^{-1}\dvg (\rho \dot u)-\dfrac{1}{\mu}\mathcal{R}\mathcal{R}(-\Delta)^{-1}\cdot\rot (\rho \dot u)\right)+\Bigl(\dfrac{1}{\nu} \mathcal{R} \mathcal{R} G\Bigr).
\end{align*}

Thanks to  the regularity of $\dot u$ provided by functionals $\mathcal{A}_1$ and $\mathcal{A}_2$, $\nabla \widetilde u\in L^1((0, t), L^\infty(\R^d))$.
Motivated by the pioneering work of Chemin \cite{chemin1991mouvement,chemin1993persistance} and Danchin, Fanelli, and Paicu \cite{danchin2020well}  which show $\mathcal{RR}G\in   L^\infty(\R^d)$  provided with extra  tangential regularity on $G$,
the $L^\infty$-bound for $\nabla u_G$ in our case  relies on the following logarithmic inequality, which is simply the application of the Sobolev embedding $L^\infty(\R^d)\cap \dot W^{1,p}(\R^d)\subset \cC^{1-\frac dp}(\R^d)$ to  Theorem 7.40 of \cite{bahouri2011fourier}. 
\begin{prop}[\cite{bahouri2011fourier}, $L^\infty$-bound for double Riesz transforms provided with tangential regularity]\label[prop]{Linftyestimate}
Let $\mathcal X=(X_\upsilon)_{1\leqslant \upsilon\leqslant m}\subset \mathbb L^{\infty,p}(\R^d)$, with $d<p<\infty$, be a \textit{non-degenerate} family of $ m\in \N^*$ vector  fields 
as in \cref{Subs:TR}. 
Let $1\leqslant q<\infty$. 

There exits a constant $C=C(m,d,p,q)>0$  such that for all  $G\in L^q(\R^d)\cap L^\infty(\R^d)\cap \mathbb{L}^{p}_{\mathcal X}(\R^d)$, the following estimate holds true:
    \begin{gather}\label{ep3.34}
    \norm{\mathcal{R} \mathcal{R}  G}_{L^\infty(\R^d)}\leqslant C\norm{G}_{L^q(\R^d)}+ C\norm{G}_{L^\infty(\R^d)}\left(1+\log \left(e+\dfrac{\norm{G}_{\mathbb{L}_{\mathcal X}^p(\R^d)}}{\norm{G}_{L^\infty(\R^d)}}\right)\right).
    \end{gather}
\end{prop} 

With the aid of the above logarithmic estimate, we can propagate tangential regularity of density and achieve Lipschitz regularity of the velocity at the same time.
\begin{lemma}[Tangential regularity for the density  and Lipschitz continuity for the velocity]\label[lemma]{Aprioriestimates4}
    Assume the initial condition \eqref{ep4.29} that $\rho_0\in \mathbb L^{p}_{\mathcal X_0}(\R^d)$ where $\mathcal X_0=(X_{0,\upsilon})_{1\leqslant \upsilon \leqslant m}\subset \mathbb L^{\infty,p}(\R^d)$ is a non-degenerate family of $ m\in \N^*$ vectors fields, with $m\geqslant d-1$, 
    with $2<p<\infty$ if $d=2$ or $3<p<6$ if $d=3$.

    Then,  the family of vector fields $\mathcal X(t)=(X_\upsilon(t))_{1\leqslant \upsilon \leqslant m}$, defined as solution of the Cauchy problem \eqref{ep4.33}, is non-degenerate and  $\mathcal X(t)\subset \mathbb L^{\infty,p}(\R^d)$.
    Moreover, $\rho(t)\in \mathbb L^p_{\mathcal X(t)}(\R^d)$ and the following bounds hold true:
        \begin{gather}\label{ep4.25}
    \begin{cases}\displaystyle
    \mathcal{A}_3(t)\leqslant   \mathcal{A}_3(0) \exp\Bigl(
    C_0 \int_0^t\Bigl[1+\sqrt{t}+\norm{\nabla u(t')}_{L^\infty(\R^d)}dt'\Bigr]\Bigr),\\
    \displaystyle
          \int_0^t \norm{\nabla u(t')}_{L^\infty(\R^d)}dt'\leqslant C_0\bigl(1+\frac{\mathcal{A}_3(0)}{I(\mathcal{X}_0)}\bigr)\exp(C_0t).
  \end{cases}
    \end{gather}
\end{lemma}
 
The proof of the above Lemma is the object of \cref{Proofs4}.  \begin{rema}[Improved time regularity]\label[rema]{rema1}
 We have the following   improved time regularity, which is required for the uniqueness result, see e.g. \cite[Equation (4.31)]{danchin2020well}: For some $t_0>0$,
$$\int^{t_0}_0\sigma(t')^s \|\nabla u(t')\|_{L^\infty(\R^d)}^2 dt'<\infty,$$
where $s=4/9$ if $d=2$ and $s=1/2$ if $d=3$.
Indeed, we apply H\"older's inequality with respect to the time variable to  \eqref{ep4.24} in the proof in \cref{Proofs4}   to obtain (noticing \eqref{ep4.26} and \eqref{ep4.27})
\begin{align*}
    \int^{t}_0\|\nabla u_G(t')\|_{L^\infty(\R^d)}^2dt'
    \leq C_0t\Bigl(1+ \frac{\mathcal{A}_3(0)}{I(\mathcal{X}_0)} +t + \int^t_0\|\nabla u(t')\|_{L^\infty(\R^d)} \Bigr)^2
\end{align*}
and similarly as in the proof of \eqref{ep3.31} and \eqref{ep3.32}, we have
 \begin{gather*} 
         \int_0^t\sigma^{\tfrac{4}{9}}\norm{\nabla \widetilde u}_{L^\infty(\R^2)}^2\leqslant 
         C_0(1+t^{\frac13}),\quad
         \int_0^t\sqrt\sigma\norm{\nabla \widetilde u}_{L^\infty(\R^3)}^2\leqslant C_0.
 \end{gather*}
 \end{rema}
 
To complete the proof of \cref{th1}, we need to construct an approximate sequence $(\rho^\delta,u^\delta)_\delta$ globally defined in time  that converges to $(\rho,u)$, the unique solution of  \eqref{ep0.1}. Once this is done, we will have obtained a sequence $(\rho^{(\nu)}, u^{(\nu)})$ of solutions to \eqref{ep0.1}, and the last step will be to justify that this sequence converges to some $(\varrho, v)$ that solves the inhomogeneous incompressible model. This is the purpose of \cref{approximate}, and, as we will see, the local solutions constructed in \cite{danchin2020well} cannot serve as building blocks. Thus, we will need to establish local well-posedness for the system \eqref{ep0.1} in \cref{local}.
\subsection{Proofs}\label[section]{proofs}

In this subsection we give the detailed proofs of \cref{Aprioriestimates1}, \cref{Aprioriestimates2}, \cref{Aprioriestimates3} and  \cref{Aprioriestimates4}. 
Before that we recall some basic  facts,  which will be used freely in the proofs below.

\subsubsection{Basic facts}\label{subss:Facts}
Under the assumption \eqref{assumption:rho,bdd}, we have    the $L^\infty((0,t),L^q(\mathbb{R}^d))$ estimate for the vibration of the pressure term  $G(t,x)=\bigl(P(\rho)\bigr)(t,x)-\widetilde P$  
given in \eqref{ep4.46} 
\begin{gather}\label{ep4.46+}
     \norm{G}_{L^\infty([0,t],L^q(\R^d))}   
    \leqslant   C^\ast ( E_0)^{\frac1q}, \quad \text{with}\quad q\in [2,\infty],
\end{gather}
where $C^\ast$ depends on $q,\rho^\ast$.
Here the case $q=\infty$ follows straightforwardly from the definition.
 Recall also the estimate \eqref{ep1.9} for $G$ by $F$: For any $l\geqslant  1$,
    \begin{gather}\label{ep1.9++}
        \dfrac{1}{\nu} \norm{G}_{L^{l+1}((0,t)\times \R^d)}^{l+1} \leqslant  C^\ast(l) E_0+ \dfrac{C}{\nu} \norm{F}_{L^{l+1}((0,t)\times \R^d)}^{l+1},\quad \forall t\in (0,\infty),
\end{gather}
where we estimated $\|H_l(\rho_0)\|_{L^1(\R^d)}$ by $C^\ast(l) E_0$.

Recall also the relations \eqref{Deltau}: $$\Delta u^j = \partial_j \dvg u+ \partial_k \rot_{jk} u,\quad j,k=1,\cdots,d,$$ and \eqref{divu,F,curlu} between $\dvg u, F, G, \rho \dot u, \rot u$:
\begin{gather}\label{ep1.6}
\dvg u=\dfrac{1}{\nu} \left(F+ G\right),\; \quad  F= -(-\Delta)^{-1}\dvg (\rho \dot u)\quad \text{and}\quad \mu \rot u= -(-\Delta)^{-1}  \rot (\rho \dot u).
\end{gather}
By use of the $L^q(\R^d)$, $q\in (1,\infty)$, $d\geqslant  2$-boundedness of Riesz operators, the following estimates follow immediately
\begin{align}
    &\|\nabla u\|_{L^q(\R^d)}\leqslant C(q,d)\bigl( \|\dvg u\|_{L^q(\R^d)}+\|\rot u\|_{L^q(\R^d)}\bigr),\label{Lp:gradient}
    \\
    &\|\nabla F\|_{L^q(\R^d)}+\mu\|\nabla\rot u\|_{L^q(\R^d)}\leqslant C(q,d)\|\rho\dot u\|_{L^q(\R^d)}.
    \label{Lp:dotu}
\end{align}

We now recall the compensated result by Coifman, Lions, Meyer, Semmes \cite{coifman1993compensated} in dimension two.
\begin{prop}[Coifman-Lions-Meyer-Semmes' Estimate for $d=2$]\label[prop]{propcom}
    Let $v,\,w\in \dot H^1(\R^2;\R)$ be two functions and let us define
    \[
    g=\det 
    \begin{pmatrix}
        \partial_1 v & \partial_2 v\\
        \partial_1 w &  \partial_2 w
    \end{pmatrix}.
    \]
    Then $g$ belongs to the Hardy space $\mathcal{H}^1(\R^2)$, whence for all $f\in BMO(\R^2)$ we have the estimate
    \[
    \bigg|\int_{\R^2}f(x) g(x)dx\bigg|\leqslant  \norm{f}_{BMO(\R^2)}\norm{\nabla v}_{L^2(\R^2)}\norm{\nabla w}_{L^2(\R^2)}.
    \]
    In particular, since $\dot H^1(\R^2)\hookrightarrow BMO(\R^2)$, for all $f\in \dot H^1(\R^2)$ we have:
    \[
    \bigg|\int_{\R^2}f(x) g(x)dx\bigg|\leqslant  \norm{\nabla f}_{L^2(\R^2)}\norm{\nabla v}_{L^2(\R^2)}\norm{\nabla w}_{L^2(\R^2)}.
    \]
\end{prop}
It helps in the proof, by virtue of  the following equalities in dimension two
\begin{gather}\label{ep2.10}
\nabla u^j \cdot \nabla u^k \partial_k u^j = \dvg u \{\abs{\nabla u}^2- \det (\nabla u)\},
\quad  \nabla u^l \cdot\partial_l u= (\dvg u)^2 -2 \det (\nabla u).
\end{gather}
Here and in the following we use Einstein's summation convention for   repeated indices, unless otherwise claimed.

\subsubsection{Proof of \cref{Aprioriestimates1} for $d=2$}\label[section]{proofs1}
This paragraph is devoted to obtaining bounds for functionals $\mathcal{A}_1$ and $\mathcal{A}_2$ as defined in \eqref{ep4.42} for $d=2$, provided with bounded density function \eqref{assumption:rho,bdd}.
The constants in the following estimates may depend on the  viscosity coefficient $\mu$ and the lower bound $\underline{\nu}$ for $\nu$, while not on the viscosity coefficient $\nu$ which will be chosen to be big.

\dem[Proof of \eqref{Aprioriestimates1:A1}]


 The functional $\mathcal{A}_1$ arises while using $\dot u$ as a test functional in the weak formulation of the momentum equation $\eqref{ep0.1}_2$. By doing so, one obtains the following equality: 
\begin{align}
    \mathcal{A}_1(t)&=\dfrac{\mu}{2}\norm{\nabla u_0}_{L^2(\R^d)}^2+\dfrac{\mu+\lambda}{2}\norm{\dvg u_0}_{L^2(\R^d)}^2
    -\mu \int_0^t \int_{\R^d} \nabla u^j \cdot \nabla u^k \partial_k u^j+\dfrac{\mu}{2} \int_0^t\int_{\R^d}\abs{\nabla u}^2 \dvg u\nonumber\\
    &+\dfrac{\mu+\lambda}{2} \int_0^t\int_{\R^d}(\dvg u)^3- (\mu+\lambda)\int_0^t\int_{\R^d}\dvg u\nabla u^l \cdot \partial_l u+\int_0^t\int_{\R^d}\nabla u^l\cdot \partial_l u\, G
    \nonumber\\
    &+\int_{\R^d}\dvg u(s) G(s)\bigg|_{s=0}^{s=t} + \int_0^t\int_{\R^d}(\rho P'(\rho)-P(\rho)+\widetilde P)(\dvg u)^2.\label{ep1.2}
\end{align}

\noindent\textbf{Step 1: Reformulation of the energy equality.} 
In the following lines, we will reformulate the terms  appearing in the right hand side above by use of \eqref{ep1.6} and \eqref{ep2.10}.

By \eqref{ep1.6} and \eqref{ep2.10},  the sum of the third and the fourth terms on the right hand side of \eqref{ep1.2} can be reduced as follows:
\begin{align}
  &\mu\int_0^t \int_{\R^2} \dvg u \left[\det (\nabla u)-\dfrac{1}{2}\abs{\nabla u}^2\right]
   \nonumber\\
  &    =\dfrac{\mu}{\nu}\int_0^t \int_{\R^2}F \det (\nabla u)-\dfrac{\mu}{2\nu}\int_0^t \int_{\R^2}F\abs{\nabla u}^2 
    +\dfrac{\mu}{\nu}\int_0^t \int_{\R^2}G\left[\det (\nabla u)
    -\dfrac{1}{2}\abs{\nabla u}^2\right],\label{ep1.10}
\end{align}
and similarly, the sum of the second and the third terms of the second line of \eqref{ep1.2}  read
\begin{align}
&-\dfrac{\mu+\lambda}{\nu}\int_0^t\int_{\R^2} F\nabla u^l \cdot\partial_l u
 + \dfrac{\mu}{\nu}\int_0^t\int_{\R^2}\nabla u^l\cdot \partial_l u\,G\nonumber\\
&=2\dfrac{\mu+\lambda}{\nu}\int_0^t\int_{\R^2} F\det (\nabla u)-\dfrac{\mu+\lambda}{\nu}\int_0^t\int_{\R^2} F (\dvg u)^2 
 + \dfrac{\mu}{\nu}\int_0^t\int_{\R^2}\nabla u^l \cdot\partial_l u\,G.\label{ep1.4}
\end{align}
Now we pack the first term of the second line of \eqref{ep1.2} and the middle term in the above \eqref{ep1.4} and use \eqref{ep1.6} to get
\begin{align}\label{F,P,quadrad}
     \dfrac{\mu+\lambda}{2} \int_0^t\int_{\R^2}(\dvg u)^3-\dfrac{\mu+\lambda}{\nu}\int_0^t\int_{\R^2} F (\dvg u)^2 =  -\dfrac{\mu+\lambda}{2\nu^2} \int_0^t\int_{\R^2}\left(F^2- G^2\right) \dvg u .
\end{align}

\noindent\textbf{Step 2: Estimates for the integrals in terms of $E_0, E_0^\nu$ and $L^4_{t,x}$-norms of $(\nabla u, G, F)$.}
We are ready to estimate  all the terms above.
\begin{itemize}
\item 
With the help of  \cref{propcom} and \eqref{Lp:dotu}, 
 the first terms of the right hand side of \eqref{ep1.10} and \eqref{ep1.4} can be estimated, as follows:
\begin{align}
    \bigg|\dfrac{3\mu+2\lambda}{\nu}\int_0^t\int_{\R^2} F\det (\nabla u)\bigg|&\leqslant C\int_0^t \norm{\nabla F}_{L^2(\R^2)}\norm{\nabla u}_{L^2(\R^2)}^2\nonumber\\
    &= C\int_0^t \norm{\rho \dot u}_{L^2(\R^2)}\norm{\nabla u}_{L^2(\R^2)}^2
    \nonumber\\
    &\leqslant \eta \int_0^t \norm{\sqrt{\rho}\dot u}_{L^2(\R^2)}^2+ \dfrac{C\rho^*}{4\eta}\int_0^t \norm{\nabla u}_{L^2(\R^2)}^4,\hbox{ by Young's inequality}, \label{ep1.15}
\end{align}
for some $\eta>0$ small enough to be determined later. 

\item By the energy balance \eqref{ep1.1}  and the upper bound $\rho^\ast$ for the density  \eqref{assumption:rho,bdd}, the last terms of  \eqref{ep1.2} and the terms involving the pressure deviation $G$ in \eqref{ep1.10} and\eqref{ep1.4}    can be bounded as follows
\begin{align}
    \bigg|\int_0^t\int_{\R^d}(\rho P'(\rho)-P(\rho)+\widetilde P)(\dvg u)^2\bigg|&+\dfrac{\mu}{\nu}\bigg|\int_0^t \int_{\R^2}G \left[\det (\nabla u)
    -\dfrac{1}{2}\abs{\nabla u}^2\right]\bigg|\nonumber\\
    &+ \dfrac{\mu}{\nu}\bigg|\int_0^t\int_{\R^2}\nabla u^l \partial_l u\,G\bigg|\leqslant \dfrac{ C^\ast}{\nu} E_0.\label{ep1.12}
\end{align}

\item 
Next, the middle term of \eqref{ep1.10} is:
\begin{align}
\dfrac{\mu}{2\nu}\bigg|\int_0^t \int_{\R^2}F\abs{\nabla u}^2\bigg| &\leqslant \dfrac{C}{\nu} \int_0^t \norm{F}_{L^4(\R^2)}\norm{\nabla u}_{L^{4}(\R^2)}\norm{\nabla u}_{L^2(\R^2)}\nonumber\\
&\leqslant C\int_0^t\left[ \dfrac{1}{\nu^{5/2}}\norm{F}_{L^4(\R^2)}^4+\dfrac{1}{\nu^{3/2}}\norm{\nabla u}_{L^4(\R^2)}^4\right]+C\int_0^t \norm{\nabla u}_{L^{2}(\R^2)}^2.
\end{align}

\item 
 The    term in  \eqref{F,P,quadrad} can be estimated as follows 
\begin{align}
\dfrac{\mu+\lambda}{2\nu^2}\bigg|\int_0^t\int_{\R^2}\left(F^2- G^2\right) \dvg u\bigg| \leqslant \dfrac{C}{\nu^3}\int_0^t\left[\norm{F}_{L^4(\R^2)}^4+\norm{G}_{L^4(\R^2)}^4\right]
 +C\nu\int_0^t\norm{\dvg u}_{L^2(\R^2)}^2.\label{ep1.5}
\end{align} 

\item 
It only remains the first term in the last line of \eqref{ep1.2} which can be bounded as follows:
\begin{align*}
    \Bigg|\int_{\R^d}\dvg u(s) G(s)\bigg|_{s=0}^{s=t}\Bigg|
    &\leqslant \eta \nu\norm{\dvg u(t)}_{L^2(\R^2)}^2+ \dfrac{C}{4\eta\nu}\norm{G(t)}_{L^2(\R^2)}^2\\
    &\qquad + \dfrac{C}{\nu}\norm{G_0}_{L^2(\R^2)}^2+C\nu\norm{\dvg u_0}_{L^2(\R^2)}^2,
\end{align*}
with the second term of the right hand side controlled by the initial energy as in \eqref{ep4.46+}.
\end{itemize}

We combine all of these estimates and we choose $\eta$ small in order to obtain the following:
\begin{multline*}
    \mathcal{A}_1(t)\leqslant C \left(1+ \dfrac{C^\ast}{\nu}\right)E_0^\nu+ C\rho^*\int_0^t \norm{\nabla u}_{L^2(\R^2)}^4
    +\dfrac{C}{\nu^{3/2}}\int_0^t \left[\norm{\nabla u}_{L^4(\R^2)}^4\right.\left.
    +\dfrac{1}{\nu}\norm{F}_{L^4(\R^2)}^4+\dfrac{1}{\nu^{3/2}}\norm{G}_{L^4(\R^2)}^4\right]
\end{multline*}
where $E_0^\nu$ is given in \eqref{E0-lambda}.
Hence Gronwall Lemma yields:
\begin{gather}\label{ep1.13}
\mathcal{A}_1(t)\leqslant  \left[  C^\ast E_0^\nu+\dfrac{C}{\nu^{3/2}}\int_0^t \left(\vphantom{\dfrac{1}{\nu^{3/2}}}\norm{\nabla u}_{L^4(\R^2)}^4\right.\left.
    +\dfrac{1}{\nu}\norm{F}_{L^4(\R^2)}^4+\dfrac{1}{\nu^{3/2}}\norm{G}_{L^4(\R^2)}^4\right)\right]\exp\left(C^* E_0\right).
\end{gather}

\noindent\textbf{Step 3: Final estimates.}
The next step is devoted to obtaining estimate for the $L^4((0,t)\times \R^2)$ norm of the velocity gradient $\nabla u$, the pressure deviation $G$ and the effective flux $F$.
\begin{itemize}
    \item $L^4$-Estimate for  $G$.
    Recall \eqref{ep1.9++} with $l=3$:
    \begin{gather}\label{ep1.9+}
      \dfrac{1}{\nu} \norm{G}_{L^4((0,t)\times \R^d)}^4\leqslant  C^\ast E_0+ \dfrac{C}{\nu} \norm{F}_{L^4((0,t)\times \R^d)}^4.
\end{gather} 

\item $L^4$-Estimate for  $F$.
The $L^4((0,t)\times \R^2)$-norm of the effective flux $F$ follows from Gagliarodo-Nirenberg's inequality $\|f\|_{L^4(\R^2)}^2\lesssim \|f\|_{L^2(\R^2)}\|\nabla f\|_{L^2(\R^2)}$ and \eqref{Lp:dotu}   as
\begin{align*}
   \norm{F}_{L^4( \R^2)}
   \leqslant C\|F\|_{ L^2(\R^2)}^{\frac12}
   \|\nabla F\|_{L^2  (\R^2) }^{\frac12} 
 \leqslant  
  C\|F\|_{ L^2(\R^2) }^{\frac12}\|\rho\dot u\|_{L^2(  \R^2)}^{\frac12}, 
\end{align*}
which can be bounded further by virtue of \eqref{ep1.6} and the definition of $E(t), \mathcal{A}_1(t)$ by
\begin{align*}
   \norm{F}_{L^4((0,t)\times \R^2)}^4
 &  \leqslant   
  C\int^t_0 (\nu\|\dvg u\|_{ L^2(\R^2) }
  +\|G\|_{ L^2(\R^2) })^2\|\rho\dot u\|_{L^2(  \R^2)}^2
  \\
  &\leqslant C^\ast(\nu\mathcal{A}_1(t)+  E_0) 
  \mathcal{A}_1(t)  .
\end{align*}

\item $L^4$-Estimate for  $\nabla u$.
Similar as above for $F$, we have $L^4$-Estimate for $\rot u$:
\begin{align*}
     \norm{\rot u}_{L^4((0,t)\times \R^2)}
 &  \leqslant   
  C\|\rot u\|_{L^\infty((0,t);L^2(\R^2))} ^{\frac12}\|\rho\dot u\|_{L^2((0,t)\times \R^2)}^{\frac12}
 \leqslant C 
  \mathcal{A}_1(t)^{\frac12}.
\end{align*}
Hence, by use of \eqref{ep1.6}-\eqref{Lp:gradient}, the following inequality holds true 
\begin{align*}
   \norm{\nabla u}_{L^4((0,t)\times \R^2)}^4
&\leqslant C\left(\norm{\dvg u}_{L^4((0,t)\times \R^2)}^4+\norm{\rot u}_{L^4((0,t)\times \R^2)}^4\right)
\\
&\leqslant \frac{C}{\nu^4}(\|F\|_{L^4((0,t)\times \R^2)}^4+\|G\|_{L^4((0,t)\times \R^2)}^4)+C\norm{\rot u}_{L^4((0,t)\times \R^2)}^4\\
&\leqslant 
\frac{C^\ast}{\nu^3} E_0
 + C^\ast \mathcal{A}_1(t)\Bigl(  E_0+\mathcal{A}_1(t)\Bigr) .
\end{align*}
\end{itemize}
Finally, we go back to \eqref{ep1.13} and we have \eqref{Aprioriestimates1:A1}.
\enddem
\dem[Proof of \eqref{Aprioriestimates1:A2}]
We turn to providing bound for the second functional $\mathcal{A}_2$ for  $d=2$. For this purpose, by \eqref{ep1.6}, we rewrite the momentum equation $\eqref{ep0.1}_2$
as follows 
\begin{gather}\label{ep1.36}
    \rho \dot u=\mu \Delta u+\dfrac{\mu+\lambda}{\nu}\nabla F-\dfrac{\mu}{\nu} \nabla G.
\end{gather}
We apply the operator $\dpt\; \cdot\;+ \dvg (\;\;\cdot \;u)$ to \eqref{ep1.36} and we obtain the following equation for the material derivative 
of the velocity:
\begin{multline}\label{ep1.23}
    \dpt (\rho \dot u^j)+\dvg (\rho \dot u^j u)
    -\mu \Delta \dot u^j
    -\dfrac{\mu+\lambda}{\nu} \partial_j \dot F\\
    = -\mu \partial_k (\nabla u^j \cdot \partial_k u)+ \mu \partial_k (\partial_k u^j \dvg u)-\mu \dvg (\partial_k u^j \partial_k u)\\
    + \dfrac{\mu+\lambda}{\nu} \partial_j (F \dvg u)-\dfrac{\mu+\lambda}{\nu}\dvg (F \partial_j u)
    \\
    + \dfrac{\mu}{\nu}\partial_j\left((\rho P'(\rho)- P(\rho)+\widetilde P)\dvg u\right)+ \dfrac{\mu}{\nu}\dvg (\partial_j u G),\quad j=1,\cdots,d.
\end{multline}

\noindent\textbf{Step 1: Formulation of the energy equality.}
To obtain the functional $\mathcal{A}_2$, it suffices to multiply the equation above by $\sigma \dot u^j$, with $\sigma=\sigma(t)=\min\{1,t\}$, sum up $j$, and integrate it in time 
and space. 
The most delicate term is $-\dfrac{\mu+\lambda}{\nu} \partial_j \dot F$ on the left hand side of \eqref{ep1.23}, which gives
\begin{align*}
      -\dfrac{\mu+\lambda}{\nu}\int_0^t \int_{\R^2}\sigma\dot u^j\partial_j \dot F 
      =\dfrac{\mu+\lambda}{\nu}\int_0^t \sigma\int_{\R^2}\dot F \dvg \dot u.
\end{align*}
We first focus on this integral for a while.
Applying material derivative to \eqref{ep1.6} gives
\begin{gather}
\dvg \dot u= \dfrac{1}{\nu}\left(\dot F-\rho P'(\rho)\dvg u\right)+\nabla u^k\cdot \partial_k u,
\end{gather}
and hence
\begin{align}
    -\dfrac{\mu+\lambda}{\nu}\int_0^t \int_{\R^2}\sigma\dot u^j\partial_j \dot F &= \dfrac{\mu+\lambda}{\nu^2}\int_0^t\sigma\norm{\dot F}^2_{L^2(\R^2)} +\dfrac{\mu+\lambda}{\nu}\int_0^t \sigma\int_{\R^2}\dot F \nabla u^k\cdot \partial_k u\nonumber\\
    &-\dfrac{\mu+\lambda}{\nu^2}\int_0^t\sigma \int_{\R^2}\dot F  \rho P'(\rho)\dvg u.\label{ep1.20}
\end{align}
To conclude, testing \eqref{ep1.23} by $\sigma\dot u$ implies  
\begin{multline}\label{ep1.21}
\dfrac{1}{2}\sigma(t)\norm{\sqrt{\rho}\dot u}_{L^2(\R^2)}^2+ \mu \int_0^t\sigma\norm{\nabla \dot u}_{L^2(\R^2)}^2+\dfrac{\mu+\lambda}{\nu^2}\int_0^t\sigma\norm{\dot F}^2_{L^2(\R^2)} 
=\dfrac{1}{2}\int_0^{\sigma(t)}\norm{\sqrt{\rho}\dot u}_{L^2(\R^2)}^2+\sum_{k=1}^4 I_k ,
\end{multline}
where
\begin{align*}
&I_1=-\dfrac{\mu+\lambda}{\nu}\int_0^t \sigma\int_{\R^2}\dot F \nabla u^k\cdot \partial_k u +\dfrac{\mu+\lambda}{\nu^2}\int_0^t\sigma \int_{\R^2}\dot F  \rho P'(\rho)\dvg u,
\\
&I_2=\mu \int_0^t \sigma \int_{\R^2}
\Bigl( \nabla u^j \cdot \partial_k  u \partial_k \dot u^j- \partial_k u^j \dvg u \partial_k \dot u^j+ \partial_k u^j \partial_ku  \cdot \nabla \dot u^j\Bigr) ,
\\
& I_3= \dfrac{\mu+\lambda}{\nu}\int_0^t \sigma\int_{\R^2} 
\Bigl(-F \dvg u\dvg \dot u+F\partial_j u^k \partial_k \dot u^j\Bigr),
\\
&I_4= -\dfrac{\mu}{\nu}\int_0^t \sigma\int_{\R^2}
\Bigl( \dvg u\dvg \dot u (\rho P'(\rho)- P(\rho)+\widetilde P)  +\partial_j u^k \partial_k \dot u^j G\Bigr).
\end{align*}

\noindent\textbf{Step 2: Estimate  for  $I_1$.}
We focus first on the first integral in $I_1$, which can be reformulated by integration by parts (noticing $\sigma(0)=0$) as
\begin{align}
\int_0^t\sigma\int_{\R^2}\dot F \nabla u^k\cdot \partial_k u
&=\sigma(t)\int_{\R^2} F(t) \nabla u^k\cdot \partial_k u(t)-\int_0^{\sigma(t)}\int_{\R^2}F \nabla u^k \cdot\partial_k u-2\int_0^t\sigma\int_{\R^2} F\partial_k u\cdot\nabla \dot u^k\nonumber \\
&+2\int_0^t \sigma \int_{\R^2} F \partial_k u\cdot \nabla u^l \partial_l u^k -\int^t_0\sigma\int_{\R^2}F \dvg u \nabla u^k\cdot\partial_k u .\label{ep1.16}
\end{align}
The first term of the right hand side in \eqref{ep1.16} above can be estimated  similarly as for the derivation of \eqref{ep1.15} and \eqref{ep1.5}, using the equality $\nabla u^l \cdot\partial_l u= (\dvg u)^2 -2 \det (\nabla u)$ in \eqref{ep2.10}:
\begin{align*}
 \bigg|\sigma(t)\int_{\R^2} F(t) (\nabla u^k\cdot \partial_k u)(t)\bigg|&\leqslant C\sigma(t)\norm{\rho \dot u}_{L^2(\R^2)}\norm{\nabla u}_{L^2(\R^2)}^2
+  \sigma(t)\dfrac{1}{\nu^2}\bigg|\int_{\R^2} F(t) (F(t)+ G(t))^2\bigg|\\
&\leqslant  \eta \sigma(t)\norm{\sqrt{\rho} \dot u}_{L^2(\R^2)}^2+ \sigma(t)\dfrac{C^2\rho^*}{4\eta}\norm{\nabla u}_{L^2(\R^2)}^4\\
&\quad +   \sigma(t)\dfrac{C}{\nu^2}\left(\norm{F(t)}_{L^3(\R^2)}^3 + \norm{G(t)}_{L^3(\R^2)}^3\right).
\end{align*}
Exactly as in the derivation of \eqref{ep1.15} and \eqref{ep1.5}, the second integral of the right hand side of
\eqref{ep1.16} can be estimated as follows:
\begin{multline}
 \bigg|\int_0^{ \sigma(t)}\int_{\R^2}F \nabla u^k \cdot\partial_k u\bigg| \leqslant C E_0+ C\int_0^{\sigma(t)} \norm{\sqrt{\rho}\dot u}_{L^2(\R^2)}^2+ C\rho^*\int_0^{\sigma(t)} \norm{\nabla u}_{L^2(\R^2)}^4\\
+\dfrac{C}{\nu^3}\int_0^{\sigma(t)} \left[\norm{F}_{L^4(\R^2)}^4+\norm{G}_{L^4(\R^2)}^4\right].
\end{multline}
In order to estimate the third term of the right hand side of \eqref{ep1.16}, we write:
\begin{gather}\label{ep1.18}
\partial_j u^k\partial_k \dot u^j=\dvg u\dvg \dot u-(\partial_1 u^1\partial_2 \dot u^2-\partial_2 u^1\partial_1 \dot u^2)-(\partial_2 u^2\partial_1 \dot u^1-\partial_1 u^2 \partial_2 \dot u^1)
\end{gather}
in such a way that, after making use of the compensated result \cref{propcom} and Young's inequality, we have:
\begin{multline}\label{ep1.19}
 \bigg|\int_0^t\sigma\int_{\R^2} F\partial_k u\cdot\nabla \dot u^k\bigg|
 \leqslant  
 \eta \int_0^t \sigma\norm{\nabla \dot u}_{L^2(\R^2)}^2+ \dfrac{C\rho^*}{\eta} \int_0^t \sigma\norm{\sqrt{\rho}\dot u}_{L^2(\R^2)}^2\norm{\nabla u}_{L^2(\R^2)}^2\\
+\dfrac{C}{\eta}\int_0^t\left[\dfrac{1}{\nu^2}\norm{F}_{L^4(\R^2)}^4+ \dfrac{1}{\nu^2}\norm{G}_{L^4(\R^2)}^4\right]. 
\end{multline}
Thanks to \eqref{ep1.6} and the equality
$\partial_k u^l \nabla u^k \cdot\partial_l u=\dvg u\{(\dvg u)^2- 3\det (\nabla u)\}
$ in \eqref{ep2.10}, we have the following estimate for the last two terms of \eqref{ep1.16}:
\begin{gather*}
2 \bigg|\int_0^t \sigma\int_{\R^2} F \partial_k u\cdot \nabla u^l \partial_l u^k \bigg|+\bigg|\int^t_0\sigma\int_{\R^2}F \dvg u \nabla u^k\cdot\partial_k u\bigg|
\\
\leqslant C\int_0^t\sigma\left[\norm{\nabla u}_{L^4(\R^2)}^4+\dfrac{1}{\nu^2}\norm{F}_{L^4(\R^2)}^4+ \dfrac{1}{\nu^2}\norm{G}_{L^4(\R^2)}^4\right].
\end{gather*}

Finally, owing to H\"older's and Young's inequalities the second integral in $I_0$ can be estimated as follows: 
\begin{gather*}
    \dfrac{\mu+\lambda}{\nu^2}\bigg|\int_0^t\sigma \int_{\R^2}\dot F  \rho P'(\rho)\dvg u\bigg|\leqslant  \eta \dfrac{\mu+\lambda}{\nu^2}\int_0^t \sigma\norm{\dot F}_{L^2(\R^2)}^2+ \dfrac{ C^\ast}{\nu^2\eta} E_0.
\end{gather*} 
Gathering the above computations, we obtain the following estimate for $I_1$:
\begin{align*}
 |I_1|   &\leqslant
    C (1+\frac{C ^\ast }{\nu^2\eta})E_0
    +\eta \sigma(t)\norm{\sqrt{\rho} \dot u}_{L^2(\R^2)}^2
    +\eta \int_0^t \sigma\norm{\nabla \dot u}_{L^2(\R^2)}^2
    +\eta \dfrac{\mu+\lambda}{\nu^2}\int_0^t \sigma\norm{\dot F}_{L^2(\R^2)}^2\\
    &+C\int_0^{\sigma(t)} \norm{\sqrt{\rho}\dot u}_{L^2(\R^2)}^2+ \sigma(t)\dfrac{C ^*}{ \eta}\norm{\nabla u}_{L^2(\R^2)}^4+  \dfrac{C\rho^*}{\eta} \int_0^t \sigma\norm{\sqrt{\rho}\dot u}_{L^2(\R^2)}^2\norm{\nabla u}_{L^2(\R^2)}^2\\
    &+ C\rho^*\int_0^{\sigma(t)} \norm{\nabla u}_{L^2(\R^2)}^4+   \sigma(t)\dfrac{C}{\nu^2}\left(\norm{F(t)}_{L^3(\R^2)}^3 + \norm{G(t)}_{L^3(\R^2)}^3\right)\\
    &+C\int_0^t\sigma\left[\norm{\nabla u}_{L^4(\R^2)}^4+\dfrac{1}{\eta \nu^2}\norm{F}_{L^4(\R^2)}^4+ \dfrac{1}{\eta \nu^2}\norm{G}_{L^4(\R^2)}^4\right].
\end{align*}

\noindent\textbf{Step 3: Final estimates.}
We now turn to the estimate of  the last terms $\displaystyle\sum_{k=2}^4 I_k$ in \eqref{ep1.21}.
By Young's inequality it is straightforward to get
\begin{align*}
    |I_2|\leqslant \eta \int_0^t \sigma\norm{\nabla \dot u}_{L^2(\R^2)}^2+ \dfrac{C}{\eta}\int_0^t \sigma\norm{\nabla u}_{L^4(\R^2)}^4.
\end{align*}
Similar as in Step 2 we have

\begin{align*}
&\abs{I_3}\leqslant \eta \int_0^t \sigma\norm{\nabla \dot u}_{L^2(\R^2)}^2
+ \dfrac{C\rho^*}{\eta} \int_0^t \sigma\norm{\sqrt{\rho}\dot u}_{L^2(\R^2)}^2\norm{\nabla u}_{L^2(\R^2)}^2
+\dfrac{C}{\nu^2\eta}\int_0^t\sigma \left[\norm{F}_{L^4(\R^2)}^4+\norm{G}_{L^4(\R^2)}^4\right],
\\
&|I_4|\leqslant \eta \int_0^t\sigma \norm{\nabla \dot u}_{L^2(\R^2)}^2+ \dfrac{C C^\ast}{\nu^2\eta} E_0.
\end{align*}
Summing up,   we have for $\eta$ small enough that
\begin{align}
    \mathcal{A}_2(t)&\leqslant 
    C (1+\frac{C ^\ast }{\nu^2 })E_0
    +C^\ast (1+E_0+\mathcal{A}_1(t))  \mathcal{A}_1(t) 
 \nonumber\\
    &+ \sigma(t)\dfrac{C}{\nu^2}\left(\norm{F(t)}_{L^3(\R^2)}^3 + \norm{G(t)}_{L^3(\R^2)}^3\right) +C\int_0^t\sigma\left[\norm{\nabla u}_{L^4(\R^2)}^4+\dfrac{1}{ \nu^2}\norm{F}_{L^4(\R^2)}^4+ \dfrac{1}{ \nu^2}\norm{G}_{L^4(\R^2)}^4\right].\label{ep3.20}
\end{align}

Recall Step 3 in the Proof of \eqref{Aprioriestimates1:A1} above, and we get similar $L^4$-Estimates with the time weight $\sigma$.
We have the following similar as  \eqref{ep1.9+} 
\begin{gather*}
  \dfrac{1}{2\nu}\int_0^t\sigma \norm{G}_{L^4(\R^2)}^4\leqslant  
  C^\ast E_0+ \dfrac{C}{\nu}\int_0^t \sigma\norm{F}_{L^4(\R^2)}^4,
\end{gather*}
which implies
\begin{gather}\label{ep3.21}
    \int_0^t\sigma\left[\norm{\nabla u}_{L^4(\R^2)}^4+\dfrac{1}{ \nu^2}\norm{F}_{L^4(\R^2)}^4+ \dfrac{1}{ \nu^2}\norm{G}_{L^4(\R^2)}^4\right] \leqslant 
    \dfrac{C^\ast}{\nu^3}E_0+C^\ast(E_0+\mathcal{A}_1(t))\mathcal{A}_1(t). 
\end{gather}
On the other hand, we have  \eqref{ep4.46+}: $\dfrac{1}{\nu^2}\norm{G(t)}_{L^3(\R^2)}^3\leqslant \dfrac{C^\ast}{\nu^2}  E_0$  and hence  the following   thanks to  Gagliardo-Nirenberg inequality and \eqref{ep1.6}-\eqref{Lp:dotu}
\begin{align}
    \dfrac{\sigma(t)}{\nu^2}\norm{F(t)}_{L^3(\R^2)}^3
    &\leqslant C\dfrac{\sigma(t)}{\nu^2} \norm{\nabla F(t)}_{L^2(\R^2)} \norm{F(t)}_{L^2(\R^2)}^2\nonumber\\
    &\leqslant \sigma(t)\norm{\rho \dot u}_{L^2(\R^2)}\left(\dfrac{1}{\nu^2} \norm{G}_{L^2(\R^2)}^2+\norm{\dvg u}_{L^2(\R^2)}^2\right)\nonumber\\
    &\leqslant \eta \sigma(t)\norm{\sqrt{\rho}\dot u}_{L^2(\R^2)}^2+ \dfrac{C ^*}{\eta\nu^2}\left(\dfrac{1}{\nu^2} E_0^2 +(\mathcal{A}_1(t))^2 
    \right).\label{ep3.22}
\end{align}
We finally combine \eqref{ep3.20}, \eqref{ep3.21} and \eqref{ep3.22}, we choose $\eta$ small  to derive \eqref{Aprioriestimates1:A2}.
\enddem
\subsubsection{Proof of \cref{Aprioriestimates2}}\label[section]{proofs2}
This section is devoted to obtaining  bounds for functionals $\mathcal{A}_1$ and $\mathcal{A}_2$ as defined in \eqref{ep4.42}, for $d=3$.
The proof is similar as in  \cref{Aprioriestimates1}, and we adapt the estimates in three dimension, for instance, the $L^4(\R^2)$-norm is replaced by $L^6(\R^3)$-norm below.
Since \cref{propcom} does not hold in three dimension anymore,   we will simply use the Sobolev embedding $\dot H^1(\R^3)\hookrightarrow L^6(\R^3)$ in the estimates.  
\dem[Proof of \eqref{Aprioriestimates2:A1}]
We recall that the first functional appears while using $\dot u$ as a test function in the weak formulation of $\eqref{ep0.1}_2$. By doing so, we obtain again \eqref{ep1.2}: 
\begin{align}
    \mathcal{A}_1(t)&=\dfrac{\mu}{2}\norm{\nabla u_0}_{L^2(\R^3)}^2+\dfrac{\mu+\lambda}{2}\norm{\dvg u_0}_{L^2(\R^3)}^2
    -\mu \int_0^t \int_{\R^3} \nabla u^j \cdot \nabla u^k \partial_k u^j+\dfrac{\mu}{2} \int_0^t\int_{\R^3}\abs{\nabla u}^2 \dvg u\nonumber\\
    &+\dfrac{\mu+\lambda}{2} \int_0^t\int_{\R^3}(\dvg u)^3- (\mu+\lambda)\int_0^t\int_{\R^3}\dvg u\nabla u^l \partial_l u
   +\int_0^t\int_{\R^3}\nabla u^l \partial_l u G\nonumber\\
    &  +\int_{\R^3}\dvg u(s) G(s)\bigg|_{s=0}^{s=t}+ \int_0^t\int_{\R^3}(\rho P'(\rho)-P(\rho)+\widetilde P)(\dvg u)^2.\label{ep3.23}
\end{align}

\noindent\textbf{Step 1: Estimates in terms of $E_0^\nu$ and   $L^2_tL^6_{x}$-norms of $(\nabla u, G, F)$.}
With the help of H\"older's inequality  the last two terms of the first line above can be straightforwardly estimated    as follows:
\begin{align}
    \bigg|-\mu \int_0^t \int_{\R^3} \nabla u^j \cdot \nabla u^k \partial_k u^j+\dfrac{\mu}{2} \int_0^t\int_{\R^3}\abs{\nabla u}^2 \dvg u\bigg|&\leqslant C \|\nabla u\|_{L^3((0,t)\times\R^3)}^3.\nonumber
\end{align}
Similarly, together with the relation between $\dvg u$ and $F, G$ as well as H\"older's inequality, the whole second line in \eqref{ep3.23} can be bounded by 
\begin{align*}
    C\int^t_0 \Bigl( \|\dvg u\|_{L^2(\R^3)}
     \|\dvg u\|_{L^3(\R^3)}
     \|(F,G)\|_{L^6(\R^3)}
     +\|\nabla u\|_{L^2(\R^3)}
     \|\nabla  u\|_{L^3(\R^3)}
     \|(F,\frac{1}{\nu}G)\|_{L^6(\R^3)}\Bigr)
\end{align*}
which is, by virtue of the   interpolation inequality $\| f\|_{L^3(\R^3)}\lesssim \|  f\|_{L^2(\R^3)}^{1/2}\|  f\|_{L^6(\R^3)}^{1/2}$, bounded by
\begin{align*}
    C\int^t_0 \Bigl( \|\sqrt{\nu}\dvg u\|_{L^2(\R^3)}^{\frac32}
     \frac{1}{\nu^{\frac54}}\|(F,G)\|_{L^6(\R^3)}^{\frac32}
     +\|\nabla u\|_{L^2(\R^3)}^{\frac32}
     \|\nabla  u\|_{L^6(\R^3)}^{\frac12}
     \|(F,\frac{1}{\nu}G)\|_{L^6(\R^3)}\Bigr).
\end{align*}
This can be further estimated by Young's inequality by the following with some small constant $\eta>0$:
\begin{align*}
    \eta \int^t_0 \Big\|\Big(\nabla u, F,\frac{1}{\nu^{\frac56}}G\Big)\Big\|_{L^6(\R^3)}^2  
    +\frac{C}{\eta}\int^t_0   \|(\sqrt{\nu}\dvg u,\nabla u)\|_{L^2(\R^3)}^6 .
\end{align*} 
By the same argument the whole third line in \eqref{ep3.23} can be bounded by 
\begin{align*}
    &\eta\nu\|\dvg u(t)\|_{L^2(\R^3)}^2+\frac{C}{\eta}\|G(t)\|_{L^2(\R^3)}^2
    +\|\dvg u_0\|_{L^2(\R^3)}\|G_0\|_{L^2(\R^3)}
    +C^\ast \|\dvg u\|_{L^2((0,t)\times\R^3)}^2\\
    &\leqslant \eta\nu\|\dvg u(t)\|_{L^2(\R^3)}^2+(\frac{C^\ast}{\eta}+\frac{C^\ast}{\sqrt{\nu}})E_0^\nu.
\end{align*}
Summing up,   we obtain the following  by further applying   the interpolation inequality and then Young's inequality to  $\|\nabla u\|_{L^3((0,t)\times\R^3)}^3$:
\begin{equation}\label{ep1.35}
\mathcal{A}_1(t)\leqslant  \frac{C^\ast}{\eta} E_0^\nu + 
 \eta\Bigl( \nu\|\dvg u(t)\|_{L^2(\R^3)}^2+\int^t_0 \Big\|\Big(\nabla u, F,\frac{1}{\nu^{\frac56}}G\Big)\Big\|_{L^6(\R^3)}^2  \Bigr)
    +\frac{C}{\eta}\int^t_0   \|(\sqrt{\nu}\dvg u,\nabla u)\|_{L^2(\R^3)}^6 ,
\end{equation}
with some small parameter $\eta\leqslant 1$ to be determined below.

\noindent\textbf{Step 2: Final Estimates.}
We now turn to the estimate of the $L^2((0,t); L^6(\R^3))$-norm $(\nabla u, G, F)$, similar as in Step 3 in the proof of \eqref{Aprioriestimates1:A1}.
\begin{itemize}
    \item $L^2_tL^6_x$-Estimate for  $G$.
Recall \eqref{ep1.40} with $l=5$:  
 \begin{gather*} 
 \dfrac{d}{dt}\norm{H_5(\rho)}_{L^1(\R^3)}+\dfrac{1}{ \nu}\norm{G}_{L^6(\R^3)}^6\leqslant \dfrac{1}{\nu} \norm{G}_{L^6(\R^3)}^5\norm{F}_{L^6(\R^3)} .
 \end{gather*}
We define the time-dependant function $h(t):=\|H_5(\rho)\|_{L^1(\R^3)}^{\frac13}$.
Thanks to the equivalence between $\norm{H_5(\rho)}_{L^1(\R^3)}$ and $\norm{G}_{L^6(\R^3)}^6$, we have 
$ \norm{G}_{L^6(\R^3)}^2\sim_{C^\ast}h$.
 \[
 3 \dfrac{d}{dt}h + \dfrac{1}{ \nu C^\ast}h \leqslant \dfrac{1 }{\nu} h^{\frac12}\norm{F}_{L^6(\R^3)},
 \]
 which together with Young's inequality yields:
 \begin{gather*}
 3\sup_{[0,t]}h+\dfrac{1}{2\nu C^\ast}\int_0^th\leqslant 3 {h(0)} + 2\dfrac{ C^\ast }{\nu}\int_0^t\norm{F}_{L^6(\R^3)}^2.
 \end{gather*}
 Finally we find the following estimate for $G$ in terms of $E_0$ and $F$:
\begin{gather}\label{ep4.47}
\dfrac{1}{\nu}\int^t_0 \norm{G}_{ L^6(\R^3)}^2\leqslant C^\ast  E_0^{\frac13}+\dfrac{C^\ast}{\nu}  \int_0^t\norm{F}_{L^6(\R^3)}^2 .
\end{gather}

\item $L^2_tL^6_x$-Estimate for  $F$.
We use  the Sobolev embedding $\| g\|_{L^6(\R^3)}\lesssim \|\nabla g\|_{L^2(\R^3)}$ to bound the $L^6_x$-norm of $F$ by \eqref{Lp:dotu}: 
 \begin{align*}
\int_0^t\norm{F}_{L^6(\R^3)}^2
    \leqslant  C\int^t_0\|\nabla F\|_{L^2(\R^3)}^2
    =C\int^t_0\|\rho\dot u\|_{L^2(\R^3)}^2.
\end{align*}

\item $L^2_tL^6_x$-Estimate for  $\nabla u$. 
 Similarly, by use of \eqref{ep1.6}-\eqref{Lp:gradient}-\eqref{Lp:dotu} 
 the following inequality holds true  
\begin{align*}
 \int^t_0\|\nabla u\|_{L^6(\R^3)}^2 
&\leqslant C\left( \int^t_0\|\dvg u\|_{L^6(\R^3)}^2 
+ \int^t_0\|\rot u\|_{L^6(\R^3)}^2 \right)
\\
&\leqslant \frac{C}{\nu^2} \int^t_0\|(F,G)\|_{L^6(\R^3)}^2 
+C\int^t_0\|\rho\dot u\|_{L^2(\R^3)}^2.
\end{align*}
\end{itemize}
Finally   \eqref{Aprioriestimates2:A1} follows from \eqref{ep1.35} when we choose $\eta$ small enough.

\enddem
\dem[Proof of \eqref{Aprioriestimates2:A2}]
Here we derive estimate for the functional $\mathcal{A}_2$ as defined in \eqref{ep4.42} for $d=3$. We recall that it appears while rewriting the equation on the form \eqref{ep1.36}, next applying the operator $\dpt \cdot +\dvg (\, \cdot\, u)$ in order to obtain \eqref{ep1.23} which we test with the material derivative of the velocity $\dot u$. By doing so, we obtain \eqref{ep1.21}:
\begin{multline}\label{ep1.38}
\frac12\sigma(t)\norm{\sqrt{\rho}\dot u}_{L^2(\R^3)}^2+ \mu \int_0^t\sigma\norm{\nabla \dot u}_{L^2(\R^3)}^2
+\dfrac{\mu+\lambda}{\nu^2}\int_0^t\sigma\norm{\dot F}^2_{L^2(\R^2)} 
=\dfrac{1}{2}\int_0^{\sigma(t)}\norm{\sqrt{\rho}\dot u}_{L^2(\R^2)}^2 +\sum_{k=1}^4 I_k ,
\end{multline}
where $I_k$, $k=1,2,3,4$ are given as in \eqref{ep1.21}.

\noindent\textbf{Step 1 Estimates for $I_1$.}
By the identity \eqref{ep1.16}  and H\"older's inequality we achieve the estimates
\begin{align*}
   \left| \int_0^t\sigma\int_{\R^3}\dot F \nabla u^k \cdot \partial_k u\right| 
   &\lesssim  \sigma(t)\norm{F}_{L^6(\R^3)}\norm{\nabla u}_{L^{2}(\R^3)}\norm{\nabla u}_{L^3(\R^3)}
   +\int_0^{\sigma(t)} \norm{F}_{L^6(\R^3)}\norm{\nabla u}_{L^2(\R^3)}\norm{\nabla u}_{L^3(\R^3)}
   \\
&   + \int_0^{t}\sigma \norm{F}_{L^6(\R^3)}\norm{\nabla u}_{L^3(\R^3)}\norm{\nabla \dot u}_{L^2(\R^3)}
   +\int_0^t\sigma\norm{F}_{L^6(\R^3)}\norm{\nabla u}_{L^6(\R^3)}^2\norm{\nabla u}_{L^2(\R^3)}.
\end{align*} 
Similar as in the proof of \eqref{Aprioriestimates2:A1} above, we use  the interpolation $\|\nabla u\|_{L^3(\R^3)}\lesssim \|\nabla u\|_{L^2(\R^3)}^{1/2}\|\nabla u\|_{L^6(\R^3)}^{1/2}$ and Young's inequality  to derive (noticing $\sigma(t)\leqslant 1$)
\begin{align*}
   &\left| \int_0^t\sigma\int_{\R^3}\dot F \nabla u^k \cdot \partial_k u\right|
    \\
& \leqslant C\bigl(\sqrt{\sigma(t)}\norm{F}_{L^6(\R^3)}\big) 
   \norm{\nabla u}_{L^{2}(\R^3)}^{\frac32}\big(\sqrt{\sigma(t)}\norm{\nabla u}_{L^6(\R^3)} \big)^{\frac12}  +\int^{\sigma(t)}_0 \norm{F}_{L^6(\R^3)} 
    \norm{\nabla u}_{L^6(\R^3)}  ^{\frac12}
    \norm{\nabla u}_{L^{2}(\R^3)}^{\frac32} 
    \\
&    +\int^t_0 
\norm{\sqrt{\sigma}F}_{L^6(\R^3)}  \|\nabla u\|_{L^6(\R^3)}^{\frac12}\|\nabla u\|_{L^2(\R^3)}^{\frac12}\|\sqrt{\sigma}\nabla\dot u\|_{L^2(\R^3)}
    +\int^t_0\sigma \norm{F}_{L^6(\R^3)}\norm{\nabla u}_{L^6(\R^3)}^2\norm{\nabla u}_{L^2(\R^3)}
   \\
   &\leqslant 
   \eta\Bigl(  \sigma(t) \|(\nabla u, F)\|_{L^6(\R^3)}^2
   +\int^{t}_0\sigma \|  \nabla\dot u\|_{L^2(\R^3)}^2\Bigr)
   \\
&\quad     +\frac{C}{\eta}\Bigl( \|\nabla u\|_{L^2(\R^3)}^6 
   +\int_0^t \sigma\bigl( \norm{F}_{L^6(\R^3)}^2\norm{\nabla u}_{L^6(\R^3)}+\norm{F}_{L^6(\R^3)}\norm{\nabla u}_{L^6(\R^3)}^2\bigr)\norm{\nabla u}_{L^2(\R^3)}\Bigr)
  \\
&\quad  +\Bigl( \int^{\sigma(t)}_0\| F \|_{L^6(\R^3)}^2\Bigr)^{\frac12}
\Bigl( \int^{\sigma(t)}_0 \bigl( \|\nabla u\|_{L^6(\R^3)}^ 2
+\|\nabla u\|_{L^2(\R^3)}  ^6\bigr)\Bigr)^{\frac12} .
\end{align*} 
By the Sobolev embedding $\dot H^1(\R^3)\hookrightarrow L^6(\R^3)$ and \eqref{ep1.6}-\eqref{Lp:gradient}-\eqref{Lp:dotu}-\eqref{ep4.47} 
we have 
\begin{align*}
&\|(F,\rot u)\|_{L^6(\R^3)}\leqslant C{\rho^\ast}\|\sqrt{\rho}\dot u\|_{L^2(\R^3)},
\\
&\frac1\nu\int^t_0\|G\|_{L^6(\R^3)}^2
\leqslant C ^\ast E_0^{\frac13}+\frac{C ^\ast}{\nu}\int^t_0\|\sqrt{\rho}\dot u\|_{L^2(\R^3)}^2
\leqslant C ^\ast E_0^{\frac13}+\frac{C ^\ast}{\nu}\mathcal
A_1(t),
\\
&    \|\nabla u\|_{L^6(\R^3)}
   \leqslant C(\|\dvg u\|_{L^6(\R^3)}+\|\rot u\|_{L^6(\R^3)})
 \leqslant C{\rho^\ast}\|\sqrt{\rho}\dot u\|_{L^2(\R^3)}
+\frac{C}{\nu} \|G\|_{L^6(\R^3)},
\end{align*}
which implies
\begin{align*}
    &\sigma(t)\|(\nabla u, F)\|_{L^6(\R^3)}^2
    \leqslant C^\ast \mathcal{A}_2(t)+\frac{C}{\nu^2}E_0^{\frac13},
    \\
    &\int_0^t\sigma \norm{F}_{L^6(\R^3)} \norm{(F,\nabla u)}_{L^6(\R^3)}^2 \norm{\nabla u}_{L^2(\R^3)}
    \leqslant C^\ast \mathcal{A}_1(t)^{\frac12}\mathcal{A}_2(t)^{\frac12}
    \bigl(\mathcal{A}_1(t)+\frac{C^\ast}{\nu}E_0^{\frac13}\bigr),
    \\
    &\Bigl( \int^{\sigma(t)}_0\| F \|_{L^6(\R^3)}^2\Bigr)^{\frac12}
\Bigl( \int^{\sigma(t)}_0 \bigl( \|\nabla u\|_{L^6(\R^3)}^ 2
+\|\nabla u\|_{L^2(\R^3)}  ^6\bigr)\Bigr)^{\frac12}
    \leqslant C^\ast \mathcal{A}_1(\sigma(t))^{\frac12}
    \bigl(  \mathcal{A}_1(\sigma(t))+\frac{1}{\nu}E_0^{\frac13}+E_0\mathcal{A}_1(\sigma(t))^2\bigr)^{\frac12} .
\end{align*}
It holds  by Young's inequality that
\[
\dfrac{\mu+\lambda}{\nu^2}\bigg|\int_0^t\sigma \int_{\R^3}\dot F  \rho P'(\rho)\dvg u\bigg|\leqslant 
\eta\dfrac{\mu+\lambda}{2\nu^2} \int_0^t \sigma\norm{\dot F}_{L^2(\R^3)}^2+\dfrac{ C^\ast}{\nu^2\eta} E_0.
\]
Gathering all of these computations and using Young's inequality, we have:
\begin{align}
    |I_1|&\leqslant  C^\ast\Bigl( \eta\mathcal{A}_2(t)   +\frac{1}{\nu}E_0^{\frac13}
    +\frac{1 }{\eta^2}\mathcal{A}_1(t)^3 
    +\frac{1}{\eta}\mathcal{A}_1(t) 
     \frac{1}{\nu^2} E_0^{\frac23} +   \mathcal{A}_1(\sigma(t)) 
     +  E_0\mathcal{A}_1(\sigma(t))^{\frac32}\Bigr).
\end{align}

\noindent\textbf{Step 2. Final estimates.}
Recall the definitions of $I_k$, $k=2,3,4$ given in \eqref{ep1.21}, and we can estimate them similarly as in Step 1 as follows: 
\begin{align*}
    |I_2+I_3+I_4|
    &\leqslant C\int^t_0\sigma\|\nabla\dot u\|_{L^2(\R^3)}\|(\nabla u,F)\|_{L^6(\R^3)}\|\nabla u\|_{L^3(\R^3)}
    +C^\ast \int^t_0\sigma \|\nabla\dot u\|_{L^2(\R^3)}\|\nabla u\|_{L^2(\R^3)}
    \\
    &\leqslant \eta\int^t_0\sigma\|\nabla\dot u\|_{L^2(\R^3)}^2+\frac{C^\ast }{\eta}\int^t_0
    \Bigl(\sigma\|(\nabla u, F)\|_{L^6(\R^3)}^3\|\nabla u\|_{L^2(\R^3)}+\|\nabla u\|_{L^2(\R^3)}^2\Bigr)
    \\
    &\leqslant \eta\mathcal{A}_2(t)+\frac{C^\ast}{\eta}\mathcal{A}_1(t)^{\frac12}\bigl(\mathcal{A}_2(t)^{\frac12}  \mathcal{A}_1(t)+\frac{E_0^{\frac23}}{\nu^2}+\frac{E_0^{\frac16}}{\nu^3}\mathcal{A}_1(t)\bigr)+\frac{C^\ast}{\eta}E_0.
\end{align*}
Gathering all of these computations, we obtain, after choosing $\eta$ small enough:
\begin{align*}
    \mathcal{A}_2(t)\leqslant C^\ast\Bigl(  \frac{1}{\nu}E_0^{\frac13}+E_0  +\frac{E_0^{\frac23}}{\nu^2}\mathcal{A}_1(t)^{\frac12}+ (1+\frac{1}{\nu^2}E_0^{\frac23})\mathcal{A}_1(t)+(E_0+\frac{E_0^{\frac16}}{\nu^3})\mathcal{A}_1(t)^{\frac32}  +\mathcal{A}_1(t)^3\Bigr)
\end{align*}
A further application of Young's inequality yields \eqref{Aprioriestimates2:A2}.
\enddem

\subsubsection{Proof of \cref{Aprioriestimates3}}\label[section]{Proofs3}
In this paragraph, we derive 
a priori estimate for the upper bound of the density in terms of $\mathcal{A}_1, \mathcal{A}_2$, under the assumption $\displaystyle\sup_{t,x}\rho(t,x)\leqslant \rho^\ast$. 
The basic facts in \cref{subss:Facts} will be used freely.
\dem 
With the help of the expression of $\dvg u$ in \eqref{ep1.6}, we begin by rewriting the mass equation $\eqref{ep0.1}_1$ in terms of $F,G$ as follows:
\[
\dpt \rho+ u\cdot \nabla \rho +\dfrac{\rho}{\nu}G=-\dfrac{\rho}{\nu} F.
\]
Due to the fact that the pressure is an increasing function of the density and $G=P(\rho)-\widetilde P$ the above equation yields:
\begin{gather}\label{ep1.43}
    \dpt \abs{\rho-\widetilde\rho}+ u\cdot \nabla \abs{\rho-\widetilde \rho}+ \dfrac{\rho }{\nu} \abs{G} =-\dfrac{\rho}{\nu} \text{sgn}(\rho-\widetilde\rho) F.
\end{gather}
It yields immediately  the following $L^\infty$ estimate for the density, which we use on the short time interval $[0,\sigma(t)]$ 
\begin{gather}\label{ep1.45}
    \sup_{[0,\sigma(t)]}\norm{\rho-\widetilde \rho}_{L^\infty(\R^d)}\leqslant \norm{\rho_0-\widetilde\rho}_{L^\infty(\R^d)}+\dfrac{\rho^*}{\nu}\int_0^{\sigma(t)} \norm{F(s)}_{L^\infty(\R^d)}ds.
\end{gather}

For larger time we would like to improve the $L^1_t$-norm for $\|F\|_{L^\infty_x}$ into $L^3_t$-norm, which requires less decay rate in time. From \eqref{ep1.43}  we have 
\begin{gather}\label{ep1.44}
    \dfrac{1}{3}\left(\dpt \abs{\rho-\widetilde\rho}^3+ u\cdot \nabla \abs{\rho-\widetilde \rho}^3\right)+ \dfrac{\rho }{\nu} \abs{G}\abs{\rho-\widetilde\rho}^2 =-\dfrac{\rho}{\nu} \text{sgn}(\rho-\widetilde\rho) F \abs{\rho-\widetilde\rho}^2.
\end{gather}
As well, since the pressure is an increasing function of the density such that
\begin{align*}
    |G||\rho-\widetilde\rho|^2\sim_{C^\ast}|\rho-\widetilde\rho|^3,
\end{align*}
  we derive  the following estimate from Young's inequality on larger time interval $[\sigma(t), t]$
\begin{gather}\label{ep1.46}
\sup_{[\sigma(t),t]}\norm{\rho-\widetilde\rho }_{L^\infty(\R^d)}^3\leqslant \norm{\rho(\sigma(t))-\widetilde\rho }_{L^\infty(\R^d)}^3+\dfrac{C^\ast}{\nu}\int_{\sigma(t)}^t\norm{F(s)}^3_{L^\infty(\R^d)}ds.
\end{gather}
Gathering  estimates \eqref{ep1.45} and \eqref{ep1.46}, we have:
\begin{gather}\label{ep1.49}
    \sup_{[0,t]}\norm{\rho-\widetilde\rho}_{L^\infty(\R^d)} \leqslant \norm{\rho_0-\widetilde\rho}_{L^\infty(\R^d)}+\dfrac{C^\ast}{\nu}\int_0^{\sigma(t)} \norm{F(s)}_{L^\infty(\R^d)}ds+\dfrac{C^\ast}{\nu^{1/3}}\left[\int_{\sigma(t)}^t\norm{F(s)}^3_{L^\infty(\R^d)}ds \right]^{1/3}.
\end{gather}
It only remains to estimate the norm of the effective flux $F$ in terms of the functionals $\mathcal{A}_1$ and $\mathcal{A}_2$ and to do so, we distinguish two cases according to the dimension.

\noindent\textbf{Case $d=2$.}
We recall that the effective flux is given by 
\begin{gather}\label{ep1.52}
F=-(-\Delta)^{-1}\dvg (\rho \dot u)=\nu \dvg u- G,\quad G=P(\rho)-\widetilde P.
\end{gather}
Hence, interpolation inequality  yields the following estimate on the small interval $[0,\sigma(t)]$
\begin{align*}
\int_0^{\sigma(t)}\norm{F}_{L^\infty(\R^2)}
\leqslant C\int_0^{\sigma(t)}\norm{F}_{L^2(\R^2)}^{\frac12}\|\nabla F\|_{L^4(\R^2)}^{\frac23
}&\leqslant \int_0^{\sigma(t)}\|(G,\nu\dvg u)\|_{L^2(\R^2)}^{\frac13}\norm{\rho\dot u}_{L^4(\R^2)}^{\tfrac{2}{3}} \\
  &\leqslant \bigl(C^\ast E_0^{\tfrac{1}{6}}+\nu^{\tfrac{1}{6}}\mathcal{A}_1(t)^{\tfrac{1}{6}}\bigr))\int_0^{\sigma(t)}\norm{\rho\dot u}_{L^4(\R^2)}^{\tfrac{2}{3}}.
\end{align*}
Since the density contains  vacuum states, we are not allowed to bound the last factor in the above inequality solely by the Gagliardo-Nirenberg inequality. On $\mathbb{T}^2$, we have at our disposal a logarithmic interpolation inequality (see \cite{ danchin2023compressible, danchin2023exponential,desjardins1997regularity}) which is not valid in the whole space. To address this issue, we prove in \cref{leminter} an interpolation inequality that will allow us to take into account the vacuum state. Thus, from \cref{leminter}, we have:
\begin{align*}
    \int_0^{\sigma(t)}\norm{\rho\dot u}_{L^4(\R^2)}^{\tfrac{2}{3}}&\leqslant C^\ast \int_0^{\sigma(t)} \norm{\sqrt{\rho}\dot u}_{L^2(\R^2)}^{\tfrac{1}{3}}\norm{\nabla\dot u}^{\tfrac{1}{3}}_{L^2(\R^2)}+ C^\ast E_0^{\tfrac{1}{18}}\int_0^{\sigma(t)} \norm{\sqrt{\rho}\dot u}_{L^2(\R^2)}^{\tfrac{2}{9}} \norm{\nabla\dot u}_{L^2(\R^2)}^{\tfrac{4}{9}}\\
    &\leqslant C^\ast \mathcal{A}_1(t)^{\tfrac{1}{6}}\mathcal{A}_2(t)^{\tfrac{1}{6}}+ C^\ast E_0^{\tfrac{1}{18}}\mathcal{A}_1(t)^{\tfrac{1}{9}}\mathcal{A}_2(t)^{\tfrac{2}{9}}
\end{align*}
and therefore:
\begin{gather}\label{ep1.50}
    \int_0^{\sigma(t)}\norm{F}_{L^\infty(\R^2)}\leqslant C^\ast\bigl(1+E_{0}^{\tfrac{1}{18}}\bigr))\bigl( E_0^{\tfrac{1}{6}}+\nu^{\tfrac{1}{6}}\mathcal{A}_1(t)^{\tfrac{1}{6}}\bigr))\bigl(\mathcal{A}_1(t)^{\tfrac{1}{3}}+\mathcal{A}_2(t)^{\tfrac{1}{3}}\bigr)).
\end{gather}
Similarly, interpolation inequality and  \cref{leminter} yield  on the larger time interval $[\sigma(t), t]$
\begin{align*}
&\int_{\sigma(t)}^t\norm{F}_{L^\infty(\R^2)}^3 \leqslant \int_{\sigma(t)}^t\|(G,\nu\dvg u)\|_{L^2(\R^2)}\norm{\rho\dot u}_{L^4(\R^2)}^{2} 
\\
&  \leqslant \bigl(C^\ast E_0^{\tfrac{1}{2}}+\nu^{\tfrac{1}{2}}\mathcal{A}_1(t)^{\tfrac{1}{2}}\bigr)
\bigl(C^\ast \int_0^t \norm{\sqrt{\rho}\dot u}_{L^2(\R^2)}\norm{\nabla\dot u}_{L^2(\R^2)}+ C^\ast E_0^{\tfrac{1}{6}}\int_0^t \norm{\sqrt{\rho}\dot u}_{L^2(\R^2)}^{\tfrac{2}{3}} \norm{\nabla\dot u}_{L^2(\R^2)}^{\tfrac{4}{3}}\bigr)
\\
&\leqslant \bigl(C^\ast E_0^{\tfrac{1}{2}}+\nu^{\tfrac{1}{2}}\mathcal{A}_1(t)^{\tfrac{1}{2}}\bigr)
\bigl(C^\ast \mathcal{A}_1(t)^{\tfrac{1}{2}}\mathcal{A}_2(t)^{\tfrac{1}{2}}+ C^\ast E_0^{\tfrac{1}{6}}\mathcal{A}_1(t)^{\tfrac{1}{3}}\mathcal{A}_2(t)^{\tfrac{2}{3}}\bigr).
\end{align*}
Hence, 
\begin{gather}\label{ep1.51}
\int_{\sigma(t)}^t\norm{F}_{L^\infty(\R^2)}^3\leqslant C^\ast\bigl(1+E_0^{\tfrac{1}{6}}\bigr)\bigl(C^\ast E_0^{\tfrac{1}{2}}+\nu^{\tfrac{1}{2}}\mathcal{A}_1(t)^{\tfrac{1}{2}}\bigr)\bigl(\mathcal{A}_1(t)+\mathcal{A}_2(t)\bigr).
\end{gather}
Finally \eqref{ep1.49}, \eqref{ep1.50} and \eqref{ep1.51}, lead to:
\begin{multline}
    \sup_{[0,t]}\norm{\rho-\widetilde\rho}_{L^\infty(\R^2)} \leqslant \norm{\rho_0-\widetilde\rho}_{L^\infty(\R^2)}+ \dfrac{C^\ast}{\nu^{\frac13}}  \bigl(1+E_{0}^{\tfrac{1}{18}}\bigr)\bigl( E_0^{\tfrac{1}{6}}+\nu^{\tfrac{1}{6}}\mathcal{A}_1(t)^{\tfrac{1}{6}}\bigr)  \bigl(\mathcal{A}_1(t)^{\tfrac{1}{3}}+\mathcal{A}_2(t)^{\tfrac{1}{3}}\bigr).
\end{multline}
\noindent\textbf{Case $d=3$.}
From the expression of the effective flux \eqref{ep1.52}, we have by Gagliardo-Nirenberg inequality:
\begin{align}
   & \int_0^{\sigma(t)} \norm{F}_{L^\infty(\R^3)}
   \leqslant C\int^{\sigma(t)}_0 \|F\|_{L^6(\R^3)}^{\frac12}\|\nabla F\|_{L^6(\R^3)}^{\frac12}\leqslant  C\int_0^{\sigma(t)} \norm{\rho\dot u}_{L^2(\R^3)}^{\tfrac{1}{2}}\norm{\rho\dot u}_{L^6(\R^3)}^{\tfrac{1}{2}}\leqslant C^\ast\mathcal{A}_1(t)^{\tfrac{1}{4}}\mathcal{A}_2(t)^{\tfrac{1}{4}},\label{ep3.25}
   \\
   & \int_{\sigma(t)}^t \norm{F}_{L^\infty(\R^3)}^3\leqslant C\int_{\sigma(t)}^t \norm{\rho\dot u}_{L^2(\R^3)}^{\tfrac{3}{2}}\norm{\rho\dot u}_{L^6(\R^3)}^{\tfrac{3}{2}}\leqslant C^\ast \mathcal{A}_1(t)^{\tfrac{1}{4}}\mathcal{A}_2(t)^{\tfrac{5}{4}}.\nonumber
\end{align} 
Finally, we get:
\begin{gather*}
    \sup_{[0,t]}\norm{\rho-\widetilde\rho}_{L^\infty(\R^3)} \leqslant \norm{\rho_0-\widetilde\rho}_{L^\infty(\R^3)}+  \dfrac{ C^\ast}{\nu^{1/3}} \bigl(\mathcal{A}_1(t)^{\tfrac{1}{2}}+\mathcal{A}_2(t)^{\tfrac{1}{2}}\bigr),
\end{gather*}
and this ends the proof of \cref{Aprioriestimates3}.
\enddem
\subsubsection{Proof of \cref{Aprioriestimates4}}\label[section]{Proofs4}
The present section is devoted to the propagation of tangential regularity of the density  
together with the derivation of the Lipschitz norm of the velocity, provided with the (time-independent) bounds $C_0$ for $\mathcal{A}_1(t),\mathcal{A}_2(t), \rho(t)$.

 We recall that the family of vector fields $\mathcal{X}(t)=(X_{\upsilon}(t))_{1\leqslant \upsilon\leqslant m}$  is defined as solution of \eqref{ep4.33}
\begin{gather}\label{ep1.53}
\begin{cases}
    \dpt X_\upsilon+ u\cdot \nabla X_\upsilon=\partial_{X_\upsilon}u=(X_\upsilon\cdot \nabla)u,\\
    {X_{\upsilon}}_{|t=0}=X_{0,\upsilon},
\end{cases}
\end{gather}
and we can estimate the norms of
\begin{align*}
    \|\mathcal{X}(t)\|_{\mathbb{L}^{\infty,p}(\R^d)}=\sup_{1\leqslant\upsilon\leqslant m}\|X_\upsilon(t)\|_{\mathbb{L}^{\infty,p}(\R^d)}
    =\sup_{1\leqslant \upsilon\leqslant m}
    \bigl( \|X_\upsilon(t)\|_{L^\infty(\R^d)}+\|\nabla X_\upsilon(t)\|_{L^p(\R^d)}\bigr) 
\end{align*} 
by use of the Lipschitz norm of the velocity.
On the other side, these norms will help to bound the Lipschitz norm of the velocity, and in particular the pressure-related part $\nabla u_G$ in the decomposition   \eqref{ep4.36} of $\nabla u$:
\begin{align}\label{nablauG}
\nabla u&= \nabla \widetilde u+\nabla u_G \\
&:=\Bigl( -\dfrac{1}{\nu} \mathcal{R} \mathcal{R} (-\Delta)^{-1}\dvg (\rho \dot u)-\dfrac{1}{\mu}\mathcal{R}\mathcal{R}(-\Delta)^{-1}\cdot\rot (\rho \dot u)\Bigr)
+\Bigl(\dfrac{1}{\nu} \mathcal{R} \mathcal{R}G\Bigr).\nonumber
\end{align} 
Finally we will close the estimates for 
$$\mathcal{A}_3(t)=\|\mathcal{X}(t)\|_{\mathbb{L}^{\infty,p}(\R^d)}+\sup_{1\leqslant \upsilon\leqslant m}\|\dvg(\rho X_\upsilon)(t)\|_{L^p(\R^d)}$$
by Gronwall's inequality.
\dem

\noindent\textbf{Step 1. Preliminary estimates for $\|\mathcal{X}\|_{\mathbb{L}^{\infty,p}(\R^d)}$.} 
From \eqref{ep1.53} we deduce easily:
\begin{gather}\label{ep1.57}
\norm{X_\upsilon(t)}_{L^\infty(\R^d)}\leqslant \norm{{X_{0,\upsilon}}}_{L^\infty}+ \int_0^t \norm{X_\upsilon(s)}_{L^\infty(\R^d)}\norm{ \nabla u(s)}_{L^\infty(\R^d)}ds.
\end{gather}

We now take derivatives in \eqref{ep1.53} and we obtain:
\begin{gather}\label{ep1.54}
\dpt \partial_k X_\upsilon^j + (u\cdot \nabla)\partial_k X_\upsilon^j=\partial_k X_\upsilon\cdot \nabla u^j- \partial_k u\cdot \nabla X_\upsilon^j+ \partial_{X_\upsilon} \partial_k u^j.
\end{gather}
We take the trace in the above equality and we make use of the expression of the divergence of the velocity
$
\dvg u=\dfrac{1}{\nu}\left(F+ G\right), 
$
 in order to obtain the following equation for $\dvg X_\upsilon$:
\[
\dpt (\dvg X_\upsilon)+ u\cdot \nabla \dvg X_\upsilon=\dfrac{1}{\nu}\partial_{X_\upsilon}G+\dfrac{1}{\nu}\partial_{X_\upsilon} F.
\] 
Hence, it is straightforward to show the following:
\begin{multline}\label{ep1.55}
\norm{\dvg X_\upsilon(t)}_{L^p(\R^d)}\leqslant \norm{\dvg X_{0,\upsilon}}_{L^p(\R^d)}+ \dfrac{1}{\nu}\int_0^t \norm{\partial_{X_\upsilon} G}_{L^p (\R^d)}+ \dfrac{1}{\nu}\int_0^t \norm{X_\upsilon}_{L^\infty(\R^d)}\norm{\nabla F}_{L^p (\R^d)}\\
+\int_0^t \norm{\dvg u(s)}_{L^\infty(\R^d)}\norm{\dvg  X_\upsilon(s)}_{L^p(\R^d)}ds.
\end{multline}
As well, we take the antisymmetric part in \eqref{ep1.54} to get the following equation for $\rot X_\upsilon$:
\begin{multline*}
\dpt (\partial_k X^j_\upsilon-\partial_j X^k_\upsilon)+ u\cdot \nabla (\partial_k X^j_\upsilon-\partial_j X^k_\upsilon)= \partial_k X_\upsilon\cdot \nabla u^j-\partial_j X_\upsilon\cdot \nabla u^k\\
+ \partial_j u\cdot \nabla X_\upsilon^k- \partial_k u\cdot \nabla X_\upsilon^j+ \partial_{X_\upsilon} \left(\partial_k u^j-  \partial_j u^k\right),
\end{multline*}
from which we deduce easily: 
\begin{multline}\label{ep1.56}
    \norm{\rot X_\upsilon(t)}_{L^p(\R^d)}\leqslant \norm{\rot X_{0,\upsilon}}_{L^p(\R^d)}+ \int_0^t \norm{\nabla X_\upsilon}_{L^p(\R^d)}\norm{\nabla u}_{L^\infty(\R^d)}\\
    +\int_0^t\norm{X_\upsilon}_{L^\infty(\R^d)}\norm{\nabla\rot u}_{L^p(\R^d)}.
\end{multline}
By combining \eqref{ep1.55} and \eqref{ep1.56}, we obtain the following estimate for the vector field gradient:
\begin{align}
    \norm{\nabla X_\upsilon(t)}_{L^p(\R^d)}&\leqslant \norm{\nabla X_{0,\upsilon}}_{L^p(\R^d)}+ \dfrac{1}{\nu}\int_0^t \norm{\partial_{X_\upsilon} G}_{L^p(\R^d)}+ \int_0^t \norm{\nabla X_\upsilon(s)}_{L^p(\R^d)} \norm{\nabla u(s)}_{L^\infty(\R^d)}ds
    \nonumber\\
    &+\int_0^t\norm{X_\upsilon(s)}_{L^\infty(\R^d)}
    \Bigl\| \Bigl(\frac{1}{\nu}\nabla F(s), \nabla\rot u(s)\Bigr)\Bigr\|_{L^p(\R^d)} ds.\label{ep1.58}
\end{align}

\noindent\textbf{Step 2. Estimates for $\mathcal{A}_3(t)=  \norm{\mathcal{X}(t)}_{\mathbb L^{\infty,p}(\R^d)}+\sup_{1\leqslant \upsilon\leqslant m}\norm{\dvg (\rho X_\upsilon)(t)}_{L^p(\R^d)}$.}
In order to estimate the $L^p(\R^d)$ norm of $\dvg (\rho X_\upsilon)$, we combine the equation on the density $\eqref{ep0.1}_1$  and the equation on the vector field $X_\upsilon$ \eqref{ep1.53} in order to obtain 
\begin{gather*}
    \dpt (\dvg (\rho X_\upsilon))+\dvg (u \dvg (\rho X_\upsilon))=0,
\end{gather*}
from which we deduce the following estimates:
\begin{gather}\label{ep1.59}
    \norm{\dvg (\rho X_\upsilon)(t)}_{L^p(\R^d)}\leqslant \norm{\dvg (\rho_0 X_{0,\upsilon})}_{L^p(\R^d)}+ \int_0^t \norm{\dvg u(s)}_{L^\infty(\R^d)}\norm{\dvg (\rho X_\upsilon)(s)}_{L^p(\R^d)}ds.
\end{gather} 

Consequently we can estimate
\[
\partial_{X_\upsilon} G= P'(\rho) \dvg (\rho X_\upsilon)- \rho P'(\rho)\dvg X_\upsilon
\]
by
\[
\norm{\partial_{X_\upsilon} G}_{L^p(\R^d)}\leqslant C_0\left(\norm{\dvg (\rho X_\upsilon)}_{L^p(\R^d)}+\norm{\dvg X_\upsilon }_{L^p(\R^d)}\right).
\]
We combine  \eqref{ep1.57}, \eqref{ep1.58} and \eqref{ep1.59} together with \eqref{Lp:dotu} to get
\begin{gather*}
    \mathcal{A}_3(t)\leqslant   \mathcal{A}_3(0)+C_0 \int_0^t \mathcal{A}_3(s)\left[\dfrac{1}{\nu}+\norm{\nabla u(s)}_{L^\infty(\R^d)}+\norm{\rho \dot u(s)}_{L^p(\R^d)}\right]ds.
\end{gather*} 
Gronwall's Lemma yields:
\begin{gather}\label{ep1.60}
\mathcal{A}_3(t)\leqslant  \mathcal{A}_3(0) \exp\Bigl(
 C_0\int_0^t\left[1+\norm{\nabla u(s)}_{L^\infty(\R^d)}+ \norm{\rho \dot u(s)}_{L^p(\R^d)}\right]ds\Bigr).
\end{gather}

\noindent\textbf{Step 3. Estimates for $\|\nabla u_G\|_{L^1_tL^\infty_x}$.}

We consider $\nabla u_G$ first.
We use \cref{Linftyestimate} to bound the second part $\nabla u_G=\frac1{\nu}\mathcal{R}\mathcal{R}G$: 
\begin{equation*}\label{ep1.62}
\norm{\nabla u_G(t)}_{L^\infty(\R^d)}\leqslant  \dfrac{C}{\nu}\norm{G(t)}_{L^2(\R^d)} 
+\dfrac{C}{\nu}\norm{G(t)}_{L^\infty(\R^d)}\left[1+\log \left(e+\dfrac{\norm{G(t)}_{\mathbb{L}_{\mathcal X(t)}^p(\R^d)}}{\norm{G(t)}_{L^\infty(\R^d)}}\right)\right].
\end{equation*}
Now we focus on the estimate for
\begin{align*}
    \norm{G(t)}_{\mathbb{L}_{ \mathcal X(t)}^p(\R^d)}
    &\sim  \dfrac{1}{I( \mathcal X(t))}\left[\norm{G(t)}_{L^\infty(\R^d)}\norm{\mathcal{X}(t)}_{\mathbb L^{\infty,p}(\R^d)}+\sup_{1\leqslant \upsilon\leqslant m}\norm{\partial_{X_\upsilon}G(t)}_{L^p(\R^d)}\right] 
     \leqslant \dfrac{C_0}{I( \mathcal X(t))}\mathcal{A}_3(t).
\end{align*}
 From \eqref{ep1.53}, it is obvious that this quantity comes with the following lower bound:
\[
I(\mathcal X(t))\geqslant I(\mathcal X_0) \exp\left(-\int_0^t \norm{\nabla u(s)}_{L^\infty(\R^d)}ds\right).
\]
We hence have 
\begin{align}\label{ep4.24}
\norm{\nabla u_G(t)}_{L^\infty(\R^d)}
\leqslant   
C_0\Bigl( 1 +\frac{\mathcal{A}_3(0)}{I(\mathcal{X}_0)} + \int_0^t\left[ 1+\norm{\rho \dot u}_{L^p(\R^d)}+\norm{\nabla u}_{L^\infty(\R^d)}\right](s)ds\Bigr).
\end{align}
\noindent\textbf{Step 4. Final estimates.}
In the following we continue with the estimates for $\norm{\rho \dot u}_{L^1_tL^p(\R^d)}$ and $\norm{\nabla \widetilde u}_{L^1_tL^\infty(\R^d)}$, taking the dimension into account.
\paragraph{Case $d=2$}
Similar as in the proof of   \cref{Aprioriestimates3}, by using our interpolation inequality in \cref{leminter}, we obtain:
\begin{gather*}\label{ep3.29}
        \int_0^t\norm{\rho \dot u}_{L^p(\R^2)}\leqslant C_0 \int_0^t\Bigl[\norm{\sqrt{\rho} \dot u}_{L^2(\R^2)}^{\tfrac{2}{p}}\norm{\nabla \dot u}_{L^{2}(\R^2)}^{\tfrac{1}{p'}-\tfrac{1}{p}} +E_0^{\tfrac{1}{p}(1-\tfrac{p'}{2})}\norm{\sqrt{\rho} \dot u}_{L^2(\R^2)}^{\tfrac{p'}{p}}\norm{\nabla \dot u}_{L^{2}(\R^2)}^{1-\tfrac{p'}{p}}\Bigr],
\end{gather*}
where  the first integral of the right hand side above is bounded as
\begin{align*}
    &\int_0^t\norm{\sqrt{\rho} \dot u}_{L^2(\R^2)}^{\tfrac{2}{p}}\norm{\nabla \dot u}_{L^{2}(\R^2)}^{\tfrac{1}{p'}-\tfrac{1}{p}} =\int_0^t\norm{\sqrt{\rho} \dot u}_{L^2(\R^2)}^{\tfrac{2}{p}}\left[\sigma \norm{\nabla \dot u}_{L^{2}(\R^2)}^2\right]^{\tfrac{1}{2}-\tfrac{1}{p}}\sigma^{\tfrac{1}{p}-\tfrac{1}{2}}\\
    &\leqslant \left[\int_0^t\norm{\sqrt{\rho} \dot u}_{L^2(\R^2)}^{2}\right]^{\tfrac{1}{p}}\left[\int_0^t\sigma \norm{\nabla \dot u}_{L^{2}(\R^2)}^2\right]^{\tfrac{1}{2}-\tfrac{1}{p}}\left[\int_0^t\sigma^{\tfrac{2}{p}-1}\right]^{\tfrac{1}{2}}\leqslant C(p)(1+\sqrt{t})\mathcal{A}_1(t)^{\tfrac{1}{p}}\mathcal{A}_2(t)^{\tfrac{1}{2}-\tfrac{1}{p}},
\end{align*}
and the second integral of the right hand side above is bounded as
\begin{align*}
    &\int_0^t\norm{\sqrt{\rho} \dot u}_{L^2(\R^2)}^{\tfrac{p'}{p}}\norm{\nabla \dot u}_{L^{2}(\R^2)}^{1-\tfrac{p'}{p}} =\int_0^t\left[\norm{\sqrt{\rho} \dot u}_{L^2(\R^2)}^2\right]^{\tfrac{p'}{2p}}\left[\sigma\norm{\nabla \dot u}_{L^{2}(\R^2)}^2\right]^{\tfrac{1}{2}-\tfrac{p'}{2p}}\sigma^{\tfrac{p'}{2p}-\tfrac{1}{2}}\\
    &\leqslant \left[\int_0^t\norm{\sqrt{\rho} \dot u}_{L^2(\R^2)}^2\right]^{\tfrac{p'}{2p}}\left[\int_0^t\sigma\norm{\nabla \dot u}_{L^{2}(\R^2)}^2\right]^{\tfrac{1}{2}-\tfrac{p'}{2p}}\left[\int_0^t\sigma^{\tfrac{p'}{p}-1}\right]^{\tfrac{1}{2}}\leqslant C(p)(1+\sqrt{t})\mathcal{A}_1(t)^{\tfrac{p'}{2p}}\mathcal{A}_2(t)^{\tfrac{1}{2}-\tfrac{p'}{2p}}.
\end{align*}
In sum, for all $2<p<\infty$: 
\begin{gather}\label{ep4.26}
    \int_0^t \norm{\rho \dot u}_{L^p(\R^2)} \leqslant C_0(1+\sqrt{t}).
\end{gather}

Now following the computations  leading to \eqref{ep1.50}, it is straightforward to obtain the following 
\begin{align}\label{ep3.31} 
\int^t_0\|\nabla \widetilde u\|_{L^\infty(\R^2)} ds
&\leqslant \dfrac{1}{\nu}\int_0^t \norm{\mathcal{R} \mathcal{R} (-\Delta)^{-1}\dvg (\rho \dot u)}_{L^\infty(\R^2)}
+\frac{1}{\mu }\int_0^t\norm{\mathcal{R}\mathcal{R}(-\Delta)^{-1}\rot_{jk} (\rho \dot u)}_{L^\infty(\R^2)}
\\
&\leqslant 
C^\ast \bigl( \dfrac{ 1}{\nu^{5/6}}+1\bigr) (1+t^{\tfrac{2}{3}}) \left(1+E_0^{\tfrac{1}{18}}\right)\left(A_1(t)^{\tfrac{1}{2}}+\mathcal{A}_2(t)^{\tfrac{1}{2}}\right)
\leqslant C_0(1+t^{\frac23}). \nonumber
\end{align}
We plug \eqref{ep4.26} into \eqref{ep4.24}, sum \eqref{ep4.24} and \eqref{ep3.31} up, and finally use Gronwall's inequality to get the estimate $\eqref{ep4.25}_2$ for $d=2$.
The estimate $\eqref{ep4.25}_1$ for $\mathcal{A}_3(t)$ follows correspondingly from \eqref{ep1.60}.
\smallbreak
\paragraph{Case $d=3$} 
Similarly as the proof of \eqref{ep4.26}, for $3<p<6$ we interpolate  the $L^p(\R^3)$ norm of $\sqrt{\rho}\dot u$ between $L^2(\R^3)$ and $L^6(\R^3)$, next we make use of the embedding $\dot H^1(\R^3)\hookrightarrow L^6(\R^3)$ to derive 
\begin{align}
    \int_0^t\norm{\rho \dot u}_{L^p(\R^3)}&\leqslant C_0\int_0^t \norm{\sqrt{\rho}\dot u}_{L^2(\R^3)}^{\tfrac{3}{p}-\tfrac{1}{2}}\left[\sigma\norm{\nabla \dot u}_{L^2(\R^3)}^2\right]^{\tfrac{3}{4}-\tfrac{3}{2p}}\sigma^{\tfrac{3}{2p}-\tfrac{3}{4}}\leqslant C_0(1+\sqrt t) 
    .\label{ep4.27}
\end{align}
Owing to interpolation inequality and Sobolev embedding, we have
\begin{gather}\label{ep3.32}
\int_0^t \norm{\nabla \widetilde u}_{L^\infty(\R^3)}\leqslant \int_0^{t} \norm{\rho\dot u}_{L^2(\R^3)}^{\tfrac{1}{2}}\norm{\rho\dot u}_{L^6(\R^3)}^{\tfrac{1}{2}}\leqslant C_0(1+\sqrt t).
\end{gather}
Similar as above, we plug \eqref{ep4.27} into \eqref{ep4.24}, sum \eqref{ep4.24} and \eqref{ep3.32} up, and finally use Gronwall's inequality to get the estimate $\eqref{ep4.25}_2$ for $d=3$.
The estimate $\eqref{ep4.25}_1$ for $\mathcal{A}_3(t)$ follows   from \eqref{ep1.60} in dimension three.
This ends the proof of \cref{Aprioriestimates4}.
\enddem
\subsection{\texorpdfstring{Proof of \cref{th1} and \cref{theo1}}{}} \label[section]{approximate}

    This section is devoted to the final step in  the proof of  the main result, \cref{th1}. 
We recall that   we are considering the Cauchy problem associated with equations  \eqref{ep0.1} and  with initial data 
\eqref{ep4.30} verifying \eqref{ep3.16} and \eqref{ep4.29}.

Usually the sequence of initial data $(\rho_0^\delta,u_0^\delta)$ is   obtained by mollifying 
$(\rho_0,u_0)$ with a smooth kernel.
This regularization procedure has the unfortunate effect of destroying the density's structure. As observed in \cite{hoff2002dynamics,zodji2023discontinuous}, the most effective approach is to construct the approximate solutions in a class that is very close to the limit.  From this point of view, the local result obtained by Danchin, Fanelli, Paicu in \cite{danchin2020well} should be appropriate. 
However, the argument of the maximum regularity of the heat equations requires   the density to be a small perturbation of a constant state, even for the local solution. 
We are therefore led to prove the local well-posedness  of equations \eqref{ep0.1} stated in   \cref{intro:local} in \cref{local}.
\dem[Proof of \cref{th1}] 
In order to apply the local-in-time well-posedness results in \cref{intro:local}, we consider a sequence of initial densities $(\rho_0^\delta)_\delta$ verifying: for all $0<\delta<1$:
\begin{gather}
        \rho_0^\delta \geqslant \delta,\quad \rho_0^\delta-\widetilde\rho \in L^2(\R^d)\cap L^\infty(\R^d)\cap \mathbb L_{\mathcal{X}_0}^p(\R^d)
\end{gather}
such that 
\[
\rho^\delta_0-\widetilde\rho \xrightarrow{\delta\to 0} \rho_0-\widetilde\rho \quad \text{ in } L^2(\R^d).
\]
The construction of the sequence of initial velocities $(u_0^\delta)_\delta\subset H^1(\R^d)$, converging strongly to $u_0$ in $H^1(\R^d)$ and satisfying the compatibility condition 
\begin{gather}
    \dvg \{2\mu \D u_0^\delta+ (\lambda \dvg u_0^\delta- P(\rho_0^\delta)+\widetilde P)\}\in L^2(\R^d)
\end{gather}
can be found in \cite[Section 3.5]{zodji2023discontinuous}. 
Now we can apply \cref{intro:local} to get the existence of a unique solution $(\rho^\delta, u^\delta)$  that satisfies: 
\begin{gather}
    \begin{cases}
        \dpt \rho^\delta+\dvg (\rho^\delta u^\delta)=0,\\
        \dpt(\rho^\delta u^\delta)+\dvg (\rho^\delta u^\delta \otimes u^\delta)+ \nabla P(\rho^\delta)=\mu \Delta u^\delta+ (\mu+\lambda)\nabla \dvg u^\delta
    \end{cases}
\end{gather}
with initial data 
\begin{gather*}
    (\rho^\delta)_{|t=0}=\rho^\delta_0\quad \text{ and }\quad (u^\delta)_{|t=0} =u^\delta_0.
\end{gather*}
The solution is defined up to a maximal time $T^\delta$ and enjoys the regularities outlined in \cref{intro:local} which are sufficient for the computations performed in the preceding sections to be meaningful, leading to  \cref{lemma4} and 
\cref{Aprioriestimates4}. In particular, all the conditions outlined in the blow-up criterion  \eqref{intro:thblowup} are satisfied, implying that $T^\delta=+\infty$. Finally, employing classical arguments involving compact embedding, Aubin-Lions Lemma and leveraging the regularity of the effective flux, one can establish the strong convergence of a subsequence of $(\rho^\delta, u^\delta)$ to  $(\rho,u)$ satisfying  the regularity in \cref{th1}. Furthermore, given the   improved-in-time Lipschitz bound of the velocity field   in \cref{rema1},
a change of variables into Lagrangian coordinates ensures the uniqueness of such a solution. We refer for example to \cite{danchin2020well} for the computations.

\enddem

\dem[Proof of \cref{theo1}]
At this level, we obtain a sequence $(\rho^{(\nu)}, u^{(\nu)})_{\nu\geqslant \underline{\nu}}$   satisfying:
\begin{gather}\label{ep4.16}
    \begin{cases}
        \dpt \rho^{(\nu)}+\dvg (\rho^{(\nu)} u^{(\nu)})=0,\\
        \dpt (\rho^{(\nu)} u^{(\nu)} )+\dvg (\rho^{(\nu)} u^{(\nu)} \otimes u^{(\nu)})-\nabla F^{(\nu)}- \mu \Delta u^{(\nu)}=-\dfrac{\mu}{\nu}\nabla F^{(\nu)} -\dfrac{\mu}{\nu}\nabla P(\rho^{(\nu)}),
    \end{cases}
\end{gather}
with initial data \eqref{ep4.30} verifying \eqref{ep3.16} and \eqref{ep4.29} and $\dvg u_0=0$. Above the effective flux $F^{(\nu)}$ solves the following elliptic equation:
\[
\Delta F^{(\nu)}= \dvg (\rho^{(\nu)} \dot u^{(\nu)}).
\]
Given that $(\rho^{(\nu)})_\nu$ is bounded in $L^\infty((0,\infty)\times \R^d)$ and $(\sqrt{\rho^{(\nu)}} \dot u^{(\nu)})$ is bounded in $L^2((0,\infty)\times \R^d)$ it follows that $(-\nabla F^{(\nu)})_\nu$ is bounded in $L^2((0,\infty)\times \R^d)$ , resulting in weak convergence  to some $\nabla \Pi\in L^2((0,\infty)\times \R^d)$. Obviously, the right hand side of $\eqref{ep4.16}_2$ converges strongly to zero in $L^\infty((0,\infty), \dot H^{-1}(\R^d))$, given that $(\nu^{-1/2} F^{(\nu)})$ is bounded in $L^\infty((0,\infty), L^2( \R^d))$ and  $(P(\rho^{(\nu)})-\widetilde P)_\nu$ is bounded in $L^\infty((0,\infty), L^2(\R^d)))$.

The regularity of the sequence $(u^{(\nu)})$ ensures that, up to a subsequence, $(u^{(\nu)})$ converges strongly in  $L^2_\loc ((0,\infty)\times \R^d)$ to some $v\in L^\infty((0,T), H^1(\R^d))$. Furthermore, since the sequence $(\rho^{(\nu)})_\nu$ is bounded in $L^\infty((0,\infty)\times \R^d)$, it converges weakly-* to some $\varrho\in L^\infty((0,\infty)\times \R^d)$. Additionally, $(\dvg u^{(\nu)})_\nu$ converges strongly to zero in $L^\infty((0,\infty), L^2( \R^d))$, since from the bound of functional $\mathcal{ A}_1$, the sequence $\left(\nu \norm{\dvg u^{(\nu)}}_{L^2(\R^d)}^2\right)_\nu$ is bounded. These convergences, along with the Aubin-Lions Lemma, are sufficient to pass to the limit in \eqref{ep4.16} and establish that $(\varrho,v)$ solves the incompressible model \eqref{ep4.35}. 
    The uniqueness result for \eqref{ep4.35} in \cite{prange2023inhomogeneous} and the uniform bounds in \cref{th1} implies the convergence of the whole sequence $(\rho^{(\nu)}, u^{(\nu)})_\nu$ (instead of some subsequence).
This completes the proof of   \cref{theo1}.
\enddem
\section*{Acknowledgment}
X. Liao is partially supported by Deutsche Forschungsgemeinschaft (DFG, German Research Foundation)-Project-ID 258734477-SFB 1173.
S. M. Zodji is funded by the European Union’s Horizon 2020 research and innovation program
under the Marie Skłodowska-Curie grant agreement No 945332.  This work is partially
supported by the project CRISIS (ANR-20-CE40-0020-01), operated by the French National Research Agency
(ANR). 

This work is partially finished during S.M. Zodji's research visit at the Karlsruhe Institute of Technology (KIT), and he acknowledges the hospitality of KIT.
The authors would like to thank   Rapha\"el Danchin, Francesco Fanelli and Cosmin Burtea and David Gérard-Varet (the two PhD supervisors of S.M. Zodji) for their fruitful
discussions and valuable comments.  
\appendix
\section{Interpolation inequality}\label{sec:inter}
\begin{lemma}[Density-weighted interpolation inequality]\label[lemma]{leminter}
    Let $v\in\dot H^1(\R^2)$, $\rho\geqslant 0$ such that $\sqrt{\rho} v\in L^2(\R^2)$ and $\rho-\widetilde \rho\in L^p(\R^2)$ for some $1<p<\infty$, with $\widetilde \rho>0$. Then $v\in L^2(\R^2)$ and there exists a constant $C>0$ depending only on $\widetilde\rho$ and $p$ such that the following estimate holds true: 
    \begin{gather}\label{ep1.47}
        \norm{v}_{L^2(\R^2)} \leqslant C\left(\norm{\rho-\widetilde\rho}_{L^p(\R^2)}^{\tfrac{p}{2}}\norm{\nabla v}_{L^2(\R^2)}+ \norm{\sqrt{\rho} v}_{L^2(\R^2)}\right).
    \end{gather}
    Moreover,  for all $2<q<\infty$, we have:
    \begin{gather}\label{ep1.48}
        \norm{\rho^{\tfrac{q'}{2q}} v}_{L^q(\R^2)}\leqslant C \left(\norm{\sqrt{\rho} v}_{L^2(\R^2)}^{\tfrac{2}{q}}\norm{\nabla v}_{L^{2}(\R^2)}^{\tfrac{1}{q'}-\tfrac{1}{q}} +\norm{\rho-\widetilde\rho}_{L^p(\R^2)}^{\tfrac{p}{q}(1-\tfrac{q'}{2})}\norm{\sqrt{\rho} v}_{L^2(\R^2)}^{\tfrac{q'}{q}}\norm{\nabla v}_{L^{2}(\R^2)}^{1-\tfrac{q'}{q}}\right).
    \end{gather}
\end{lemma}
\dem
In the first step of the proof, we emulate the approach in \cite[Proposition A.6]{prange2023inhomogeneous} by expressing:
\begin{gather}
\widetilde\rho \abs{v}^2= (\widetilde\rho -\rho)\abs{v}^2+ \rho \abs{v}^2.
\end{gather}
 Due to the assumption that $\sqrt{\rho} v\in L^2(\R^2)$, we only need to compute the integral of the first term of the right hand above. With the help of interpolation, H\"older and Young inequalities, we have:
\begin{align*}
\int_{\R^2} (\widetilde\rho -\rho)\abs{v}^2&\leqslant \norm{\rho-\widetilde\rho}_{L^p(\R^2)}\norm{v}_{L^{2p'}(\R^2)}^2\\
                                           &\leqslant C_p \norm{\rho-\widetilde\rho}_{L^p(\R^2)}\norm{v}_{L^{2}(\R^2)}^{\tfrac{2}{p'}}\norm{\nabla v}_{L^2(\R^2)}^{\tfrac{2}{p}}\\
                                           &\leqslant \dfrac{1}{2}\widetilde\rho\norm{v}_{L^2(\R^2)}^2+ C_{p,\widetilde\rho}\norm{\rho-\widetilde\rho}_{L^p(\R^2)}^p\norm{\nabla v}_{L^2(\R^2)}^2,
\end{align*}
and  \eqref{ep1.47} just follows. Next, for all $2<q<\infty$,  H\"older, Gagliardo Nirenberg inequalities yield:
\begin{align*}
\norm{\rho^{\tfrac{q'}{2q}} v}_{L^q(\R^2)}&\leqslant \norm{\sqrt{\rho} v}_{L^2(\R^2)}^{\tfrac{q'}{q}}\norm{v}_{L^{2q}(\R^2)}^{1-\tfrac{q'}{q}}\\
&\leqslant C\norm{\sqrt{\rho} v}_{L^2(\R^2)}^{\tfrac{q'}{q}}\norm{\nabla v}_{L^{2}(\R^2)}^{\tfrac{1}{q'}-\tfrac{1}{q}}\norm{v}_{L^2(\R^2)}^{\tfrac{1}{q}-\tfrac{q'}{q^2}}\\
&\leqslant C\norm{\sqrt{\rho} v}_{L^2(\R^2)}^{\tfrac{q'}{q}}\norm{\nabla v}_{L^{2}(\R^2)}^{\tfrac{1}{q'}-\tfrac{1}{q}}\norm{v}_{L^2(\R^2)}^{\tfrac{2}{q}-\tfrac{q'}{q}}.
\end{align*}
Hence \eqref{ep1.48} holds true while replacing the $L^2(\R^2)$ norm of the velocity by \eqref{ep1.47}.
\enddem

\section{Local well-posedness}\label[appendix]{local}
In this section, we prove the local well-posedness result in \cref{intro:local} of the Navier-Stokes equations for a compressible fluid with an initial density having  tangential regularity. Our method relies on a change of variables into Lagrangian coordinates, followed by the study of the linearized system and  the full nonlinear system, in a similar way as in \cite{zodji2023well}.
In particular, we do not require the density to be a small perturbation around a constant state in $L^\infty(\mathbb{R}^d)$.

More precisely, we consider the Cauchy problem of the system \eqref{ep0.1}:
\begin{gather}\label{ep4.7}
    \begin{cases}
        \dpt \rho +\dvg (\rho u)=0,\\
        \dpt (\rho u)+\dvg (\rho u\otimes u)+\nabla P(\rho)=\mu\Delta u+ (\mu+\lambda)\nabla \dvg u,
    \end{cases}
\end{gather}
equipped with initial data \eqref{ep4.30}:
\begin{gather}\label{ep4.8}
\rho_{|t=0}=\rho_0\quad \text{ and } \quad (\rho u)_{|t=0}=\rho_0 u_0,
\end{gather}
verifying \eqref{ep3.16}, \eqref{ep4.29} and \eqref{ep4.44}:
\begin{gather}\label{ep4.9}
    \begin{cases}\vspace{0.2cm}
        0<\underline \rho \leqslant \rho_0(x),\quad \rho_0-\widetilde \rho\in L^2(\R^d)\cap L^\infty(\R^d)\cap \mathbb L^p_{\mathcal X_0}(\R^d),\\
        u_0\in H^1(\R^d),\quad \mu \Delta u_0+ (\mu+\lambda)\nabla \dvg u_0- \nabla P(\rho_0)\in L^2(\R^d).
    \end{cases}
\end{gather}
Above $\widetilde\rho>0$ is a constant and $\mathcal X_0=(X_{0,\upsilon})_{1\leqslant \upsilon\leqslant m}\subset \mathbb L^{\infty,p}(\R^d)$, $d<p<\infty$, is  a non-degenerate family of $m\in \N^*$ vectors fields, with $m\geqslant d-1$. 

Here, we present the main lines of the proof of \cref{intro:local}; detailed computations can be found in \cite{zodji2023well}, in the more involved case of density-dependent viscosity.

\noindent\textbf{Step 1. Lagrangian coordinates.}
Let $0<T\leqslant \infty$,  $u$ be a Lipschitz vector field such that $\nabla u\in L^1((0,T), L^{\infty}(\R^d))$ and let us consider its flow map $\mathscr{X}$ given by 
\[
\mathscr{X}_{\overline{u}}(t,y)= y+ \int_0^t u(\tau,\mathscr{X}_{\overline{u}}(\tau,y))d\tau=:y+ \int_0^t \overline{u}(\tau,y)d\tau
\]
where, hereafter, for all $g=g(t,x)$, we define $\overline{g}=g(t,y)$ by:
\[
\overline{g}(t,y)= g(t,\mathscr{X}(t,y)).
\]
By performing this change of variables, the equations \eqref{ep4.7} reads:
\begin{gather}\label{ep4.10}
    \begin{cases}
        \dpt (\overline{\rho} J_{\overline{u}})=0,\\
        \rho_0 \dpt \overline{u}=\dvg \left( \adj(D \mathscr{X}_{\overline{u}})\{2\mu \D_{A_{\overline{u}}} \overline{u}+ (\lambda\dvg_{A_{\overline{u}}} \overline{u}- P(\overline{\rho})+\widetilde P)\} \right),
    \end{cases}
\end{gather}
where 
\[
J_{\overline{u}}=\det (D \mathscr{X}_{\overline{u}}),\quad A_{\overline{u}}=  (D \mathscr{X}_{\overline{u}})^{-1},\quad \dvg_{A_{\overline{u}}} w= A^T_{\overline{u}} \colon \nabla w= D w\colon A_{_{\overline{u}}},\quad 2\D_{A_{\overline{u}}} w=  D w \cdot A_{\overline{u}}+ A^T_{\overline{u}}\cdot \nabla w.
\]

\noindent\textbf{Step 2. Well-posedness of the linearized system.}
Motivated by \eqref{ep4.10}, we are led to consider the linear system
\begin{gather}\label{ep4.11}
\begin{cases}
    \rho_0 \dpt v-\mu \Delta v-(\mu+\lambda)\nabla \dvg v=\dvg f,\\
    v_{|t=0}=v_0,
\end{cases}
\end{gather}
where the source term $f$ and the initial datum $v_0$ belong to the following space:
\begin{align*}
Y_T:= &\left\{(f,v_0)\in  L^\infty((0,T), L^2(\R^d))\times H^1(\R^d)  \colon  f, \dpt f, \sigma \partial_{tt} f\in L^2((0,T)\times \R^d);\right.\\
&\left.\sqrt{\sigma}\partial_{t} f\in L^\infty((0,T), L^2(\R^d)),\;
\mu \Delta v_0+ (\mu+\lambda)\nabla \dvg v_0+ \dvg f_{|t=0} \in L^2(\R^d)\right\}.
\end{align*}
The solution of the linearized system \eqref{ep4.11} is constructed in the following space:
\begin{align*}
        Z_T:=\left\{v\in L^\infty((0,T), H^1(\R^d))\colon\right.&\left. \nabla v,\; \dpt v,\; \nabla \dpt v, \; \sqrt{\sigma}\partial_{tt} v,\; \sigma \nabla \partial_{tt}v\in L^2((0,T)\times \R^d);\;\;
        \right.\\
        &\left.
        \dpt v,\; \sqrt{\sigma}\nabla \dpt v,\; \sigma \partial_{tt} v\in L^\infty((0,T), L^2(\R^d)) 
 \right\}.
\end{align*}
It is straightforward  to check that for $T<\infty$ every $v\in Z_T$ satisfies:
\[
v\in \cC([0,T], H^1(\R^d)),\quad \dpt v\in \cC((0,T], L^2(\R^d)).
\]
The well-posedness result for the linearized system \eqref{ep4.11} reads as follows.

\begin{prop}\label[prop]{prop1}
    Let $0<T\leqslant \infty$. For all $(f,v_0)\in Y_T$, there exists a unique solution  $v\in Z_T$ of the Cauchy problem  \eqref{ep4.11}. Moreover, the following estimates holds true for $v$.
    \begin{enumerate}
        \item Basic energy estimates:
        \begin{multline}\label{ep4.12}
            \sup_{[0,T]}\norm{\bigl(\sqrt{\rho_0}v,\; \sqrt{\rho_0}\dpt v,\; \nabla v\bigr)}_{L^2(\R^d)}^2+\int_0^T \norm{\bigl(\nabla v,\; \sqrt{\rho_0} \dpt v,\; \nabla \dpt v\bigr)}_{L^2(\R^d)}^2 \\
            \lesssim \norm{\bigl(\sqrt{\rho_0} v_0,\; \sqrt{\rho_0} \dpt v_{|t=0},\; \nabla v_0\bigr)}_{L^2(\R^d)}^2+ \sup_{[0,T]}\norm{f}_{L^2(\R^d)}^2+\int_0^T \norm{\bigl(f,\, \dpt f\bigr)}_{L^2(\R^d)}^2.
        \end{multline}
        \item  Higher energy estimates:
        \begin{multline}\label{ep4.13}
             \sup_{[0,T]}\norm{\bigl(\sqrt{\sigma}\nabla \dpt v,\; \sigma \sqrt{\rho_0} \partial_{tt} v\bigr)}_{L^2(\R^d)}^2+\int_0^T \norm{\bigl(\sqrt{\sigma\rho_0} \partial_{tt} v,\; \sigma\nabla (\partial_{tt} v)\bigr)}_{L^2(\R^d)}^2\\
             \lesssim \int_0^T \norm{\nabla \dpt v}_{L^2(\R^d)}^2+\sup_{[0,T]}\sigma\norm{\dpt f}_{L^2(\R^d)}^2+\int_0^T\sigma^2\norm{\partial_{tt}f}_{L^2(\R^d)}^2.
        \end{multline}
    \end{enumerate}
    The constant appearing in the above estimates does not depend on the upper or lower bound of the density $\rho_0$.
\end{prop} 
The proof of \cref{prop1} is not part of the classical theory of parabolic systems due to the roughness of the density. However, it can be achieved by a regularization process followed by a compactness argument. We refer, for example, to \cite[Theorem 3.1]{zodji2023well} for the derivation of estimates \eqref{ep4.12} and \eqref{ep4.13}. 

\noindent\textbf{Step 3. Further estimates of the linearized system.} 
By interpolating the estimates \eqref{ep4.12} and \eqref{ep4.13}, we observe that  the following  estimates hold true for the velocity gradient and its time derivative.
\begin{coro}\label[coro]{coro1}
    The following estimates hold true.
    \begin{enumerate}
        \item Assuming that $f\in L^r((0,T), L^p(\R^d))$ for $2\leqslant r\leqslant \infty$ and $2<p<\infty$ if $d=2$ or $2<p\leqslant 6$ if $d=3$,
          we have:
        \begin{gather}\label{ep4.20}
            \norm{\nabla v}_{L^r((0,T),L^p(\R^d))}^2\lesssim \norm{(f,v_0)}_{Y_T}^2+ \norm{f}_{L^r((0,T),L^p(\R^d))}^2.
        \end{gather}
        The same estimate holds also true  if $d=3$, for $6<p<\infty$, and $2\leqslant r\leqslant 4p/(p-6)$.
        \item \label{ep4.17} For all $2\leqslant r\leqslant \infty$ and $2<p<\infty$ if $d=2$ and $2<p\leqslant 6$ if $d=3$ we have:
         \begin{gather}\label{ep4.21}
             \norm{\sigma^{s}\nabla \dpt v}_{L^r((0,T),L^p(\R^d))}^2\lesssim \norm{(f,v_0)}_{Y_T}^2+ \norm{\sigma^{s} \dpt f}_{L^r((0,T),L^p(\R^d))}^2
         \end{gather}
         where 
         \begin{gather*}
         s=
             \begin{cases}\vspace{0.2cm}
                 1-\dfrac{1}{p}-\dfrac{1}{r} \quad \text{ if } \quad d=2,\\
                 \dfrac{5}{4}-\dfrac{3}{2p}-\dfrac{1}{r} \quad \text{ if } \quad d=3.
             \end{cases}
         \end{gather*}
         For $d=3$, the same estimate also holds true for all $6<p<\infty$ and $2\leqslant r\leqslant 4p/(p-6)$.
        \item 
        Let $\mathcal X_0=(X_{0,\upsilon})_{1\leqslant \upsilon\leqslant m}\subset \mathbb L^{\infty,p}(\R^d)$, $d<p<\infty$, be    a non-degenerate family of $m\in \N^*$ vectors fields, with $m\geqslant d-1$.

        \vspace{0.2cm}
        
        \begin{enumerate}
            \item Assuming that $f\in L^r((0,T), L^\infty(\R^d)\cap \mathbb L^p_{\mathcal X_0}(\R^d))$, with $2\leqslant r\leqslant 8$ if $d=2$ and $ 2\leqslant r\leqslant 32/9$ if $d=3$, then $ \nabla v\in L^r((0,T), L^\infty(\R^d))$ and the following estimate holds true.
         \begin{gather}\label{ep4.22}
            \norm{\nabla v}_{L^r((0,T), L^\infty(\R^d))}^2\lesssim \norm{(f,v_0)}_{Y_T}^2+\norm{f}_{L^r((0,T), L^2(\R^d)\cap L^\infty(\R^d)\cap \mathbb L^p_{\mathcal X_0}(\R^d))}^2.
        \end{gather}
            \item Assuming that $f\in L^r((0,T), L^\infty(\R^d)\cap \mathbb L^p_{\mathcal X_0}(\R^d))$, with $2<p<\infty $, $2\leqslant r\leqslant 2p/(p-2)$ if $d=2$ and $3<p\leqslant 6$, $2\leqslant r\leqslant 4p/(3p-6)$ if $d=3$.  
           Then we have:
        \begin{multline}\label{ep4.14}
            \norm{\partial_{\mathcal X_0} \nabla v}_{L^r((0,T), L^p(\R^d))}^2\lesssim \norm{\mathcal X_0}_{L^\infty(\R^d)}^2\norm{(f,v_0)}_{Y_T}^2
            +\norm{\partial_{\mathcal X_0} f}_{L^r((0,T),L^p(\R^d))}^2\\+\norm{\nabla\mathcal X_0}_{L^p(\R^d)}^2\norm{f}_{L^r((0,T), L^2(\R^d)\cap L^\infty(\R^d)\cap \mathbb L^p_{\mathcal X_0}(\R^d))}^2.
        \end{multline}
        Let $d=3$ and  $6<p<\infty$, $ 2\leqslant r\leqslant \infty$. If $\sigma^s f\in L^r((0,T), L^\infty(\R^3)\cap \mathbb L^p_{\mathcal X_0}(\R^3))$ and $\sigma^{\tfrac{3}{4}-\tfrac{1}{r}} \dpt f\in L^3(\R^3)$ with 
        \[
        s=\dfrac{3}{4}-\dfrac{1}{r}-\dfrac{3}{2p}
        \]
        then, from \eqref{ep4.17}, we have $\sigma^{\tfrac{3}{4}-\tfrac{1}{r}} \nabla \dpt v\in L^3(\R^3)$ and
        \begin{multline}\label{ep4.18}
            \norm{\sigma^s\partial_{\mathcal X_0} \nabla v}_{L^r((0,T), L^p(\R^3))}^2\lesssim \norm{\mathcal X_0}_{L^\infty(\R^3)}^2\left(\norm{(f,v_0)}_{Y_T}^2+\norm{\sigma^{\tfrac{3}{4}-\tfrac{1}{r}} \dpt f}_{L^r((0,T),L^3(\R^3))}^2\right)\\
            +\norm{\sigma^s\partial_{\mathcal X_0} f}_{L^r((0,T),L^p(\R^3))}^2+\norm{\nabla \mathcal X_0}_{L^p(\R^3)}^2\norm{\sigma^sf}_{L^r((0,T), L^2(\R^3)\cap L^\infty(\R^3)\cap \mathbb L^p_{\mathcal X_0}(\R^3))}^2.
        \end{multline}
        \end{enumerate}
    \end{enumerate}
    The constant appearing in the above estimates depends on both the lower and upper bounds of the density.
\end{coro} 
Indeed, all the estimates in \cref{coro1} are based in the following expression of the velocity gradient:
\begin{align}
    \nabla v &= \nabla \mathcal{P} v+ \nabla \mathcal{Q} v\nonumber\\
             &=-\dfrac{1}{\mu} (-\Delta)^{-1}\nabla \mathcal{P}(\rho_0\dpt v)-\dfrac{1}{\nu}(-\Delta)^{-1}\nabla\mathcal{Q}(\rho_0\dpt v)\nonumber\\
             &\quad +\dfrac{1}{\mu} (-\Delta)^{-1}\nabla \mathcal{P}\dvg f+\dfrac{1}{\nu}(-\Delta)^{-1}\nabla \mathcal{Q}\dvg f.\label{ep4.19}
\end{align}
The first two terms associated with $\dpt v$ exhibit regularity due to the regularity of $\dpt v$. In particular, their 
$L^r((0,T),L^p(\R^d))$ norm estimates can be obtained by interpolating  the estimate \eqref{ep4.12}. The $L^r((0,T), L^p(\R^d))$ norm estimate for the last two terms in the expression of the velocity gradient \eqref{ep4.19} follows from the continuity of Riesz operators 
 on $L^p(\R^d)$ for all $1<p<\infty$. These computations lead to \eqref{ep4.20}.
 
The estimate 
\eqref{ep4.22} is obtained  similarly: the $L^r((0,T), L^\infty(\R^d))$ norm of the terms associated with $\dpt v$ can be estimated
by interpolating the estimate \eqref{ep4.12}, while  the norm of the remaining terms is obtained  using \cref{Linftyestimate} since Riesz operators fail to be continuous on $L^\infty(\R^d)$.

To derive the estimate \eqref{ep4.21}, we  take time derivative of  \eqref{ep4.19} and apply the continuity of Riesz transforms  on $L^p(\R^d)$, $1<p<\infty$,
to obtain norms for the terms associated with $\dpt f$. The norm of the first two terms, associated with $\partial_{tt} v$,  can be obtained by interpolating estimate \eqref{ep4.13}.

For the last estimates \eqref{ep4.14}
and \eqref{ep4.18}, we take the derivative along $\mathcal X_0$ in \eqref{ep4.19} and we obtain:
\begin{align*}
    \partial_{\mathcal X_0} \nabla v&=-\dfrac{1}{\mu}\mathcal X_0\cdot \nabla(-\Delta)^{-1}\nabla\mathcal{P}(\rho_0\dpt v)-\dfrac{1}{\nu}\mathcal X_0\cdot\nabla (-\Delta)^{-1}\nabla\mathcal{Q}(\rho_0\dpt v)\\
                           &\quad + \dfrac{1}{\mu}\partial_{\mathcal X_0} (-\Delta)^{-1}\nabla\mathcal{P}\dvg f+\dfrac{1}{\nu}\partial_{\mathcal X_0}(-\Delta)^{-1}\nabla \mathcal{Q}\dvg f.
\end{align*}
Once again, the norms of the first two terms are obtained using H\"older's inequality and by interpolating estimates \eqref{ep4.12} and \eqref{ep4.13}. For the remaining terms, we use  Lemma A.1 of \cite{danchin2020well}. This completes this step of the study of the 
linear  system \eqref{ep4.11}.

\noindent\textbf{Step 4. Final conclusion.}
 Once we conclude the study of the linear system associated with \eqref{ep4.10}, the next step is to define a map that is contracting for some small time $T>0$, such that it admits a unique fixed point, which serves as a solution to \eqref{ep4.7} after reverting to Eulerian coordinates. With \cref{prop1} and \cref{coro1} in mind, we can verify that the unique solution of the full nonlinear system can be constructed in
\begin{gather*}
    \widetilde Z_T:=\left\{v\in Z_T \colon \nabla v\in  L^2((0,T), L^\infty(\R^d)\cap \mathbb L^p_{\mathcal X_0}(\R^d))\right\}
\end{gather*}
for  $2<p<\infty$ if $d=2$ and  $3<p\leqslant 6$ if $d=3$ by following the steps outlined in \cite[Section 4]{zodji2023well}. The
only argument we need to specify is the following lemmas.
\begin{lemma} 
    Let $v$ be a vector field verifying $\nabla v\in L^1((0,t), L^\infty(\R^d))$ and $\partial_{\mathcal X_0} \nabla v\in L^1((0,t), L^p(\R^d))$ for some $t>0$. Assuming that 
    \begin{gather}\label{ep4.23}
    V:=\int_0^t\left[\norm{\nabla v}_{L^\infty(\R^d)}+\norm{\partial_{\mathcal X_0}\nabla v}_{L^p(\R^d)}\right]<1,
    \end{gather}
    then, there exists a constant $K=K(V)$ such that the following estimate holds true:
    \begin{gather*}
        \norm{\partial_{\mathcal X_0}\adj (D \mathscr X_v(t)),\; \partial_{\mathcal X_0} A_{v}(t),\; \partial_{\mathcal X_0}J^{\pm 1}_{v}(t)}_{L^p(\R^d)}\leqslant K\norm{\partial_{\mathcal X_0}D v}_{L^1((0,t), L^p(\R^d))}.
    \end{gather*}
    Moreover, we have for all $Dw\in \mathbb L^p_{\mathcal X_0}(\R^d)$,
    \begin{multline*}
            \norm{\partial_{\mathcal X_0}(\adj(D \mathscr{X}_v(t))D_{A_v(t)}w)- \partial_{\mathcal X_0}D w }_{L^p(\R^d)}+\norm{\partial_{\mathcal X_0}(\adj(D \mathscr{X}_v(t))\dvg_{A_v(t)}w)- \partial_{\mathcal X_0}\dvg w }_{L^p(\R^d)}\\
            \leqslant K\left(\norm{D w}_{L^\infty(\R^d)}+\norm{\partial_{\mathcal X_0}D w}_{L^p(\R^d)}\right)\int_0^t\left(\norm{D v}_{L^\infty(\R^d)}+\norm{\partial_{\mathcal X_0} D v}_{L^p(\R^d)}\right).
    \end{multline*}
\end{lemma}

\begin{lemma}
    Let $v_1$ and $v_2$ two vector fields verifying \eqref{ep4.23}: $V_1,V_2<1$, and let $\delta v:= v_2-v_1$. Then, there exists a constant $K=K(V_1,V_2)$ such that the following estimate holds true:
    \begin{multline*}
        \norm{\bigl(\partial_{\mathcal X_0}A_{v_2}(t)-\partial_{\mathcal X_0}A_{v_1}(t),\; \partial_{\mathcal X_0}\adj (D \mathscr{X}_{v_2}(t))-\partial_{\mathcal X_0}\adj (D \mathscr{X}_{v_1}(t)),\; \partial_{\mathcal X_0}J_{v_2}^{\pm1}(t)-\partial_{\mathcal X_0}J_{v_1}^{\pm1}(t)\bigr)}_{L^p(\R^d)}\\
        \leqslant K\int_0^t\left(\norm{D \delta v}_{L^\infty(\R^d)}+\norm{\partial_{\mathcal X_0} D \delta v}_{L^p(\R^d)}\right).
    \end{multline*}
\end{lemma}

The particular case of $3<p\leqslant 6$ for $d=3$ is sufficient for constructing blocks for the global solution of \cref{th1}. For $6<p<\infty$ in three dimension, the fixed point theorem may be applied in a closed subset of the following space:
\begin{align*}
    \widetilde Z_T:=\left\{v\in Z_T \colon \sigma^{\tfrac{3}{4}-\tfrac{1}{r}}\nabla \dpt v\in L^r((0,T), L^3(\R^d));\;\sigma^{\tfrac{3}{4}-\tfrac{1}{r}-\tfrac{3}{2p}}\nabla v\in  L^r((0,T), L^\infty(\R^d)\cap \mathbb L^p_{\mathcal X_0}(\R^d))\right\}.
\end{align*}
 This ends the sketchy proof of \cref{intro:local}.

{\small 
\bibliographystyle{acm}
\bibliography{Biblio.bib}

\begin{thebibliography}{10}

\bibitem{bahouri2011fourier}
{\sc Bahouri, H., Chemin, J., and Danchin, R.}
\newblock {\em Fourier Analysis and Nonlinear Partial Differential Equations}.
\newblock Grundlehren der mathematischen Wissenschaften. Springer Berlin
  Heidelberg, 2011.

\bibitem{bertozzi1993}
{\sc Bertozzi, A.~L., and Constantin, P.}
\newblock {Global regularity for vortex patches}.
\newblock {\em Communications in Mathematical Physics 152}, 1 (1993), 19 -- 28.

\bibitem{bresch2021extension}
{\sc Bresch, D., and Burtea, C.}
\newblock {Extension of the Hoff solutions framework to cover compressible
  {N}avier-{S}tokes equations with possible anisotropic viscous tensor}.
\newblock {\em {Indiana University Mathematics Journal} 72}, 5 (2023),
  2145--2189.

\bibitem{charve2010global}
{\sc Charve, F., and Danchin, R.}
\newblock A global existence result for the compressible {N}avier-{S}tokes
  equations in the critical ${L}^p$ framework.
\newblock {\em Archive for Rational Mechanics and Analysis 198}, 1 (2010),
  233--271.

\bibitem{chemin1991mouvement}
{\sc Chemin, J.~Y.}
\newblock Sur le mouvement des particules d'un fluide parfait incompressible
  bidimensionnel.
\newblock {\em Inventiones mathematicae 103\/} (1991), 599--629.

\bibitem{chemin1993persistance}
{\sc Chemin, J.-Y.}
\newblock Persistance de structures g{\'e}om{\'e}triques dans les fluides
  incompressibles bidimensionnels.
\newblock In {\em Annales scientifiques de l'Ecole normale sup{\'e}rieure\/}
  (1993), vol.~26, pp.~517--542.

\bibitem{chen2010global}
{\sc Chen, Q., Miao, C., and Zhang, Z.}
\newblock Global well-posedness for compressible {N}avier-{S}tokes equations
  with highly oscillating initial velocity.
\newblock {\em Communications on Pure and Applied Mathematics 63}, 9 (2010),
  1173--1224.

\bibitem{coifman1993compensated}
{\sc Coifman, R., Lions, P.-L., Meyer, Y., and Semmes, S.}
\newblock Compensated compactness and {H}ardy spaces.
\newblock {\em Journal de math{\'e}matiques pures et appliqu{\'e}es 72}, 3
  (1993), 247--286.

\bibitem{Raphael_Danchin_Persistance}
{\sc Danchin, R.}
\newblock Persistance de structures g\'eom\'etriques et limite non visqueuse
  pour les fluides incompressibles en dimension quelconque.
\newblock {\em Bulletin de la Soci\'et\'e Math\'ematique de France 127}, 2
  (1999), 179--227.

\bibitem{danchin2020well}
{\sc Danchin, R., Fanelli, F., and Paicu, M.}
\newblock A well-posedness result for viscous compressible fluids with only
  bounded density.
\newblock {\em Analysis \& PDE 13}, 1 (2020), 275--316.
\newblock \href{https://arxiv.org/pdf/1804.09503}{available online}.

\bibitem{danchin2012lagrangian}
{\sc Danchin, R., and Mucha, P.~B.}
\newblock A lagrangian approach for the incompressible {N}avier-{S}tokes
  equations with variable density.
\newblock {\em Communications on Pure and Applied Mathematics 65}, 10 (2012),
  1458--1480.

\bibitem{danchin2013incompressible}
{\sc Danchin, R., and Mucha, P.~B.}
\newblock Incompressible flows with piecewise constant density.
\newblock {\em Archive for Rational Mechanics and Analysis 207\/} (2013),
  991--1023.

\bibitem{danchin2016compressible}
{\sc Danchin, R., and Mucha, P.~B.}
\newblock Compressible {N}avier-{S}tokes system : large solutions and
  incompressible limit, 2016.

\bibitem{Danchin_2019}
{\sc Danchin, R., and Mucha, P.~B.}
\newblock From compressible to incompressible inhomogeneous flows in the case
  of large data.
\newblock {\em Tunisian Journal of Mathematics 1}, 1 (Jan. 2019), 127–149.

\bibitem{danchin2019incompressible}
{\sc Danchin, R., and Mucha, P.~B.}
\newblock The incompressible {N}avier-{S}tokes equations in vacuum.
\newblock {\em Communications on Pure and Applied Mathematics 72}, 7 (2019),
  1351--1385.

\bibitem{danchin2023compressible}
{\sc Danchin, R., and Mucha, P.~B.}
\newblock Compressible {N}avier-{S}tokes equations with ripped density.
\newblock {\em Communications on Pure and Applied Mathematics 76}, 11 (2023),
  3437--3492.

\bibitem{danchin2023exponential}
{\sc Danchin, R., and Wang, S.}
\newblock Exponential decay for inhomogeneous viscous flows on the torus.
\newblock {\em Zeitschrift f{\"u}r angewandte Mathematik und Physik 75}, 2
  (2024), 62.
\newblock \href{https://arxiv.org/pdf/2311.01072}{available online}.

\bibitem{danchin2017persistence}
{\sc Danchin, R., and Zhang, X.}
\newblock On the persistence of h{\"o}lder regular patches of density for the
  inhomogeneous {N}avier-{S}tokes equations.
\newblock {\em Journal de l’{\'E}cole polytechnique-Math{\'e}matiques 4\/}
  (2017), 781--811.

\bibitem{denisova2001evolution}
{\sc Denisova, I.~V.}
\newblock Evolution of a closed interface between two liquids of different
  types.
\newblock In {\em European Congress of Mathematics: Barcelona, July 10--14,
  2000 Volume II\/} (2001), Springer, pp.~263--272.

\bibitem{Denisova2008}
{\sc Denisova, I.~V.}
\newblock Global solvability of the problem on the motion of two fluids without
  surface tension.
\newblock {\em Journal of Mathematical Sciences 152}, 5 (Aug. 2008), 625--637.

\bibitem{desjardins1997regularity}
{\sc Desjardins, B.}
\newblock Regularity of weak solutions of the compressible isentropic
  {N}avier-{S}tokes equations.
\newblock {\em Communications in Partial Differential Equations 22}, 5-6
  (1997), 977--1008.

\bibitem{fanelli2012}
{\sc Fanelli, F.}
\newblock Conservation of geometric structures for non-homogeneous inviscid
  incompressible fluids.
\newblock {\em Communications in Partial Differential Equations 37}, 9 (2012),
  1553--1595.

\bibitem{feireisl2001existence}
{\sc Feireisl, E., Novotn{\`y}, A., and Petzeltov{\'a}, H.}
\newblock On the existence of globally defined weak solutions to the
  {N}avier—{S}tokes equations.
\newblock {\em Journal of Mathematical Fluid Mechanics 3}, 4 (2001), 358--392.

\bibitem{GSR}
{\sc Gamblin, P., and Saint~Raymond, X.}
\newblock On three-dimensional vortex patches.
\newblock {\em Bull. Soc. Math. France 123}, 3 (1995), 375--424.

\bibitem{gancedo2023global}
{\sc Gancedo, F., and Garc{\'\i}a-Ju{\'a}rez, E.}
\newblock Global regularity of 2d {N}avier-{S}tokes free boundary with small
  viscosity contrast.
\newblock {\em Annales de l'Institut Henri Poincar{\'e} C\/} (2023).

\bibitem{haspot2011existence}
{\sc Haspot, B.}
\newblock Existence of global strong solutions in critical spaces for
  barotropic viscous fluids.
\newblock {\em Archive for Rational Mechanics and Analysis 202\/} (2011),
  427--460.

\bibitem{hoff1995global}
{\sc Hoff, D.}
\newblock Global solutions of the navier-stokes equations for multidimensional
  compressible flow with discontinuous initial data.
\newblock {\em Journal of Differential Equations 120}, 1 (1995), 215--254.
\newblock
  \href{https://www.sciencedirect.com/science/article/pii/S0022039685711114}{available
  online}.

\bibitem{hoff1995strongpoly}
{\sc Hoff, D.}
\newblock Strong convergence to global solutions for multidimensional flows of
  compressible, viscous fluids with polytropic equations of state and
  discontinuous initial data.
\newblock {\em Archive for rational mechanics and analysis 132\/} (1995),
  1--14.

\bibitem{hoff2002dynamics}
{\sc Hoff, D.}
\newblock Dynamics of singularity surfaces for compressible, viscous flows in
  two space dimensions.
\newblock {\em Communications on Pure and Applied Mathematics: A Journal Issued
  by the Courant Institute of Mathematical Sciences 55}, 11 (2002), 1365--1407.

\bibitem{hoff2008lagrangean}
{\sc Hoff, D., and Santos, M.~M.}
\newblock Lagrangean structure and propagation of singularities in
  multidimensional compressible flow.
\newblock {\em Archive for rational mechanics and analysis 188}, 3 (2008),
  509--543.

\bibitem{hoff1985solutions}
{\sc Hoff, D., and Smoller, J.}
\newblock Solutions in the large for certain nonlinear parabolic systems.
\newblock In {\em Annales de l'Institut Henri Poincar{\'e} C, Analyse non
  lin{\'e}aire\/} (1985), vol.~2, Elsevier, pp.~213--235.

\bibitem{hu2020optimal}
{\sc Hu, X., and Wu, G.}
\newblock Optimal decay rates of isentropic compressible {N}avier-{S}tokes
  equations with discontinuous initial data.
\newblock {\em Journal of Differential Equations 269}, 10 (2020), 8132--8172.

\bibitem{itaya1970existence}
{\sc Itaya, N.}
\newblock The existence and uniqueness of the solution of the equations
  describing compressible viscous fluid flow.
\newblock {\em Proceedings of the Japan Academy 46}, 4 (1970), 379--382.

\bibitem{itaya1971cauchy}
{\sc Itaya, N.}
\newblock On the cauchy problem for the system of fundamental equations
  describing the movement of compressible viscous fluid.
\newblock In {\em Kodai Mathematical Seminar Reports\/} (1971), vol.~23,
  Department of Mathematics, Tokyo Institute of Technology, pp.~60--120.

\bibitem{leray1934mouvement}
{\sc Leray, J.}
\newblock Sur le mouvement d'un liquide visqueux emplissant l'espace.
\newblock {\em Acta mathematica 63\/} (1934), 193--248.

\bibitem{liao2019global}
{\sc Liao, X., and Liu, Y.}
\newblock Global regularity of three-dimensional density patches for
  inhomogeneous incompressible viscous flow.
\newblock {\em Science China Mathematics 62\/} (2019), 1749--1764.

\bibitem{liao2016global}
{\sc Liao, X., and Zhang, P.}
\newblock On the global regularity of the two-dimensional density patch for
  inhomogeneous incompressible viscous flow.
\newblock {\em Archive for Rational Mechanics and Analysis 220\/} (2016),
  937--981.

\bibitem{liao2019globallow}
{\sc Liao, X., and Zhang, P.}
\newblock Global regularity of 2d density patches for viscous inhomogeneous
  incompressible flow with general density: low regularity case.
\newblock {\em Communications on Pure and Applied Mathematics 72}, 4 (2019),
  835--884.

\bibitem{lions1996mathematicalv1}
{\sc Lions, P.}
\newblock {\em Mathematical Topics in Fluid Mechanics: Volume 1: Incompressible
  Models}.
\newblock Mathematical Topics in Fluid Mechanics. Oxford University Press,
  Incorporated, 1996.

\bibitem{lions1996mathematical}
{\sc Lions, P.}
\newblock {\em Mathematical Topics in Fluid Mechanics: Volume 2: Compressible
  Models}.
\newblock Mathematical Topics in Fluid Mechanics. Clarendon Press, 1996.

\bibitem{matsumura1980initial}
{\sc Matsumura, A., and Nishida, T.}
\newblock The initial value problem for the equations of motion of viscous and
  heat-conductive gases.
\newblock {\em Journal of Mathematics of Kyoto University 20}, 1 (1980),
  67--104.

\bibitem{nash1962probleme}
{\sc Nash, J.}
\newblock Le probl{\`e}me de cauchy pour les {\'e}quations diff{\'e}rentielles
  d'un fluide g{\'e}n{\'e}ral.
\newblock {\em Bulletin de la Soci{\'e}t{\'e} Math{\'e}matique de France 90\/}
  (1962), 487--497.

\bibitem{paicu2020striated}
{\sc Paicu, M., and Zhang, P.}
\newblock Striated regularity of 2-d inhomogeneous incompressible
  {N}avier-{S}tokes system with variable viscosity.
\newblock {\em Communications in Mathematical Physics 376}, 1 (2020), 385--439.

\bibitem{prange2023inhomogeneous}
{\sc Prange, C., and Tan, J.}
\newblock Inhomogeneous incompressible viscous flows in vacuum: the whole-space
  case.
\newblock {\em arXiv preprint arXiv:2310.09288\/} (2023).

\bibitem{serre1991variations}
{\sc Serre, D.}
\newblock Variations de grande amplitude pour la densite d'un fluide visqueux
  compressible.
\newblock {\em Physica D: Nonlinear Phenomena 48}, 1 (1991), 113--128.

\bibitem{shibata2023global}
{\sc Shibata, Y., and Zhang, X.}
\newblock Global wellposedness of the 3d compressible {N}avier-{S}tokes
  equations with free surface in the maximal regularity class.
\newblock {\em Nonlinearity 36}, 7 (2023), 3710.

\bibitem{solonnikov1980solvability}
{\sc Solonnikov, V.}
\newblock Solvability of the initial-boundary-value problem for the equations
  of motion of a viscous compressible fluid.
\newblock {\em Journal of Soviet Mathematics 14\/} (1980), 1120--1133.

\bibitem{tani1976existence}
{\sc Tani, A.}
\newblock The existence and uniqueness of the solution of equations describing
  compressible viscous fluid flow in a domain.
\newblock {\em Proceedings of the Japan Academy 52}, 7 (1976), 334--337.

\bibitem{zodji2023discontinuous}
{\sc Zodji, S.~M.}
\newblock Discontinuous solutions for the {N}avier-{S}tokes equations with
  density-dependent viscosity.
\newblock {\em arXiv preprint arXiv:2312.07578\/} (2023).

\bibitem{zodji2023well}
{\sc Zodji, S.~M.}
\newblock A well-posedness result for the compressible two-fluid model with
  density-dependent viscosity.
\newblock {\em arXiv preprint arXiv:2310.12525\/} (2023).

\end{thebibliography}
}
\end{document}